\documentclass[12pt,notitlepage]{amsart}
\usepackage{amsmath} 
\usepackage{amscd}
\usepackage{amssymb}
\usepackage{enumerate}
\usepackage{indentfirst}
\usepackage[dvips]{graphicx}

\advance\voffset by -0.6 truein
\advance\hoffset by -0.8 truein
 \font\tenbi=cmmi14
 \font\sevenbi=cmmi10 \font\fivebi=cmmi7
 \newfam\bifam \def\bi{\fam\bifam} \textfont\bifam=\tenbi
 \scriptfont\bifam=\sevenbi \scriptscriptfont\bifam=\fivebi
 \mathchardef\variablemega="7121 \def\bigomega{{\bi\variablemega}}

\renewcommand{\l}{\ell}
\newcommand{\Lcal}{\mathcal{L}}
\newcommand{\Ocal}{\mathcal{O}}
\newcommand{\Mcal}{\mathcal{M}}
\newcommand{\Ncal}{\mathcal{N}}
\newcommand{\wh}{\widehat}
\newcommand{\wt}{\widetilde}
\newcommand{\ol}{\overline}
\renewcommand{\:}{\colon}

\newcommand{\Z}{\mathbb{Z}}
\renewcommand{\P}{\mathbb{P}}
\newcommand{\Q}{\mathbb{Q}}
\newcommand{\R}{\mathbb{R}}
\newcommand{\w}{\bigomega}
\newcommand{\V}{\mathbb V}
\newcommand{\G}{\mathbb G}
\renewcommand{\O}{\mathbb O}

\DeclareMathOperator{\Grass}{Grass}
\DeclareMathOperator{\Image}{Im}

\newtheorem{theorem}{Theorem}[section]
\newtheorem{lemma}[theorem]{Lemma}
\newtheorem{corollary}[theorem]{Corollary}
\newtheorem{proposition}[theorem]{Proposition}

\theoremstyle{definition} 
\newtheorem{definition}[theorem]{Definition}
\newtheorem{remark}[theorem]{Remark}
\newtheorem{example}[theorem]{Example}

\newtheorem{subsct}[theorem]{}
\theoremstyle{plain}

\font\smallrm=cmr8
\font\midrm=cmr10
\font\smallsc=cmcsc10
\font\smallsl=cmsl10

\begin{document}
\author[{\smallrm Esteves and Medeiros}]
{Eduardo Esteves \and Nivaldo Medeiros}
\thanks{First author supported by PRONEX, Conv\^enio 41/96/0883/00, CNPq, 
Proc. 300004/95-8, and FAPERJ, Proc. E-26/170.418/2000-APQ1}
\title[{\smallrm Limit canonical systems on curves with two components}]
{Limit canonical systems \\ on curves with two components}
\begin{abstract}
In the 80's D. Eisenbud and 
J. Harris considered the following problem: 
``What are the limits of Weierstrass points in families 
of curves degenerating to stable curves?'' But for the case of stable 
curves of compact type, treated by them, this problem remained wide 
open since then. In the present article, we propose a concrete approach to 
this problem, and give a quite explicit 
solution for stable curves with 
just two irreducible components meeting at points in general position.
\end{abstract}
\maketitle

\tableofcontents

\section{Introduction} 
\begin{subsct}\setcounter{equation}{0}
\emph{Goal.}
In \cite{EH2} D. Eisenbud and 
J. Harris asked: ``What are the limits of Weierstrass points in families 
of curves degenerating to stable curves \emph{not} of compact type?'' 
As announced in \cite{EM}, 
in the present article we answer their question for stable curves with 
just two irreducible components meeting at points in general position.

For a glimpse of our results, consider the following particular case: Let 
$C$ be a nodal curve with two elliptic components $E_1$ and $E_2$ meeting 
at $\delta$ points $p_1,\dots,p_\delta$ in general position. View 
$C$ in $\mathbb P^{3\delta-1}$ under the bicanonical map,
and let $H\subseteq\mathbb P^{3\delta-1}$ be the projective subspace 
spanned by $p_1,\dots,p_\delta$. For each 
hyperplane $L\subseteq H$ let $R_1(L)$ and $R_2(L)$ be the ramification 
divisors of the projections of $E_1$ and $E_2$ with center $L$, and put
\[
W(L):=R_1(L)+R_2(L)+(\delta+1)(\delta-2)(p_1+\cdots+p_\delta).
\]
It follows 
from our Theorem~10.7 
(see Example~11.3 as well) 
that each $W(L)$ 
is limit of Weierstrass divisors in 
one--parameter 
families of smooth 
curves degenerating to $C$, and there are no other such limits.  
\end{subsct}

\begin{subsct}\setcounter{equation}{0}
\emph{Some history.} 
Limits of ramification points and linear systems were 
studied 
by Eisenbud and Harris in the 80's, when they developed the theory of 
limit linear series for curves $C$ of compact 
type; see~\cite{EH1}. 
(As 
observed in \cite[p. 256]{HM}, the theory works more generally when $C$ 
is \emph{treelike}, that is, when the normalization of $C$ at all its 
irreducible nodes is of compact type.) 
Many important applications of their theory were found; a survey 
is given in \cite{EH1}. 

Since many applications were obtained from a theory applicable only to 
curves of a special type, it's natural to expect 
more 
applications 
from a more general theory. In fact, Eisenbud and Harris wrote in 
\cite[p. 220]{EH3} that 
``\ldots there is probably a small gold mine awaiting a general insight.'' 
Despite the potential for applications and the vast interest the topic 
generated in the 80's, the theory of limit linear 
series could not be extended to stable curves of more general type. 
There was unpublished work of Z. Ran on degenerations of linear series 
\cite{Ran}, but the crucial relationship with degenerations of ramification 
points was yet to be established.

It was only recently that this relationship was established by the first 
author in 
\cite{Esteves2}, 
building up on the work of Ran 
(or rather rediscovering it) and using the substitutes for the sheaves of 
(relative) principal parts discovered in 
\cite{Esteves1}. It is apparently possible to generalize the 
theory of limit linear series of \cite{EH1} starting from results in 
\cite{Esteves2}.

Though a more general theory of limit linear series is yet to be developed, 
it 
is already 
possible to apply the results in 
\cite{Esteves2} 
to answer the 
question raised by Eisenbud and 
Harris; see~1.1. For starters 
we describe in the present article a complete answer for 
stable curves with just two irreducible components
meeting at points in general position.
This is a substantial 
breakthrough, as the question was 
not completely answered even for
the simplest 
case of 
stable 
curves with two elliptic components meeting at two points
in general position.

In fact, it was only very recently that partial progress was made by 
M. Coppens and 
L.~Gatto
towards answering the question. In 
\cite{CG} Coppens and Gatto 
consider a stable curve $C$ with just two irreducible 
components meeting at $\delta$ points in general position, and show that 
all the points in each collection of $\delta-1$ smooth points of $C$ are 
limits of Weierstrass points in a suitable family 
of smooth curves degenerating to $C$. (Their result is 
partially
recovered by 
our Proposition~11.1.) Besides this general result, 
they study the case where one of the components of $C$ is elliptic and 
obtain partial results. In contrast with 
the methods used in the present article, 
they use 
admissible 
covers.

Besides 
\cite{Esteves2}, 
the recent work by L. Main\`o on enriched structures 
in her Ph.D. thesis~\cite{Maino} was important for us. Actually, as we 
deal here only with limits of canonical systems, we used only a fraction of 
what is available in \cite{Maino}.

There are 
further 
recent works on (generalizations of) 
limit linear series worth noticing, 
even though they are not completely related to what we do here. First, 
there is a draft by 
R.~Pandharipande 
\cite{P} 
who takes a Geometric Invariant Theory 
approach to proving the existence of what he 
calls generalized linear series for certain degenerations. Second, there is 
A. Bruno's thesis \cite{B1}, where he studies 
degenerations of linear systems to a stable curve with two components. 
It's particularly interesting the conditions he gives 
for such a curve to be a limit of smooth plane quintics. Finally, 
there is \cite{B2}, where applications to (variants of) a 
problem of Severi's are given.

We thank A. Bruno, M. Coppens, D. Eisenbud, J. Harris, S. Kleiman and 
specially L.~Gatto and L. Main\`o for helpful conversations. 
\end{subsct}

\begin{subsct}\setcounter{equation}{0}
\emph{The general approach.} 
Rather than producing 
directly a space parameterizing limits of Weierstrass divisors on a 
stable curve, our approach consists of producing a variety $\V$ of 
limits of canonical systems and a formula that 
gives for each point 
of $\V$ the corresponding limit of Weierstrass divisors. 
This method is applicable to any stable curve, as we shall 
describe below.

Let $C$ be a connected, local complete intersection, projective curve of 
arithmetic genus $g>0$ defined over an 
algebraically closed field $k$ of characteristic 0. 
Let $B:=\text{Spec}(k[[t]])$, and 
denote by $o$ its special point and $\eta$ its generic point. 
A \emph{smoothing} of $C$ is a 
projective and flat map $\pi\:S\to B$ such that $S_\eta$ is smooth and 
$S_o=C$. For each smoothing $\pi\:S\to B$ of $C$ let 
$W_\pi:=\ol W_\eta\cap C$, where $W_\eta\subseteq S$ is the 
Weierstrass 
subscheme
of $S_\eta$. We call the associated cycle $[W_\pi]$ 
the \emph{limit Weierstrass divisor of $\pi$}.

Let $C_1,\dots,C_n$ denote the irreducible components of $C$. 
A smoothing $\pi\:S\to B$ of $C$ is called \emph{regular} if $S$ is 
regular everywhere but possibly at the singularities of $C$ that lie on 
just one component $C_i$. 
If $\pi\:S\to B$ is a regular smoothing of $C$, then $C_1,\dots,C_n$ are 
Cartier divisors on $S$.

Assume from now on that the irreducible components of $C$ meet 
at nodes, that is, ordinary double points. Let $\Delta$ 
denote the set of 
points of intersection between distinct components of $C$. 
For each $\mu$ in $\Z_\Delta^+:=\textstyle\prod_{p\in\Delta}\Z^+$ let 
$\wt C$ be the curve obtained from $C$ by splitting the 
branches of $C$ at each $p\in\Delta$ and 
connecting them by a chain of $\mu_p-1$ rational, smooth curves
(see Figure~1).
\begin{figure}[htb]
\begin{center}
\includegraphics[angle=0]{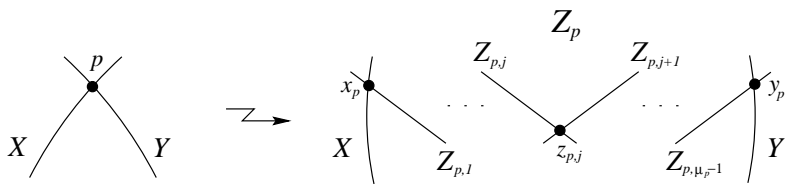}
\caption{Semi--stable reduction.}
\end{center}
\end{figure}
We say that $\wt C$ is the \emph{$\mu$--semi--stable model of $C$}. 
The smoothings of $C$ correspond to the regular smoothings of 
the semi--stable models of $C$. This correspondence is called 
\emph{semi--stable reduction}; see~2.7. We shall view the irreducible 
components of $C$ as irreducible components of $\wt C$.

For each $i=1,\dots,n$ let $d_i\in\Z$ and put 
$\mathbb L_i:=\w_i((1+d_i)\textstyle\sum_{p\in\Delta_i}p)$, where 
$\w_i$ is the dualizing sheaf of $C_i$ and $\Delta_i:=\Delta\cap C_i$.

Let $\mu\in\Z_\Delta^+$ and consider a regular smoothing 
$\wt\pi\:\wt S\to B$ of 
the $\mu$--semi--stable model $\wt C$ of $C$. Let $\w_{\wt\pi}$ 
be the (relative) dualizing sheaf of $\wt\pi$. 
By Theorem~2.2, for each $i=1,\dots,n$ there is a unique 
Cartier divisor 
$D_i$
on $\wt S$ that is 
supported in the union of all irreducible components of $\wt C$ but $C_i$ 
and such that the following two properties hold for 
$\Lcal_i:=\w_{\wt\pi}(D_i)$.
\begin{enumerate}
\item The natural map 
$\varrho_i\:H^0(\Lcal_i)(o)\to H^0(\Lcal_i|_{C_i})$
is injective.
\item The natural map 
$\varrho_{i,j}\:H^0(\Lcal_i)(o)\to H^0(\Lcal_i|_{C_j})$ 
is not zero if $j\neq i$.
\end{enumerate}
Let $V_{\wt\pi,i}:=\Image(\varrho_i)$. We call 
$\Lcal_i$
the 
\emph{canonical sheaf of $\wt\pi$ with focus on $C_i$} 
and 
$(V_{\wt\pi,i},\Lcal_i|_{C_i})$
the \emph{limit canonical aspect of $\wt\pi$ with focus on $C_i$}. 
For each 
$p\in\Delta_i$ 
denote by $\ell_i(p)$ the coefficient 
in $D_i$ of 
the irreducible component of 
$\wt C$ other than $C_i$ containing $p$. We call $\ell_i(p)$ the 
\emph{correction number for 
$\Lcal_i$ 
at $p$}. Let 
$d_{\wt\pi,i}:=\max\{\ell_i(p)\,|\, p\in\Delta_i\}$. If 
$d_i\geq d_{\wt\pi,i}$ then $\Lcal_i|_{C_i}\subseteq\mathbb L_i$, and thus
$V_{\wt\pi,i}$ can be viewed as a 
$g$--dimensional subspace of 
$H^0(\mathbb L_i)$.

Assume from now on that 
$d_i\geq d_{\wt\pi,i}$ for every regular smoothing $\wt\pi$ of every 
semi--stable model of $C$
and every $i=1,\dots,n$.
(This condition holds if 
there are no loosely connected rational tails in $C$; see Remark~2.5.) 
Let 
\[
\G:=\Grass_g(H^0(\mathbb L_1))\times\cdots\times
\Grass_g(H^0(\mathbb L_n)).
\]
Given a regular smoothing $\wt\pi$ of a semi--stable model 
of $C$ let 
\[
\nu_{\wt\pi}:=(V_{\wt\pi,1},\dots,V_{\wt\pi,n})\in\G.
\]
For each $\mu\in\Z_\Delta^+$ let 
\[
\V_{\mu}:=\{\nu_{\wt\pi}\,|\, \wt\pi\text{ is a regular smoothing of the 
$\mu$--semi--stable model of $C$}\}.
\]
Put
\begin{equation}
\V:=\bigcup_{\mu\in\Z_\Delta^+}\V_\mu\subseteq\G.
\end{equation}
We call $\V$ the \emph{variety of limit canonical systems of $C$}.

For each $\nu=(V_1,\dots,V_n)\in\G$ let
\[
W_\nu:=R_{\nu,1}+\cdots+R_{\nu,n}+\sum_{i<j}
\sum_{p\in C_i\cap C_j}g(g-1-d_i-d_j)p,
\]
where $R_{\nu,i}$ is the ramification 
divisor of $(V_i,\mathbb L_i)$ for each $i=1,\dots,n$. 
By Theorem~2.8 (see Remark~2.9 as well), if 
$\wt\pi$ is a regular smoothing of a semi--stable model of $C$ then 
$W_{\nu_{\wt\pi}}$
is the limit Weierstrass divisor of the corresponding 
smoothing of $C$.

So, in order to answer the question posed by Eisenbud and Harris (see~1.1), 
it's enough to understand $\V$ well. We propose the following problems:
\begin{enumerate}
\item Describe explicitly each $\V_\mu$ and the interrelations among 
the various $\V_\mu$.
\item Find properties of $\V$. For instance, is $\V$ projective?
Is $\V$ connected? What is the dimension of $\V$? When is $\V$ irreducible?
What is the number of irreducible components of $\V$?
\end{enumerate}
\end{subsct}

\begin{subsct}\setcounter{equation}{0}
\emph{Results.} 
Solving the problems raised at the 
end of 1.3 can be 
difficult without further hypotheses. In the present 
article we apply the approach of 1.3 for curves with 
just two irreducible components meeting at points in general position. 
We describe our results below.

Let $C$ be a connected, local complete intersection, 
projective curve of arithmetic genus $g>0$ defined over an 
algebraically closed field $k$ of characteristic 0. Assume that 
$C$ has only two irreducible components, denoted $X$ and $Y$, and that the 
points of $X\cap Y$ are nodes of $C$. Let 
$\Delta:=X\cap Y$ and $\delta:=|\Delta|$. Let $g_X$ and $g_Y$ be the 
arithmetic genera of $X$ and $Y$. Then $g=g_X+g_Y+\delta-1$. 
Assume that $\delta>1$ or 
$g_Xg_Y>0$. Let $\w_X$ and $\w_Y$ denote the 
dualizing sheaves of $X$ and $Y$.

Assume from now on 
that $\Delta$ lies in general position on $X$ and $Y$. More 
precisely, assume that each effective divisor $D$ on $X$ (resp. $Y$) 
supported in $\Delta$ and of degree at most $g_X$ (resp. $g_Y$) 
imposes independent conditions on $\w_X$ (resp. $\w_Y$).

Let $\V$ be the variety of limit canonical systems of $C$, 
as described in 1.3. Our Theorem~10.6 asserts that 
$\V$ is connected, projective and of pure dimension. 
Moreover, $\dim\V=\delta-1$ unless $g_X=g_Y=0$, 
in which case $\V$ is a point. Our Theorem~11.2 gives a lower bound for 
the number $N(\V)$ of irreducible components of $\V$. Using this lower 
bound we show that $\V$ is irreducible if and only if $g_X,g_Y\leq 1$. If 
$\delta=2$ then $N(\V)=g-\gcd(g_X+1,g_Y+1)$; see 11.5.

We solve Problem~1 of 1.3 as well. In fact, our solution to Problem~2 
follows from our quite explicit solution to Problem~1. We proceed as 
follows. To each $\mu$ in $\Z_\Delta^+:=\prod_{p\in\Delta}\Z^+$ 
we associate non--empty 
subsets $I_\mu,J_\mu\subseteq\Delta$ and $\delta$--tuples 
$\alpha_\mu,\beta_\mu\in\Z_\Delta^+$ such that $\alpha_{\mu,p}\leq g_Y$ 
and $\beta_{\mu,p}\leq g_X$ for each $p\in\Delta$; see~4.1. 
These associated data $\alpha_\mu,\beta_\mu,I_\mu,J_\mu$ 
are obtained numerically from $\mu$ alone, but they have geometric meaning 
as well. In fact, let $\wt C$ be the $\mu$--semi--stable model of $C$ and 
$\phi\:\wt C\to C$ the induced map. 
Let $\Lcal$ be the canonical sheaf with focus on $X$ (resp. $Y$) of a 
regular smoothing of $\wt C$.
Let $p\in\Delta$. Then
$p\in I_\mu$ (resp. $p\in J_\mu$) if and only if $\phi_*(\Lcal|_{\wt C})$ 
is invertible at $p$; see Remark~5.3. In addition,
$\alpha_{\mu,p}$ (resp. $\beta_{\mu,p}$) is the correction number 
for $\Lcal$ 
at the point of $\wt C$ 
lying on $X$ (resp. $Y$) and above $p$; 
see Theorem~5.2. 
So $\V_\mu\subseteq\G$, where
\[
\G:=\Grass_g(H^0(\w_X((1+g_Y)\textstyle\sum_{p\in\Delta}p)))\times
\Grass_g(H^0(\w_Y((1+g_X)\textstyle\sum_{p\in\Delta}p))).
\]

It's interesting to note that, since $I_\mu$ and $J_\mu$ are non--empty, 
the canonical sheaves with foci on $X$ and $Y$ of every regular smoothing 
of every semi--stable model of $C$ push forward to (relatively) 
simple sheaves on the corresponding smoothing of $C$. (Recall that a sheaf 
$L$ on $C$ is \emph{simple} if 
$\text{End}(L)=k$.) 
So we need only consider 
sheaves that are simple on $C$ to obtain a projective variety $\V$. This 
is the same situation found in \cite{Esteves3}, where a compactification of 
the (relative) Jacobian of a family of curves 
was constructed using only simple sheaves.

Theorem~10.1 asserts the importance of the data 
$\alpha_\mu,\beta_\mu,I_\mu,J_\mu$:
they determine $\V_\mu$. In fact, an 
explicit description of $\V_\mu$ from these data is given by 
Theorem~8.2. As a corollary we obtain Theorem~8.5, which asserts 
that each $\V_\mu$ is isomorphic to a torus and gives a formula for 
the dimension of $\V_\mu$.

Theorem~10.1 asserts as well that any two $\V_\mu$ and $\V_{\mu'}$ are 
either equal or disjoint. Now, there are 
only 
finitely many $\V_\mu$ as there 
are 
only 
finitely many associated data $\alpha_\mu,\beta_\mu,I_\mu,J_\mu$. 
So the covering (1.3.1) can actually be thought of as a finite 
stratification of $\V$. 

It remains to 
describe
the relationship between the several 
$\V_\mu$. First, Theorem~10.1 says that $\V_{r\mu}=\V_\mu$ 
for each positive integer $r$. This 
fact is not surprising as each regular 
smoothing of the $\mu$--semi--stable model of $C$ 
induces a regular smoothing of the $r\mu$--semi--stable model after 
base change of degree $r$ 
and resolution of singularities. At any rate, this fact allows us to 
define $\V_\mu$ for each $\mu$ in 
$\Q_\Delta^+:=\textstyle\prod_{p\in\Delta}\Q^+$ in the obvious way. Then we 
obtain a remarkable relation between the Zariski topology of $\V$ and the 
Euclidean topology of $\Q_\Delta^+$: Theorem~10.5 asserts that 
there is a neighborhood $U_\mu\subseteq\Q_\Delta^+$ of each 
$\mu\in\Q_\Delta^+$ such that
the closure $\ol\V_\mu\subseteq\G$ satisfies
\[
\ol\V_\mu=\bigcup_{\ol\mu\in U}\V_{\ol\mu}
\]
for every open neighborhood $U\subseteq U_\mu$ of $\mu$.
\end{subsct}

\begin{subsct}\setcounter{equation}{0}
\emph{An overview.} Let $C$ be a connected, local complete intersection, 
projective curve of arithmetic genus $g>0$ defined over an 
algebraically closed field of characteristic 0. Assume 
that the irreducible components of $C$ meet at nodes.

In Section~2 we review the 
fundamentals of the theory of limit linear systems and limit ramification 
points on $C$ developed in 
\cite{Esteves2}, 
as applied to limits of canonical 
systems and Weierstrass points. In contrast with 
\cite{Esteves2} 
we show 
here how to handle non--regular smoothings of $C$ 
using semi--stable reduction; see~2.7 and Theorem~2.8. 
Semi--stable reduction can be incorporated to the 
general theory of limit linear systems and limit ramification 
points, but we refrain from carrying this out here. 

In 
Section~3
we prove the purely numerical Lemma~3.1. 

In Section~4 we define our basic set--up, introducing the terminology that 
shall be used 
in the rest of the article. Except in Section~6, it is assumed 
from here on 
that $C$ has just two 
irreducible components $X$ and $Y$. Let $\Delta:=X\cap Y$. 
We state as well 
the general--position conditions, (4.3.1) and (4.3.2), that will be 
assumed later. Using 
Lemma~3.1,
we define for each 
$\mu$ in $\Z_\Delta^+:=\prod_{p\in\Delta}\Z^+$ the numerical data 
$\alpha_\mu,\beta_\mu,I_\mu,J_\mu$ which will a posteriori 
determine the stratum $\V_\mu$ of the variety of limit canonical 
systems $\V$. 

In Section~5 we determine for each regular smoothing of a semi--stable model 
$\wt C$ of $C$ the associated canonical sheaves and limit canonical aspects 
with foci on $X$ and $Y$. Our main result here is Theorem~5.2.

In Section~6 we review the part of the theory of enriched 
structures of Main\`o's 
which 
is necessary for us. 
Unfortunately, Main\`o's
results do not seem to meet directly our needs. So we
follow an approach that is slightly different from hers, and
give independent proofs.
For instance, for us an enriched structure
is a collection of invertible fractional ideals on $C$ satisfying certain 
properties, rather than just the collection of their associated (abstract) 
invertible sheaves; see~6.4 and Remark~6.7. The main result of 
this section, Theorem~6.5, gives a criterion for when a collection 
of fractional ideals on $C$ comes from an enriched structure. 
It would actually be more useful to have instead
a similar 
criterion valid for collections of invertible sheaves.

In Section~7 we begin our study of the regeneration process. More 
precisely, we fix $\mu\in\Z_\Delta^+$, and determine which 
invertible sheaves 
on the $\mu$--semi--stable model $\wt C$ of $C$ 
are restrictions of canonical sheaves of 
regular smoothings of $\wt C$ with foci on $X$ and/or $Y$. Theorem~6.5 
is our main tool.

In Section~8 we apply 
the theorem of 
Section 7 to give a characterization 
of the stratum $\V_\mu$ in terms of the numerical data 
$\alpha_\mu,\beta_\mu,I_\mu,J_\mu$. Our main results here are 
Theorems~8.2~and~8.5.

In Section~9 we study boundary points of certain tori orbits in 
products of Grassmannians. This section 
does not rely on
the rest of the paper and might be of independent interest. In fact, 
Lemma~9.1 is probably well known.
We give a proof of it as the proof 
of Lemma~9.2 is based on that proof.

In Section 10 we prove our most important results, 
Theorems~10.1,~10.5,~10.6~and~10.7. Here we tackle 
all strata $\V_\mu$ of $\V$ together and determine which lie in the 
closure in $\V$ of which. Lemma~9.2 is needed for this purpose.

In Section~11 we collect a few results: Proposition~11.1 
recovers in part the main result of \cite{CG} and Theorem~11.2
gives estimates on the number of irreducible
components of $\V$.
In addition, we describe
how to represent $\V$ graphically (see~11.4) and describe this 
representation for $|\Delta|=2$ and $|\Delta|=3$; see~11.5--6. 
\end{subsct}

\begin{subsct}\setcounter{equation}{0}
\emph{Notation.} 
If $F$ and $G$ are sets, let 
$G_F:=\prod_{p\in F}G$. If $G$ is a group then so is $G_F$; in addition, 
if $E\subseteq F$ 
we shall 
view $G_E$ as a subgroup of $G_F$ 
in the natural way.

Let $F$ be a finite set. Given 
$b,c\in\Z_F$ write $b\leq c$ if $b_p\leq c_p$ for every $p\in F$. 
For each $b\in\Z_F$ let $|b|:=\textstyle\sum_{p\in F}b_p$. 
View each $p\in F$ as an element of $\Z_F$, defined 
by $p_p:=1$ and $p_q:=0$ if $q\neq p$. In addition, view each 
$b\subseteq F$ as the sum $\sum_{p\in b}p$ inside $\Z_F$. 

Let $k$ be a field. Given $b\in\Z_F$ and $t\in k^*_F$, 
let $t^b:=\prod_{p\in F}t_p^{b_p}$. The group $k^*_F$ 
acts diagonally on $k_F$: given $t\in k^*_F$ and $v\in k_F$, 
define $t\cdot v\in k_F$ by $(t\cdot v)_p:=t_pv_p$ for each 
$p\in F$.

If $\mathcal F$ is a coherent sheaf and $D$ is a 
Cartier divisor on a scheme $X$, let $\Ocal_X(D)$ 
denote the invertible sheaf associated to $D$, and put 
$\mathcal F(D):=\mathcal F\otimes\Ocal_X(D)$. If $X$ is a scheme and $Y$ 
a 
closed
subscheme, let $[Y]$ denote the associated cycle on $X$. 

In the present article, a \emph{curve} is a connected, reduced, 
\emph{local complete intersection}, projective 
scheme of pure dimension 1 defined 
over an algebraically closed field of characteristic 0. 
A singularity of a curve is called 
\emph{irreducible} 
if it 
lies on 
just one irreducible component, and 
\emph{reducible} 
if not. It is called 
a \emph{node} if it is an ordinary double point. In the present article, 
a curve is called \emph{nodal} if all its \emph{reducible} singularities 
are nodes. In addition, a curve is called \emph{treelike} if 
its normalization at all irreducible singularities is of compact type. 
\end{subsct}

\section{Limit canonical systems and Weierstrass divisors}

\begin{subsct}\setcounter{equation}{0} 
\emph{Limit canonical systems.} 
Let $C$ be a 
curve and $C_1,\dots,C_n$ its irreducible components. 
A \emph{smoothing} of $C$ is a projective and flat map 
$\pi\:S\to B$, where $B$ is the spectrum of a discrete valuation ring, 
the special fiber of $\pi$ is $C$ and the generic fiber is 
smooth. If $\pi$ is a smoothing 
of $C$ let $\w_\pi$ denote the (relative) dualizing sheaf 
of $\pi$. Then $\w_\pi$ is an invertible extension 
to $S$ of the canonical sheaf on the generic fiber of $\pi$. 

A smoothing $\pi\:S\to B$ of $C$ is 
said to be \emph{regular} if $S$ is regular everywhere \emph{but} possibly 
at the 
irreducible 
singularities of $C$. 

Let $\pi\:S\to B$ be a regular smoothing of $C$ and $o$ the special 
point of $B$. Then $C_i$ is a Cartier divisor on $S$ for each $i=1,\dots,n$. 
(In fact, $C_i$ is principal 
away from the irreducible singularities of $C$ 
because $S$ is regular 
there, and 
at each irreducible singularity of $C$ in $C_i$ 
because $C_i$ is 
locally the pull--back of $o$.) 
In addition, $C_1+\cdots+C_n$ is the pull--back of $o$, 
hence linearly equivalent to 0. Since $C$ is connected, a 
$\Z$--linear combination of $C_1,\dots,C_n$ is linearly equivalent to 0 
if and only if it is a multiple of $C_1+\cdots+C_n$.

If $\Lcal$ is an invertible extension to $S$ 
of the canonical sheaf on 
the generic fiber of 
$\pi$, then
\[
\Lcal\cong\bigomega_\pi(t_1C_1+\cdots+t_nC_n)
\]
for certain integers $t_1,\dots,t_n$. 
We say that $\Lcal$ is a \emph{canonical sheaf of $\pi$}. 
In particular, $\w_\pi$ is a canonical sheaf.

For each canonical sheaf $\Lcal$ 
of
$\pi$, 
the base--change map 
$\varrho_\Lcal\:H^0(\Lcal)(o)\to H^0(\Lcal|_C)$ 
is injective. We say that $(\Image(\varrho_\Lcal),\Lcal|_C)$ is a 
\emph{limit canonical system of $\pi$}. 
\end{subsct}

\begin{theorem}\setcounter{equation}{0}
Let $C$ be a nodal curve of arithmetic genus $g>0$ and 
$C_1,\dots,C_n$ its irreducible components. 
Let 
$\pi\:S\to B$ 
be a regular smoothing of $C$ and $o$ 
the special point 
of $B$. Then for each $i=1,\dots,n$ there is a unique canonical sheaf 
$\Lcal_i$ of $\pi$ 
meeting 
the following two conditions.
\begin{enumerate}
\item The natural map $\varrho_i\:H^0(\Lcal_i)(o)\to H^0(\Lcal_i|_{C_i})$ 
is injective.
\item The natural map $\varrho_{i,j}\:H^0(\Lcal_i)(o)\to H^0(\Lcal_i|_{C_j})$ 
is not zero if $j\neq i$.
\end{enumerate}
\end{theorem}

\begin{proof}
As in \cite[Thm. 1, p. 23 and Prop. 2, p. 26]{Esteves2}. 
See \cite{Ran} as well.
\end{proof}

\begin{definition}\setcounter{equation}{0} 
Keep the set--up of Theorem~2.2. For each $i=1,\dots,n$ we say that 
$\Lcal_i$ is the \emph{canonical sheaf of $\pi$ with focus on $C_i$}. 
The associated limit canonical system is also said to have focus on $C_i$. 
In addition, we 
call the linear system $(\Image(\varrho_i),\Lcal_i|_{C_i})$ the 
\emph{limit canonical aspect of $\pi$ with focus on $C_i$}. 

For each $i=1,\dots,n$ there are unique integers $t_{i,1},\dots,t_{i,n}$
such that $t_{i,i}=0$ and
\[
\Lcal_i\cong\w_\pi(\sum_{j=1}^{n}t_{i,j}C_j).
\]
For each reducible 
node $p$ of $C$ in $C_i$ let $\ell_i(p):=t_{i,j}$,
where $j$ is the unique integer such that $j\neq i$ and $p\in C_j$. 
Call $\ell_i(p)$ the \emph{correction number for $\Lcal_i$ at $p$}.
\end{definition}

\begin{proposition}\setcounter{equation}{0} 
Let $C$ be a nodal curve of arithmetic genus $g>0$ and 
$C_1,\dots,C_n$ its irreducible components. 
Let $\pi\:S\to B$ be a regular smoothing of $C$. 
If $\Mcal$ is 
a canonical sheaf of $\pi$ such that the restriction 
$H^0(\Mcal)\to H^0(\Mcal|_{C_i})$ is not zero for every 
$i=1,\dots,n$, then
for each $i=1,\dots,n$ there are unique non--negative integers 
$t_{i,1},\dots,t_{i,n}$ such that $t_{i,i}=0$ and
\[
\Lcal_i \cong \Mcal(\sum_{j=1}^{n}t_{i,j}C_j),
\]
where $\Lcal_i$ is the canonical sheaf of $\pi$ with focus on $C_i$. 
\end{proposition}

\begin{proof} 
As in \cite[Prop.~4, p. 27]{Esteves2}.
\end{proof}

\begin{remark}\setcounter{equation}{0} 
Keep the set--up of Proposition~2.4. For each 
$i=1,\dots,n$ let $e_i$ be the minimum non--negative integer such that 
the restriction,
\[
H^0(\w_\pi(-e_iC_i))\longrightarrow H^0(\w_\pi(-e_iC_i)|_{C_i}),
\]
is not zero. Since forming $\w_\pi$ commutes with changing $B$, we 
have $e_i>0$ if and only if $C_i$ is smooth and rational and 
there are as many connected components of $\ol{C-C_i}$ as
points in $C_i\cap\ol{C-C_i}$. (In this case, $C_i$ is called 
a \emph{loosely connected rational tail}; see \cite[Def.~3.2,~p.~75]{C}.) 

Let 
\[
\Mcal:=\w_\pi(-\sum_{i=1}^ne_iC_i).
\]
By 
Proposition~2.4, 
for each $i=1,\dots,n$, there are non--negative integers 
$t_{i,1},\dots,t_{i,n}$ such that $t_{i,i}=0$ and
\[
\Lcal_i\cong\Mcal(\sum_{j=1}^{n}t_{i,j}C_j).
\]
If $p\in C_i\cap C_j$ for $j\neq i$ then $t_{i,j}\leq 2g-2$ because 
$\deg\Lcal_i|_{C_i}\leq 2g-2$. Observe that $\ell_i(p)=t_{i,j}-e_j+e_i$. 
Thus $\ell_i(p)\leq 2g-2$ if $C_i$ is not a loosely connected rational 
tail.
\end{remark}

\begin{definition}\setcounter{equation}{0} 
Let $\pi\:S\to B$ be a smoothing of a curve $C$. The Weierstrass 
scheme on the generic fiber of $\pi$ extends to a unique 
closed subscheme 
$\mathcal W\subseteq S$ that is flat over $B$. We say that the special 
fiber $W$ of $\mathcal W$ is the \emph{limit Weierstrass scheme of 
$\pi$}, and call 
the associated cycle $[W]$ the 
\emph{limit Weierstrass divisor of $\pi$}.
\end{definition}

\begin{subsct}\setcounter{equation}{0}
\emph{Semi--stable reduction.}
Let $C$ be a nodal curve and 
$\pi\:S\to B$ 
a smoothing of $C$. 
If $p$ is a node of $C$, then the local equation for $\pi$ at $p$ is 
of the form $t^{\mu_p}=h_1h_2$, where $\mu_p$ is a positive integer, 
$t$ is a local parameter of $B$ at the special point, and 
$h_1$ and $h_2$ are the local equations of the branches of $C$ at $p$. 
The integer $\mu_p$ is called the 
\emph{singularity type of $\pi$ at $p$}.
The smoothing is regular if and only if its singularity type at 
every reducible node of $C$ is 1.

If $S$ is not regular, then we resolve its singularities at the 
reducible nodes of $C$ by blowing up successively. We obtain a regular 
smoothing $\wt\pi\:\wt S\to B$ 
whose generic fiber is equal to that of 
$\pi$, and whose special fiber $\wt C$ is the curve obtained from $C$ by 
splitting the branches of $C$ at each reducible node $p$, and 
connecting them by a chain $Z_p$ of 
$\mu_p-1$ rational, smooth curves $Z_{p,1},\dots,Z_{p,\mu_p-1}$,
where $\mu_p$ is the singularity type of $\pi$ at $p$
(see Figure~1).
We say that the curve $\wt C$ described above is 
a \emph{semi--stable model of $C$},
and call $\wt\pi$ the \emph{semi--stable reduction of $\pi$}. 
We'll view the irreducible components of $C$ as irreducible components 
of $\wt C$.

Conversely, let $\wt\pi'\:\wt S'\to B'$ be any regular 
smoothing of the curve $\wt C$ described above. Blow down the 
chains $Z_p$ of rational, smooth curves, and let $S'$ denote 
the ensuing surface. Since the only curves contracted lie on the special 
fiber of $\wt\pi$, the map $\wt\pi'$ descends to a 
map $\pi'\:S'\to B'$. The self--intersection of 
each $Z_{p,j}$ in $\wt S'$ is $-2$ because $\wt S'$ is regular 
along $Z_{p,j}$. So, each chain $Z_p$ 
is blown down to a node, and we recover the curve $C$ as 
the special fiber of $\pi'$. Thus $\pi'$ is a smoothing of $C$.

Let 
$\phi\:\wt S\to S$ 
denote the birational map. Let $\w$ and $\wt\w$ 
denote the (relative) dualizing sheaves of $\pi$ and $\wt\pi$, 
respectively. Since $C$ is Gorenstein, 
$\w$ and $\wt\w$ are 
invertible. Since $\wt\w$ and $\phi^*\w$ coincide on the generic fiber of 
$\wt\pi$, and on every irreducible component of $\wt C$, 
the two sheaves are equal. Hence, by the projection formula, 
$\w=\phi_*\wt\w$. So, if $W$ and $\wt W$ are the limit Weierstrass 
schemes of $\pi$ and $\wt\pi$ then $[W]=\phi_*[\wt W]$.

Therefore, to determine all limit Weierstrass divisors 
on $C$ we need only determine the limit Weierstrass divisor 
of each regular smoothing 
of each semi--stable model 
of $C$.
\end{subsct}

\begin{theorem}
\setcounter{equation}{0}
Let $C$ be a nodal curve of arithmetic genus $g>0$ 
and $C_1,\dots,C_n$ its irreducible components. 
Let $N$ be the set of 
reducible
nodes of $C$. 
Let $\pi$ be a smoothing of $C$ and $\wt\pi$ 
its semi--stable reduction. Let $\wt C$ be the special fiber of 
$\wt\pi$. For each $i=1,\dots,n$, let $R_i$ be the ramification 
divisor of the limit canonical aspect of $\wt\pi$ with focus on $C_i$. 
For each $i=1,\dots,n$ and each $p\in N\cap C_i$ let 
$\ell_i(p)$ be the correction number 
for the canonical sheaf of $\wt\pi$ with focus on $C_i$
at the point of $\wt C$ in $C_i$ 
lying above 
$p$. 
Then the limit Weierstrass scheme $W$ of $\pi$ satisfies
\[
[W]=\sum_{i=1}^nR_i+\sum_{i<j}\sum_{p\in C_i\cap C_j}
g(g-1-\ell_i(p)-\ell_j(p))p.
\]
\end{theorem}

\begin{proof} Let $p\in N$, say $p\in C_i\cap C_j$ 
where $1\leq i<j\leq n$. Let $\mu_p$ be the singularity type 
of $\pi$ at $p$. Let $Z_{p,1},\dots,Z_{p,\mu_p-1}$ be the chain of 
rational,
smooth curves on $\wt C$ that 
are
contracted to $p$ 
in $C$.
Assume that $Z_{p,1}$ intersects $C_i$ and $Z_{p,\mu_p-1}$ intersects $C_j$. 
Put $Z_{p,0}:=C_i$ and $Z_{p,\mu_p}:=C_j$. For each $h=0,\dots,\mu_p-1$ let 
$z_{p,h}$ be the point of intersection of $Z_{p,h}$ and 
$Z_{p,h+1}$ (see Figure~1).

For each $h=0,\dots,\mu_p$ 
let $\Ncal_{p,h}$ denote the canonical sheaf of $\wt\pi$ with focus on 
$Z_{p,h}$. For each $h=1,\dots,\mu_p-1$ let $\lambda_{p,h}$ and $\nu_{p,h}$ 
be the correction numbers for $\Ncal_{p,h}$ at $z_{p,h-1}$ and $z_{p,h}$, 
respectively. Since $Z_{p,h}$ is smooth and rational, 
and its self--intersection is $-2$, it follows that 
$\Ncal_{p,h}|_{Z_{p,h}}$ has degree $\lambda_{p,h}+\nu_{p,h}$. 
Hence, by 
Pl\"ucker's formula, the 
ramification divisor $R_{p,h}$ of the limit canonical aspect of $\wt\pi$ 
with focus on $Z_{p,h}$ satisfies
\begin{equation}
\deg R_{p,h}=g(\lambda_{p,h}+\nu_{p,h}-(g-1)).
\end{equation}
Put $\nu_{p,0}:=\ell_i(p)$ and $\lambda_{p,\mu_p}:=\ell_j(p)$. Then 
\[
\Ncal_{p,h+1}\cong
\Ncal_{p,h}((\lambda_{p,h+1}+\nu_{p,h})Z_{p,h}+0 Z_{p,h+1}+\cdots)
\]
for $h=0,\dots,\mu_p-1$.
 
By 
\cite[Thm. 7, p. 30]{Esteves2}
the limit Weierstrass scheme $\wt W$ of $\wt\pi$ 
satisfies
\begin{equation}
[\wt W]=\sum_{i=1}^nR_i+\sum_{p\in N}\Big(\sum_{h=1}^{\mu_p-1}R_{p,h}+
\sum_{h=0}^{\mu_p-1}g(g-1-\lambda_{p,h+1}-\nu_{p,h})z_{p,h}\Big).
\end{equation}
Now, $[W]$ is the push--forward of $[\wt W]$ to $C$. Combining (2.8.1) and 
(2.8.2) we get the expression for $[W]$ claimed.
\end{proof}

\begin{remark}\setcounter{equation}{0} 
Keep the set--up of Theorem~2.8. 
For each $i=1,\dots,n$ let $\w_i$ be the dualizing sheaf of $C_i$. 
Assume that for each $i=1,\dots,n$ 
there is an 
integer $d_i$ such that $\ell_i(p)\leq d_i$ for each $p\in N\cap C_i$. 
(For instance, we may 
let
$d_i:=2g-2$ for each $i=1,\dots,n$ if no 
irreducible component of $C$ is a loosely connected rational tail; see 
Remark~2.5.) 
Then, instead of viewing the limit canonical aspect of $\wt\pi$ with focus 
on $C_i$ as a system with sections in 
$\w_i(\textstyle{\sum}_{p\in N\cap C_i}(1+\ell_i(p))p)$, we may view it as 
having sections in $\w_i((1+d_i)\textstyle{\sum}_{p\in N\cap C_i}p)$.
Let $R_i(d_i)$ 
be the corresponding ramification divisor. Then
\[
[W]=\sum_{i=1}^nR_i(d_i)+\sum_{i<j}g(g-1-d_i-d_j)\sum_{p\in C_i\cap C_j}p.
\]
\end{remark}

\section{A numerical lemma}

\begin{lemma}\setcounter{equation}{0} 
Let $\Delta$ be a 
non--empty 
finite set and $\mu\in\Q_\Delta^+$. For 
each
$\upsilon\in\Z$ there are unique $\alpha\in\Z_\Delta$ and 
$\rho\in\Q_\Delta$ satisfying the following four conditions.
\begin{subequations}
\begin{align}
&\text{$0<\rho_p\leq\mu_p$ for every $p\in\Delta$,}\\
&\text{$I:=\{p\in\Delta\,|\,\rho_p=\mu_p\}$ is non--empty,}\\
&\text{$\upsilon\leq|\alpha|<\upsilon+|I|$,}\\
&\text{$\mu_p(\alpha_p+1)-\rho_p=\mu_q(\alpha_q+1)-\rho_q$ 
for all $p,q\in\Delta$.}
\end{align}
\end{subequations}
If $\upsilon\geq 1-|\Delta|$ then $\alpha\geq 0$. 
If $\alpha\geq 0$ but $\alpha_p=0$ for some $p\in I$ then $\alpha=0$ and 
$I=\Delta$.
\end{lemma}

\begin{proof} The existence of 
$\alpha\in\Z_\Delta$ 
and 
$\rho\in\Q_\Delta$ 
satisfying Conditions~(3.1.1) is equivalent to the existence of 
$\beta\in\Z_\Delta$
and 
$\sigma\in\Q_\Delta$
satisfying the following four conditions.
\begin{subequations}
\begin{align}
&\text{$0\leq\sigma_p\leq\mu_p$ for every $p\in\Delta$,}\\
&\text{$J:=\{p\in\Delta\,|\, \sigma_p=0\}$ is non--empty,}\\
&\text{$|\beta|=\upsilon$,}\\
&\text{$\mu_p\beta_p+\sigma_p=\mu_q\beta_q+\sigma_q$ for all 
$p,q\in\Delta$.}
\end{align}
\end{subequations}
In fact, suppose that 
$\beta\in\Z_\Delta$
and 
$\sigma\in\Q_\Delta$
satisfy Conditions (3.1.2). Define 
$\alpha\in\Z_\Delta$
and 
$\rho\in\Q_\Delta$
by letting for each $p\in\Delta$,
\[
(\alpha_p,\rho_p):=
\begin{cases}
(\beta_p+1,\mu_p),&\text{if $\sigma_p=\mu_p$,}\\
(\beta_p,\mu_p-\sigma_p),&\text{otherwise.}
\end{cases}
\]
Then $\alpha$ and $\rho$ satisfy Conditions (3.1.1). Conversely, 
suppose that 
$\alpha\in\Z_\Delta$ 
and 
$\rho\in\Q_\Delta$
satisfy Conditions (3.1.1). By (3.1.1c), there is a 
proper 
subset $K$ of $I$ 
with $|K|=|\alpha|-\upsilon$. Define 
$\beta\in\Z_\Delta$ 
and 
$\sigma\in\Q_\Delta$ 
by letting for each $p\in\Delta$,
\[
(\beta_p,\sigma_p):=
\begin{cases}
(\alpha_p-1,\mu_p),&\text{if $p\in K$,}\\
(\alpha_p,\mu_p-\rho_p),&\text{otherwise.}
\end{cases}
\]
Then $\beta$ and $\sigma$ satisfy Conditions (3.1.2).

Let's prove 
now the existence of $\beta\in\Z_\Delta$ and $\sigma\in\Q_\Delta$ satisfying 
Conditions (3.1.2). 
Let $\Gamma$ be the set of points 
$x\in\R_\Delta$ 
such that
\begin{equation}
x_p=\mu_p\beta_p+\sigma_p\quad\text{for each }p\in\Delta,
\end{equation}
where 
$\beta\in\Z_\Delta$ 
and 
$\sigma\in\R_\Delta$ 
satisfy
Conditions (3.1.2a--c). Let 
$$
D:=\{y\in
\R_\Delta
\,|\, y_p=y_q\text{ for all }p,q\in\Delta\}.
$$
If $\Gamma$ meets $D$, then there are 
$\beta\in\Z_\Delta$ 
and 
$\sigma\in\R_\Delta$
satisfying Conditions (3.1.2). In addition, 
$\sigma\in\Q_\Delta$ 
by (3.1.2b) and (3.1.2d). 
Thus, it is enough to show that 
$\Gamma$ meets $D$. 

If $x\in\Gamma$ then 
$\upsilon\leq\sum_{p\in\Delta}x_p/\mu_p<\upsilon+\delta$ by 
(3.1.2a--c) and (3.1.3). Since $D\not\subseteq P$, where
\[
P:=\{y\in
\R_\Delta
\,|\,\textstyle\sum_{p\in\Delta}y_p/\mu_p=0\},
\] 
there is a point $x\in\Gamma$ closest to $D$. Fix this 
point $x$,
and let 
$\beta\in\Z_\Delta$ and 
$\sigma\in\R_\Delta$ 
satisfying (3.1.2a--c) and 
(3.1.3). We'll show that $x\in D$.

Suppose $x\not\in D$. 
Let 
$d\:\R_\Delta\to\R$ 
be the distance function to $D$. Since the 
coordinate axes are not parallel to 
$D$ and $x\not\in D$, we have 
$\partial_p d(x)\neq 0$ 
for every $p\in\Delta$. If $0<\sigma_p<\mu_p$ for a certain 
$p\in\Delta$, 
then we 
would be able to produce a point of 
$\Gamma$ closer to $D$ by making $\sigma_p$ vary. So, either 
$\sigma_p=0$ or $\sigma_p=\mu_p$ for each 
$p\in\Delta$. 
In addition, 
$\partial_p d(x)>0$ if $\sigma_p=0$ and $\partial_p d(x)<0$ if 
$\sigma_p=\mu_p$.

Suppose that $\sigma_p=0$ and $\sigma_q=\mu_q$ for certain 
$p,q\in\Delta$.
Define 
$\beta'\in\Z_\Delta$ 
and 
$\sigma'\in\Q_\Delta$ 
by letting for each 
$r\in\Delta$,
\[
(\beta'_r,\sigma'_r):=
\begin{cases}
(\beta_p-1,\mu_p),&\text{if $r=p$,}\\
(\beta_q+1,0),&\text{if $r=q$,}\\
(\beta_r,\sigma_r),&\text{otherwise.}
\end{cases}
\]
Then $x_r=\mu_r\beta'_r+\sigma'_r$ for every $r\in\Delta$. However, 
since $\partial_q d(x)<0$ and $\sigma'_q=0$, we would get a 
point of $\Gamma$ closer to $D$ by increasing $\sigma'_q$, 
reaching a contradiction.
Therefore, either 
$\sigma=0$ or $\sigma=\mu$. By (3.1.2b), $\sigma=0$, and thus 
$\partial_p d(x)>0$ for every $p\in\Delta$. 

Now, let $y\in D$ be the closest point to $x$, and $z:=x-y$. Then 
$\sum_{p\in\Delta}z_p=0$. But $z_p\geq 0$ 
for every $p\in\Delta$ because $\partial_p d(x)>0$.
Hence $z=0$ and 
$x\in D$, reaching a contradiction. Thus $x\in D$, and hence $\beta$ and 
$\sigma$ satisfy Conditions~(3.1.2).

Let's prove now uniqueness of 
$\alpha\in\Z_\Delta$ 
and 
$\rho\in\Q_\Delta$ 
meeting Conditions~(3.1.1). Let 
$\alpha'\in\Z_\Delta$ 
and 
$\rho'\in\Q_\Delta$ 
satisfying the same 
conditions as $\alpha$ and $\rho$.
Suppose 
$\alpha_p>\alpha'_p$ for some $p\in\Delta$. 
So,
\begin{equation}
\mu_q(\alpha'_q+1)-\rho'_q=\mu_p(\alpha'_p+1)-\rho'_p<
\mu_p(\alpha_p+1)-\rho_p=\mu_q(\alpha_q+1)-\rho_q,
\end{equation}
and hence $\alpha_q\geq\alpha'_q$ for every 
$q\in\Delta$,
with equality only if 
$\rho'_q>\rho_q$. Thus $\alpha_q>\alpha'_q$ for every $q\in I$. Then
\[
|\alpha|\geq|\alpha'|+|I|\geq\upsilon+|I|,
\]
reaching a contradiction with (3.1.1c). So $\alpha\leq\alpha'$. 
Interchanging 
the roles of $\alpha$ and $\alpha'$, we get $\alpha=\alpha'$. 

Suppose now that $\rho'_p>\rho_p$ for some 
$p\in\Delta$. 
So (3.1.4) 
holds, and hence $\rho'_q>\rho_q$ for every 
$q\in\Delta$.
Since 
$\rho_q=\mu_q$ and $\mu_q\geq\rho'_q$ 
for every $q\in I$, we 
reach 
a 
contradiction. So, $\rho'_p\leq\rho_p$ for every 
$p\in\Delta$. 
Interchanging the 
roles of $\rho$ and $\rho'$, we get $\rho=\rho'$.

Finally, if $\alpha_p<0$ for some 
$p\in\Delta$, 
then $\alpha_q<0$ for every 
$q\in\Delta$ 
by Condition~(3.1.1d). If so, $\upsilon\leq -|\Delta|$ by 
Condition~(3.1.1c). Now, if $\alpha\geq 0$ but $\alpha_p=0$ for some 
$p\in I$, then $\alpha=0$ and $I=\Delta$ by Condition~(3.1.1d).
\end{proof}

\begin{definition}\setcounter{equation}{0} 
Keep the set--up of Lemma~3.1. We call $(\alpha,\rho,I)$, and sometimes 
only $(\alpha,I)$, the 
\emph{numerical data associated to $\mu$ and $\upsilon$}.
\end{definition}

\section{Set--up} 

\begin{subsct}\setcounter{equation}{0}
\emph{Two--component curves.} 
Let $C$ be a nodal curve with only two 
irreducible components, denoted $X$ 
and $Y$. Let $\Delta:=X\cap Y$ and $\delta:=|\Delta|$. To avoid 
exceptional cases, we shall 
assume that $\delta>1$ or $g_Xg_Y>0$.
(If also the 
irreducible 
singularities of $C$ are 
nodes, then we are assuming that $C$ is semi--stable.)

For each $p\in\Delta$, let $x_p$ and $y_p$ denote the points of $X$ and 
$Y$ lying over $p$. Let $g_X$ and $g_Y$ be the arithmetic genera of $X$ 
and $Y$. Then the arithmetic genus $g$ of $C$ is given by
\[
     g = g_X + g_Y + \delta-1.
\]
Let $\w_X$ and $\w_Y$ be the dualizing sheaves of $X$ and $Y$, and $\w$ 
that of $C$.
Let
\[
\mathbb L:=\w_X((1+g_Y)\sum_{p\in\Delta}x_p)\quad \text{ and } \quad
\mathbb M:=\w_Y((1+g_X)\sum_{p\in\Delta}y_p).
\]
Let
$\mathbb G_X:=\Grass_g(H^0(\mathbb L))$ and 
$\mathbb G_Y:=\Grass_g(H^0(\mathbb M))$
denote the Grassmannians of $g$--dimensional 
vector subspaces of $H^0(\mathbb L)$ 
and $H^0(\mathbb M)$. Set $\mathbb G:=\mathbb G_X\times\mathbb G_Y$.

For each 
$\mu\in\Q_\Delta^+$,
let $(\alpha_\mu,\rho_\mu,I_\mu)$ be the 
numerical data associated with $\mu$ and $g_Y$, and 
$(\beta_\mu,\sigma'_\mu,J_\mu)$ the numerical data associated 
with $\mu$ and $g_X$; see Definition~3.2. 
Set $\sigma_\mu:=\mu-\sigma'_\mu$.

Note that 
$\rho_\mu,\sigma_\mu\in\Z_\Delta$
if 
$\mu\in\Z_\Delta^+$.
In this case, let $\gamma_\mu:=\alpha_{\mu,p}\mu_p$ for 
(any) $p\in I_\mu$ and $\epsilon_\mu:=\beta_{\mu,p}\mu_p$ for (any)
$p\in J_\mu$. If $g_Y>0$ then 
$\gamma_\mu>0$ by Lemma~3.1. Analogously, if 
$g_X>0$ then $\epsilon_\mu>0$.

If $g_Xg_Y>0$ let $t\in\Z^+$ be such that
$t\mu\in\Z^+_\Delta$ and define
$\wt\alpha_\mu:=\gamma_{t\mu}/\gcd(\gamma_{t\mu},\epsilon_{t\mu})$ and 
$\wt\beta_\mu:=\epsilon_{t\mu}/\gcd(\gamma_{t\mu},\epsilon_{t\mu})$. 
Note that $\wt\alpha_\mu$ and $\wt\beta_\mu$
do not depend on the choice of $t$.
In addition, if $p\in I_\mu\cap J_\mu$ then 
$\wt\alpha_\mu=\alpha_{\mu,p}/\gcd(\alpha_{\mu,p},\beta_{\mu,p})$ 
and $\wt\beta_\mu=\beta_{\mu,p}/\gcd(\alpha_{\mu,p},\beta_{\mu,p})$.
\end{subsct}

\begin{subsct}\setcounter{equation}{0}
\emph{Semi--stable reduction.} 
Preserve~4.1. 
Fix $\mu\in\Z_\Delta^+$. 
Abbreviate $\alpha_\mu$, $\rho_\mu$, $I_\mu$, 
$\beta_\mu$, $\sigma_\mu$, $J_\mu$, $\gamma_\mu$, $\epsilon_\mu$, 
$\wt\alpha_\mu$ and $\wt\beta_\mu$ 
by $\alpha$, $\rho$, $I$, $\beta$, $\sigma$, 
$J$, $\gamma$, $\epsilon$, 
$\wt\alpha$ and $\wt\beta$.
Let
\begin{equation}
\begin{cases}
L_X:=\w_X(\textstyle\sum_{p\in\Delta}(1+\alpha_p)x_p),&\\
L_Y:=
\w_Y(\textstyle\sum_{p\in I}y_p-\textstyle\sum_{p\in\Delta}\alpha_py_p),&\\
M_Y:=\w_Y(\textstyle\sum_{p\in\Delta}(1+\beta_p)y_p),&\\
M_X:=\w_X(\textstyle\sum_{p\in J}x_p-\textstyle\sum_{p\in\Delta}\beta_px_p).
\end{cases}
\end{equation}
By Lemma~3.1,
$\alpha_p\geq 0$ for every $p\in\Delta$. In addition, 
$\alpha_p\leq g_Y$ for every $p\in\Delta$. Indeed, if $g_Y=0$ then 
$\alpha=0$. Now, if $g_Y>0$ then $\alpha_p>0$ for every $p\in I$; 
hence $\alpha_p\leq g_Y$ for every $p\in\Delta$ because $|\alpha|<g_Y+|I|$. 
Likewise, $0\leq\beta_p\leq g_X$ for every $p\in\Delta$. So we may (and 
eventually will) view $L_X$ and $M_Y$ as subsheaves of $\mathbb L$ and 
$\mathbb M$, respectively.

Denote by $\wt C$ the curve obtained from $C$ by splitting its branches 
at each $p\in\Delta$ and connecting them by a chain,
\[
Z_p:=Z_{p,1}\cup\cdots\cup Z_{p,\mu_p-1},
\]
of $\mu_p-1$ rational, smooth curves $Z_{p,j}$, as depicted in 
Figures~1~and~2. 
\begin{figure}[ht]
\begin{center}
\includegraphics[angle=0,height=5.5cm]{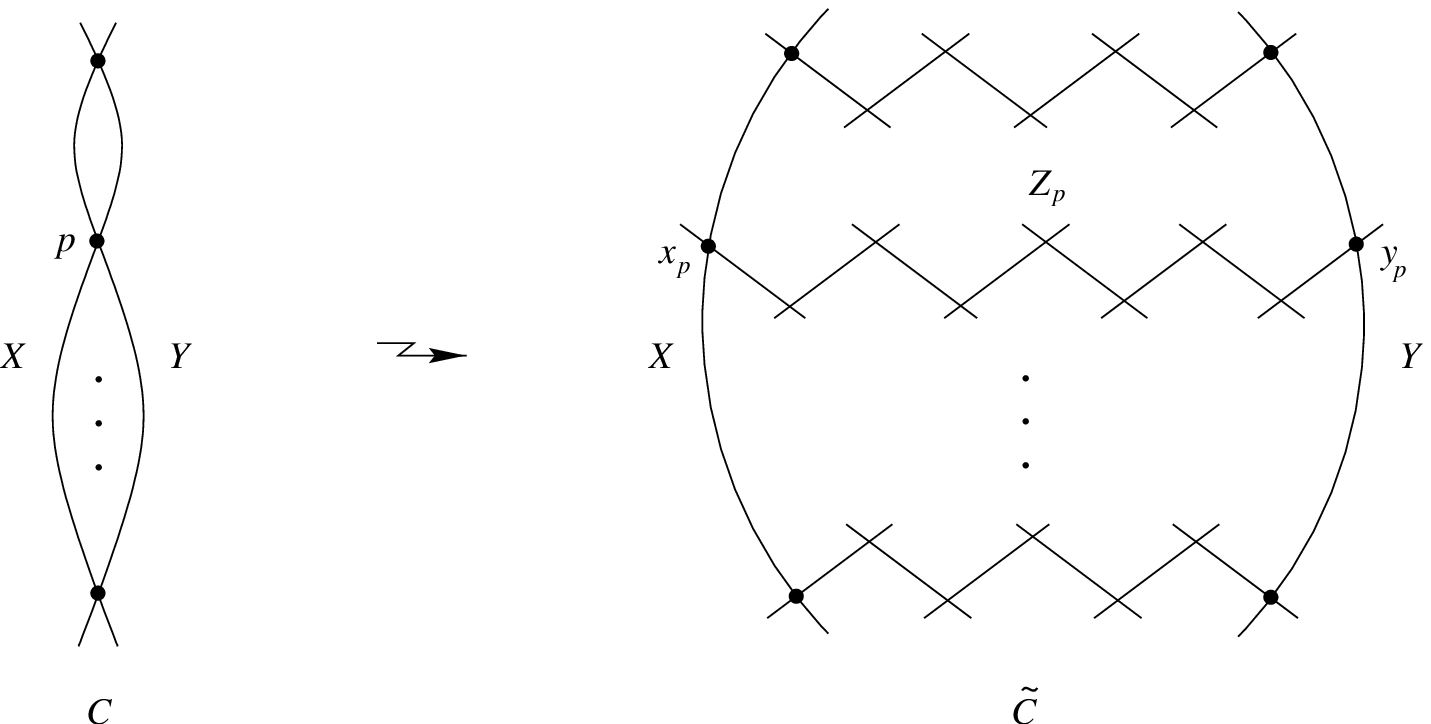}
\caption{The curve $\wt C$.}
\end{center}
\end{figure}
By convention, the leftmost curve in each chain $Z_p$ is $Z_{p,1}$ 
and the rightmost is $Z_{p,\mu_p-1}$. In addition, set $Z_{p,0}:=X$ and 
$Z_{p,\mu_p}:=Y$ for 
each
$p\in\Delta$. For 
each
$p\in\Delta$ and 
each
$j=0,\dots,\mu_{p-1}$, let $z_{p,j}$ be the point of 
$Z_{p,j}\cap Z_{p,j+1}$, as depicted in Figure~1.
So $z_{p,0}=x_p$ and $z_{p,\mu_p-1}=y_p$ for every $p\in\Delta$. 
If $\mu_p=1$ for a certain $p\in\Delta$, 
then $Z_p$ is empty; in other words, the branches of $C$ at $p$ are 
not split in $\wt C$.

Let $\wt\w$ denote the dualizing sheaf 
of
$\wt C$. As observed in 2.7, 
$\wt\w$ is the pull--back of $\w$ to $\wt C$ and $\w$ is the push--forward 
of 
$\wt\w$ to $C$.

For each $p\in\Delta$, 
each 
$j=1,\dots,\mu_p-1$ and 
each
$m\in\Z$, define the 
formal sums:
\begin{align*}
&Z_p^{(m)}:=mZ_{p,1}+2mZ_{p,2}+\cdots+(\mu_p-1)mZ_{p,\mu_p-1},\\
&\wh Z_p^{(m)}:=(\mu_p-1)mZ_{p,1}+(\mu_p-2)m Z_{p,2}+\cdots+mZ_{p,\mu_p-1},\\
&Z_p^{(m,j)}:=Z_p^{(m)}+Z_{p,j+1}+2Z_{p,j+2}+\cdots+
(\mu_p-1-j)Z_{p,\mu_p-1},\\
&\wh Z_p^{(m,j)}:=\wh Z_p^{(m)}+Z_{p,j-1}+2Z_{p,j-2}+\cdots+(j-1)Z_{p,1}.
\end{align*}
Set $Z_p^{(m,\mu_p)}:=Z_p^{(m)}$ and $\wh Z_p^{(m,0)}:=\wh Z_p^{(m)}$.
\end{subsct}

\begin{subsct}\setcounter{equation}{0}
\emph{General--position conditions.}
We shall eventually assume general--position conditions on effective 
divisors $D$ of $C$ supported in $\Delta$. Typically, if 4.1 is preserved 
we shall consider both of the following two conditions:
\begin{subequations}
\begin{align}
&h^0(\w_X(-D))=0\text{ if }D\text{ is effective},\  \deg D=g_X\text{ and }
\text{Supp}(D)\subseteq\{x_p|\  p\in\Delta\},\\
&h^0(\w_Y(-D))=0\text{ if }D\text{ is effective},\  \deg D=g_Y\text{ and }
\text{Supp}(D)\subseteq\{y_p|\  p\in\Delta\},
\end{align}
\end{subequations}
whereas if also 4.2 is preserved 
we shall consider one or both of the following 
conditions:
\begin{subequations}
\begin{align}
&h^0(L_Y(-D))=0\text{ if }\deg D=g_Y+|I|-|\alpha|\text{ and }
0\leq D\leq\textstyle\sum_{p\in I}y_p,\\
&h^0(M_X(-D))=0\text{ if }\deg D=g_X+|J|-|\beta|\text{ and }
0\leq D\leq\textstyle\sum_{p\in J}x_p.
\end{align}
\end{subequations}
If Conditions (4.3.1) hold, so do Conditions (4.3.2).
\end{subsct}

\section{Degeneration} 

\begin{lemma}\setcounter{equation}{0}
Preserve~{\rm 4.1--2}. Let $\wt\pi$ be a regular smoothing of 
$\widetilde C$ and $\w_{\wt\pi}$ the dualizing sheaf of $\wt\pi$. Put
\[
\Lcal:=\w_{\wt\pi}(\sum_{p\in\Delta}Z_p^{(\alpha_p,\rho_p)}+\gamma Y).
\]
If {\rm (4.3.2a)} holds, 
so do the following five assertions.
\begin{subequations}
\begin{align}
&h^0(\Lcal|_Y)>0,\\
&h^0(\Lcal|_Y(-\textstyle\sum_{p\in I}y_p))=0,\\
&h^0(\Lcal|_X(-\textstyle\sum_{p\in I}x_p))=g-h^0(\Lcal|_Y),\\
&h^0(\Lcal|_{Z_p}(-x_p-y_p))=0\text{ for every }p\in\Delta,\\
&h^0(\Lcal|_{Z_p}(-x_p))=h^0(\Lcal|_{Z_p}(-y_p))=0\text{ for every }p\in I.
\end{align}
\end{subequations}
\end{lemma}

\begin{proof} Note first that
\begin{equation}
\Lcal|_X\cong L_X\text{ and }\Lcal|_Y\cong L_Y,
\end{equation}
where the second isomorphism holds because of (3.1.1a,d).
Now, for each $p\in\Delta$ and each $j=1,\dots,\mu_p-1$, let $\w_{Z_{p,j}}$ 
denote the dualizing sheaf of $Z_{p,j}$. Then
\[
\Lcal|_{Z_{p,j}}\cong
\begin{cases}
\w_{Z_{p,j}}((1+\alpha_p)z_{p,j}+(1-\alpha_p)z_{p,j-1})& 
\text{if $j<\rho_p$,}\\
\w_{Z_{p,j}}((2+\alpha_p)z_{p,j}+(1-\alpha_p)z_{p,j-1})& 
\text{if $j=\rho_p$,}\\
\w_{Z_{p,j}}((2+\alpha_p)z_{p,j}-\alpha_pz_{p,j-1})& 
\text{if $j>\rho_p$.}
\end{cases}
\]
Since $Z_{p,j}$ is smooth and rational, 
$\w_{Z_{p,j}} \cong \Ocal_{Z_{p,j}}(-2)$. So,
\begin{equation}
\Lcal|_{Z_{p,j}}\cong
\begin{cases}
\Ocal_{Z_{p,j}}&\text{if $p\in I$ or $j\neq\rho_p$},\\
\Ocal_{Z_{p,j}}(1)&\text{if $p\not\in I$ and $j=\rho_p$.}
\end{cases}
\end{equation}
Equations~(5.1.1d--e) follow immediately.

Since $\alpha\geq 0$ by 
Lemma~3.1,
it follows from (5.1.2) that
\begin{equation}
h^0(\Lcal|_X)=g_X+|\alpha|+\delta-1.
\end{equation}
In addition, 
since $|\alpha|\geq g_Y$ by (3.1.1c), Equation~(5.1.1b) follows from 
(4.3.2a).

By 
Lemma~3.1,
if there is $p\in I$ such 
that $\alpha_p=0$, then $\alpha=0$, and hence $g_Y=0$. Conversely, if 
$g_Y=0$ then $\alpha=0$. So, since $|\alpha|<g_Y+|I|$ 
by (3.1.1c), by (5.1.2) and (4.3.2a),
\begin{equation}
h^0(\Lcal|_Y)=
\begin{cases}
g_Y+|I|-|\alpha|&\text{if $g_Y>0$,}\\
\delta-1&\text{if $g_Y=0$.}
\end{cases}
\end{equation}
In any case, we get (5.1.1a). 

Since $\alpha\geq 0$, with equality only if $g_Y=0$, by (5.1.2) and 
Riemann--Roch,
\begin{equation}
h^0(\Lcal|_X(-\textstyle\sum_{p\in I}x_p))=
\begin{cases}
g_X+|\alpha|-|I|+\delta-1&\text{if $g_Y>0$,}\\
g_X&\text{if $g_Y=0$.}
\end{cases}
\end{equation}
Putting together (5.1.5) and (5.1.6) we get (5.1.1c).
\end{proof}

\begin{theorem}\setcounter{equation}{0}
Preserve~{\rm 4.1--2}. Let $\wt\pi$ be a regular smoothing 
of $\wt C$ and $\pi$ the induced smoothing of $C$. Let 
$\w_{\wt\pi}$
be the dualizing sheaf of $\wt\pi$. Set
\begin{align}
\Lcal:=&\w_{\wt\pi}(\sum_{p\in\Delta}Z_p^{(\alpha_p,\rho_p)}+\gamma Y),\\
\Mcal:=&\w_{\wt\pi}(\sum_{p\in\Delta}\wh Z_p^{(\beta_p,\sigma_p)}+
\epsilon X).
\end{align}
Let $V_X$ and $V_Y$ denote the images of the restriction maps,
\[
\tau_X\:H^0(\Lcal|_{\wt C})\to H^0(\Lcal|_X)\quad\text{and}\quad
\tau_Y\:H^0(\Mcal|_{\wt C})\to H^0(\Mcal|_Y).
\]
Then the following three statements hold.
\begin{enumerate}
\item If {\rm (4.3.2a)} holds, then $\Lcal$ is the canonical sheaf of 
$\wt\pi$ with 
focus on $X$, the map $\tau_X$ is injective, $(V_X,\Lcal|_X)$ is the 
limit canonical aspect of $\wt\pi$ 
with focus on $X$, 
the point $x_p$ is not a base point of the aspect for any $p\in\Delta$,
\[
\text{\rm codim}(V_X,H^0(\Lcal|_X))=|\alpha|-g_Y
\quad\text{and}\quad 
V_X\supseteq H^0(\Lcal|_X(-\textstyle\sum_{p\in I}x_p)).
\]
\item If {\rm (4.3.2b)} holds, then $\Mcal$ is the canonical 
sheaf of $\wt\pi$ with focus on $Y$, the map $\tau_Y$ is injective, 
$(V_Y,\Mcal|_Y)$ is the limit canonical aspect of $\wt\pi$ 
with focus on $Y$, 
the point $y_p$ is not a base point of the aspect for any $p\in\Delta$,
\[
\text{\rm codim}(V_Y,H^0(\Mcal|_Y))=|\beta|-g_X
\quad\text{and}\quad 
V_Y\supseteq H^0(\Mcal|_Y(-\textstyle\sum_{p\in J}y_p)).
\]
\item If {\rm (4.3.2a--b)} hold, then the limit Weierstrass scheme $W$ 
of $\pi$ satisfies
\[
[W]=R_X+R_Y+\sum_{p\in\Delta}g(g-1-\alpha_p-\beta_p)p,
\]
where $R_X$ and $R_Y$ are the ramification divisors of $(V_X,\Lcal|_X)$ and 
$(V_Y,\Mcal|_Y)$.
\end{enumerate}
\end{theorem}

\begin{proof} Let's prove Statement~1. Since (4.3.2a) holds, 
Lemma~5.1 applies.
Let $L:=\Lcal|_{\wt C}$. Let's check first 
that $\tau_X$ is injective. In fact, 
let $s\in H^0(L)$ such that $s|_X=0$. Then
$s|_{Z_p}\in H^0(L|_{Z_p}(-x_p))$ for every $p\in\Delta$. 
So $s|_{Z_p}=0$ for every $p\in I$ by (5.1.1e), and hence 
$s|_Y\in H^0(L|_Y(-\sum_{p\in I}y_p))$. By (5.1.1b), 
$s|_Y=0$ as well. So $s|_{Z_p}\in H^0(L|_{Z_p}(-x_p-y_p))$, and thus 
$s|_{Z_p}=0$ for every $p\in\Delta$ by (5.1.1d). Hence $s=0$. So 
$\tau_X$ is injective.

Let's prove now that $\Lcal$ is the canonical sheaf of $\wt\pi$ 
with focus on $X$. 
We must check that Conditions~1~and~2 of 
Theorem~2.2 are verified. 
Condition~1 holds because $\tau_X$ is injective.

Let's check Condition~2. The limit canonical system of $\wt\pi$ 
associated to 
$\Lcal$ has (affine) rank $g$. Hence, to show that Condition~2 holds it is 
enough to 
prove 
that the kernel 
$N_E$ of the restriction map $H^0(L) \to H^0(L|_E)$ has dimension 
strictly less than $g$ for every irreducible component $E$ of $\wt C$ 
other than $X$. 
We divide the proof in three steps.

\emph{Step 1}:
We show that $\dim N_Y<g$. In fact, if $s\in N_Y$ then 
$s|_{Z_p}\in H^0(L|_{Z_p}(-y_p))$ for each $p\in\Delta$. 
So $s|_{Z_p}=0$ for each $p\in I$ by (5.1.1e),
and hence
$s|_X\in H^0(L|_X(-\sum_{p\in I}x_p))$. If $s|_X=0$
then 
$s|_{Z_p}=0$ for every $p\in\Delta$ by (5.1.1d), and hence $s=0$. So,
\begin{equation}
\dim N_Y\leq h^0(L|_X(-\textstyle\sum_{p\in I}x_p)).
\end{equation}
Thus $\dim N_Y<g$ by (5.1.1a) and (5.1.1c).

\emph{Step 2}: Let $p\in\Delta$ and $j\in\Z$ with $0<j<\mu_p$ 
such that either $p\in I$ or $\rho_p\leq j<\mu_p$.
We show that 
$\dim N_{Z_{p,j}}<g$.
In fact, if $s\in N_{Z_{p,j}}$ then 
$s|_Y\in H^0(L|_Y(-y_p))$ by (5.1.3).
By~(5.1.1e), if $s|_Y=0$ then 
$s|_{Z_q}=0$ for every $q\in I$, and thus 
$s|_X\in H^0(L|_X(-\sum_{q\in I}x_q))$. If $s|_X=0$ as well, 
then $s=0$ by (5.1.1d). So,
\[
\dim N_{Z_{p,j}}\leq h^0(L|_Y(-y_p))+h^0(L|_X(-\textstyle\sum_{q\in I}x_q)).
\]
Since $\Lcal|_Y\cong L_Y$, it follows from (4.3.2a) that 
$h^0(L|_Y(-y_p))=h^0(L|_Y)-1$.
By (5.1.1c), $\dim N_{Z_{p,j}}<g$.

\emph{Step 3}: Let $p\in\Delta-I$ and $j\in\Z$ such that $0<j<\rho_p$. 
We show that 
$\dim N_{Z_{p,j}}<g$.
In fact, if $s\in N_{Z_{p,j}}$ then  
$s|_X\in H^0(L|_X(-x_p))$ by (5.1.3).
If $s|_Y=0$, then $s|_{Z_q}=0$ for every $q\in I$ by (5.1.1e), and thus 
$s|_X\in H^0(L|_X(-x_p-\sum_{q\in I}x_q))$. If $s|_X=0$ as well, then 
$s=0$ by (5.1.1d). So,
\begin{equation}
\dim N_{Z_{p,j}}\leq h^0(L|_Y)+h^0(L|_X(-x_p-\textstyle\sum_{q\in I}x_q)).
\end{equation}
By 
Lemma~3.1, 
$\alpha\geq 0$, with equality only if 
$I=\Delta$. 
So $\alpha\neq 0$ because $p\in\Delta-I$. Since $\Lcal|_X\cong L_X$, we have
\begin{equation}
h^0(L|_X(-x_p-\textstyle\sum_{q\in I}x_q))=
h^0(L|_X(-\textstyle\sum_{q\in I}x_q))-1.
\end{equation}
So $\dim N_{Z_{p,j}}<g$ by (5.1.1c) and (5.2.4--5).

Condition~2 is checked. So $\Lcal$ is the canonical sheaf of $\wt\pi$ 
with focus on $X$.

Let's prove now that 
$(V_X,L|_X)$
is the limit canonical aspect of 
$\wt\pi$ with focus on $X$. 
Since the aspect has rank $g$, it is enough to 
show that $h^0(L)\leq g$. By definition of $N_Y$, we have 
$h^0(L)\leq\dim N_Y+h^0(L|_Y)$. Hence, 
\[
h^0(L)\leq h^0(L|_X(-\textstyle\sum_{p\in I}x_p))+h^0(L|_Y)
\]
by (5.2.3). So $h^0(L)\leq g$ by (5.1.1c). So 
$(V_X,L|_X)$ 
is the 
limit canonical aspect of $\wt\pi$ with focus on $X$.

Consider the remaining assertions in Statement~1. Let's check first that 
\begin{equation}
V_X\supseteq H^0(L|_X(-\textstyle\sum_{p\in I} x_p)).
\end{equation}
In fact, if $s\in N_Y$ then $s|_{Z_p}=0$ for every $p\in I$ by (5.1.1e). 
Therefore,
\begin{equation}
\tau_X(N_Y)\subseteq H^0(L|_X(-\textstyle\sum_{p\in I}x_p)).
\end{equation}
Now, since $h^0(L)=g$, from the definition of $N_Y$ we get 
$\dim N_Y\geq g-h^0(L|_Y)$, and hence 
$\dim N_Y\geq h^0(L|_X(-\sum_{p\in I}x_p))$ by (5.1.1c). So 
equality holds in (5.2.7), and thus (5.2.6) follows.

Let's check now that 
$x_p$ is not a base point of $(V_X,L|_X)$ for any $p\in\Delta$. 
In fact, if there were $p\in I$ such that 
$V_X\subseteq H^0(L|_X(-x_p))$, then the 
restriction $H^0(L)\to H^0(L|_{Z_p})$ would be zero by (5.1.1e), 
thus contradicting Condition~2 of Theorem~2.2. In 
addition, if $V_X\subseteq H^0(L|_X(-x_p))$ for a certain $p\in\Delta-I$, 
then
\[
h^0(L|_X(-\textstyle\sum_{q\in I}x_q))=
h^0(L|_X(-x_p-\textstyle\sum_{q\in I}x_q))
\]
by (5.2.6), thus contradicting (5.2.5). So, no point of 
$\{x_p|\  p\in\Delta\}$ 
is a 
base point of $(V_X,L|_X)$. 

Finally, 
$\text{codim}(V_X,H^0(L|_X))=|\alpha|-g_Y$ because $\dim V_X=g$, and 
because of (5.1.4). The proof of Statement~1 is complete. Statement~2 
follows by analogy.

Let's prove now Statement~3. By Statements~1~and~2, the sheaves 
$\Lcal$ and $\Mcal$ are the canonical sheaves and $(V_X,\Lcal|_X)$ and 
$(V_Y,\Mcal|_Y)$ are the limit canonical aspects of $\wt\pi$ with foci on 
$X$ and $Y$, respectively. By (5.2.1), for each $p\in\Delta$ 
the correction number for $\Lcal$ at $x_p$ is $\alpha_p$. By (5.2.2), 
for each $p\in\Delta$ the correction number for $\Mcal$ at $y_p$ is 
$\beta_p$. Then Statement~3 follows from Theorem~2.8.
\end{proof}

\begin{remark}\setcounter{equation}{0} 
Keep the set--up of Theorem~5.2. Let 
$L_C$ 
be the push--forward to $C$ of $\Lcal|_{\wt C}$ and 
$M_C$
that of 
$\Mcal|_{\wt C}$. It follows from (5.1.3) that the sheaf 
$L_C$ 
is 
torsion--free of generic rank 1, and fails to be invertible precisely 
along $\Delta-I$. Likewise, 
$M_C$
is 
torsion--free of generic rank 1, and fails to be invertible precisely 
along $\Delta-J$.
\end{remark}

\section{Enriched structures}

\begin{subsct}\setcounter{equation}{0}
\emph{Deformation theory.} Let $C$ be a nodal curve
defined over a field $k$.
Since $C$ is one--dimensional, generically smooth, and a local 
complete intersection, $\text{Ext}^2_C(\Omega^1_C,\Ocal_C)=0$ by 
\cite[Lemma~1.3, p. 79]{DM}.
By \cite{Sch}, there exists a versal formal 
deformation $\mathcal V$ of $C$ over the base scheme 
\[
M:=\text{Spec}(k[[t_1,\dots,t_e]]),\text{ where }
e:=\dim_k\text{Ext}^1_C(\Omega^1_C,\Ocal_C).
\]
Since $C$ is a curve, $H^2(C,\Ocal_C)=0$. By \cite{Sch}, the formal scheme 
$\mathcal V$ is the formal completion of a scheme $V$ projective and flat 
over $M$.

Let $p$ be a node of $C$. So $\wh\Ocal_{C,p}\cong k[[u,v]]/(uv)$. 
The versal formal deformation of $\wh\Ocal_{C,p}$ is a complete, 
local $k[[t]]$--algebra $\wh\Ocal$ such that 
$\wh\Ocal/t\wh\Ocal\cong\wh\Ocal_{C,p}$. More precisely, there is an 
isomorphism of $k[[t]]$--algebras,
\[
\wh\Ocal\cong k[[t,u,v]]/(uv-t).
\]

Let $p_1,\dots,p_\delta$ denote the 
nodes of $C$. For 
$i=1,\dots,\delta$, let 
$\wh\Ocal_i$ denote the versal formal deformation of $\wh\Ocal_{C,p_i}$ 
over $k[[t]]$. Since global deformations induce local ones, there is a 
map 
$\phi_i\:k[[t]]\to k[[t_1,\dots,t_e]]$ 
such that 
\[
\wh\Ocal_{V,p_i}\cong\wh\Ocal_i\otimes_{k[[t]]} k[[t_1,\dots,t_e]].
\]
Since $C$ is a curve, by 
\cite[Prop. 1.5, p. 81]{DM}, 
the tensor--product map,
\[
\phi_1\otimes\cdots\otimes\phi_\delta\:k[[t_1,\dots,t_\delta]]
\longrightarrow k[[t_1,\dots,t_e]],
\]
is formally smooth. Up to changing variables, we may assume that 
$\phi_1\otimes\cdots\otimes\phi_\delta$ is the inclusion map. Thus, there 
is an isomorphism of $k[[t_1,\dots,t_e]]$--algebras,
\[
\wh\Ocal_{V,p_i}\cong k[[t_1,\dots,t_e,u_i,v_i]]/(u_iv_i-t_i),
\]
for $i=1,\dots,\delta$.
\end{subsct}

\begin{definition}\setcounter{equation}{0} 
Let $\Upsilon$ be a set and 
$\lambda\in\Z_\Upsilon$. A subset $P$ of $\Upsilon$ is called 
\emph{$\lambda$--balanced} if 
$\lambda_E=\lambda_F$ 
for all $E,F\in P$. 
If $P$ is $\lambda$--balanced, set $\lambda_P:=\lambda_E$ 
for (any) $E\in P$.
\end{definition}

\begin{subsct}\setcounter{equation}{0}\emph{Certain fractional ideals.} 
Let 
$\wt C$
be a nodal curve
defined over a field $k$.
Let $\Upsilon$ denote the set of irreducible 
components of 
$\wt C$.
For 
each 
$P\subseteq\Upsilon$ set
\[
\wt C_P:=\bigcup_{E\in P}E.
\]
We say that 
$\wt C_P$ is a \emph{subcurve of $\wt C$}. 
Let 
$\ol C$ denote the disjoint union of the 
irreducible 
components of 
$\wt C$.

For 
each pair $(E,F)$ of distinct irreducible components of $\wt C$ 
let $\Delta_{E,F}$ denote 
the reduced Weil divisor with support $E\cap F$. We 
shall
view $\Delta_{E,F}$ as a 
Cartier
divisor on every subcurve of 
$\wt C$
or 
$\ol C$ containing either $E$ or $F$ but not both. 

For each $\lambda\in\Z_\Upsilon$ and each $P\subseteq\Upsilon$ 
define
\[
D_{\lambda,P}:=\sum_{E\in P}\sum_{F\not\in P}
(\lambda_F-\lambda_E)\Delta_{E,F}
\quad \text{ and }\quad 
N_{\lambda,P}:=\Ocal_{\wt C_P}(D_{\lambda,P}).
\]
If $\mathbb P$ is a partition of $\Upsilon$, let $\ol C_{\mathbb P}$ be 
the disjoint union of the subcurves 
$\wt C_P$ 
for 
$P\in\mathbb P$, and set
\[
\ol N_{\lambda,\mathbb P}:=\bigoplus_{P\in\mathbb P}N_{\lambda,P}.
\]
Let $\mathbb T$ be the total partition of $\Upsilon$. 
Then $\ol C_{\mathbb T}=\ol C$. 
For each $\lambda\in\Z_\Upsilon$ set 
$\ol N_\lambda:=\ol N_{\lambda,\mathbb T}$. 
Then, for each 
partition $\mathbb P$ of $\Upsilon$ 
in $\lambda$--balanced subsets, 
there is a natural injection 
$\ol N_{\lambda,\mathbb P}\hookrightarrow\ol N_{\lambda}$ satisfying 
$\ol N_{\lambda,\mathbb P}\Ocal_{\ol C}=\ol N_{\lambda}$. 

Note that
\[
\ol N_{\l_1\lambda_1+\l_2\lambda_2}=
\ol N_{\lambda_1}^{\l_1}\ol N_{\lambda_2}^{\l_2}
\]
for all $\l_1,\l_2\in\Z$ and all $\lambda_1,\lambda_2\in\Z_\Upsilon$. So,
$\ol{\mathbb E}:=\{\ol N_\lambda\,|\, \lambda\in\Z_\Upsilon\}$ is 
a subgroup of the group 
$\ol{\mathbb U}$ 
of invertible sheaves of 
fractional ideals 
of
$\ol C$.
\end{subsct}

\begin{subsct}\setcounter{equation}{0}\emph{Enriched structures.}
Preserve~6.3. Let $\wt\pi:\wt S\to B$ be a regular smoothing of $\wt C$. 
Let $t$ be a parameter of $B$ at its special point.
For each $\lambda\in\Z_\Upsilon$ define 
$\Ncal_\lambda:=\Ocal_{\wt S}(\sum_{E\in\Upsilon}\lambda_EE)$. Then, for
each $\lambda$--balanced 
subset $P\subseteq\Upsilon$ define $\rho_{\lambda,P}$ as the composition
\[
\begin{CD}
\rho_{\lambda,P}\:\Ncal_\lambda @>\cdot t^{\lambda_P}>> 
\Ocal_{\wt S}(\sum_{E\not\in P}(\lambda_E-\lambda_P)E)
@>>> N_{\lambda,P},
\end{CD}
\]
where the second map is the natural surjection. 

For each $\lambda\in\Z_\Upsilon$ and 
each partition 
$\mathbb P$ of $\Upsilon$ 
in
$\lambda$--balanced subsets let
\[
\rho_{\lambda,\mathbb P}:=\bigoplus_{P\in\mathbb P}\rho_{\lambda,P}|_{\wt C}
\:\Ncal_\lambda|_{\wt C}\longrightarrow\ol N_{\lambda,\mathbb P}.
\]
Let $N_{\lambda,\mathbb P}$ be the image of $\rho_{\lambda,\mathbb P}$. 
Since each $\rho_{\lambda,P}$ is surjective, 
$\rho_{\lambda,\mathbb P}$ is an isomorphism onto 
$N_{\lambda,\mathbb P}$
and 
$N_{\lambda,\mathbb P}\Ocal_{\ol C_{\mathbb P}}=\ol N_{\lambda,\mathbb P}$. 
Note that $N_{\lambda,\mathbb P}$ does not depend on the choice of $t$.

For each $\lambda\in\Z_\Upsilon$ set
$N_\lambda:=N_{\lambda,\mathbb T}$. 
Then, for each 
partition $\mathbb P$ of $\Upsilon$ 
in
$\lambda$--balanced subsets, 
$\rho_{\lambda,\mathbb T}$ factors as $\rho_{\lambda,\mathbb P}$ 
followed by the natural injection 
$\ol N_{\lambda,\mathbb T}\hookrightarrow\ol N_{\lambda}$. So
$N_{\lambda,\mathbb P}=N_\lambda$.

Let 
$\wt{\mathbb U}$ 
denote the group of invertible sheaves of 
fractional ideals 
of $\wt C$. 
Let 
\[
\mathbb E_{\wt\pi}:=\{N_{\lambda}\,|\, \lambda\in\Z_\Upsilon\}\subseteq
\wt{\mathbb U}.
\]
Note that 
$N_{\l_1\lambda_1+\l_2\lambda_2}=N_{\lambda_1}^{\l_1}N_{\lambda_2}^{\l_2}$ 
for all $\l_1,\l_2\in\Z$ and all $\lambda_1,\lambda_2\in\Z_\Upsilon$. So 
$\mathbb E_{\wt\pi}$ 
is the image of the group homomorphism,
\[
e_{\wt\pi}\:\Z_\Upsilon\longrightarrow\wt{\mathbb U},
\text{ defined by }\lambda\mapsto
N_{\lambda}.
\]
We say that 
$e_{\wt\pi}$ is the \emph{enriched structure on $\wt C$ given by $\wt\pi$}.
\end{subsct}

\begin{theorem}\setcounter{equation}{0}
Preserve~{\rm 6.3}. 
For each $m=1,\dots,d$, let $\lambda_m\in\Z_\Upsilon$ and let 
$N_m$ be an invertible 
$\Ocal_{\wt C}$--submodule 
of $\ol N_{\lambda_m}$ 
satisfying 
$N_m\Ocal_{\ol C}=\ol N_{\lambda_m}$. For each $\l\in\Z^d$ set 
$\tau^{(\l)}:=\sum_{i=1}^d\l_i\lambda_i$. If
\begin{equation}
(N_1^{\l_1}\cdots N_d^{\l_d})\Ocal_{E\cup F}=
N_{\tau^{(\l)},\{E,F\}}
\end{equation}
for each $\l\in\Z^d$ and all distinct $E,F\in\Upsilon$ such that 
$\tau^{(\l)}_E=\tau^{(\l)}_F$, then there exists a regular smoothing 
$\wt\pi$ of $\wt C$ whose enriched structure $e_{\wt\pi}$ 
satisfies
\[
N_m=e_{\wt\pi}(\lambda_m)\text{ for }m=1,\dots,d.
\]
\end{theorem}

\begin{proof} Let $p_1,\dots,p_\delta$ denote the 
reducible 
nodes 
of 
$\wt C$. 
Let $V/M$ be the versal formal deformation of 
$\wt C$.
Then $M=\text{Spec}(R)$, where 
\[
R:=k[[t_1,\dots,t_\delta,s_1,\dots,s_e]]
\]
for a certain integer $e$.
In addition, we may assume that there is an isomorphism of $R$--algebras,
\begin{equation}
\wh O_{V,p_h}\cong R[[u_h,v_h]]/(u_hv_h-t_h),
\end{equation}
for each $h=1,\dots,\delta$, where $u_h$ and $v_h$ are the local equations 
of the branches 
of $\wt C$ 
meeting at $p_h$; see~6.1.

For each $h=1,\dots,\delta$ choose $a_h\in k^*$. 
Consider the surjection 
$\phi\:R\to k[[t]]$, 
sending each 
$t_h$ to $a_ht$ and each $s_i$ to $t$. 
Let $B:=\text{Spec}(k[[t]])$ and 
set $S:=V\times_M B$, where $B$ is viewed 
as a closed subscheme of $M$ by means of $\phi$. Since all $a_h$ are 
non--zero, the surface $S$ is regular. So $S/B$ is a regular 
smoothing of 
$\wt C$. 
In addition, the isomorphism~(6.5.2) restricts to 
an isomorphism of $k[[t]]$--algebras,
\[
\wh O_{S,p_h}\cong k[[t,u_h,v_h]]/(u_hv_h-a_ht),
\]
for each $h=1,\dots,\delta$.

Recall the notation in 6.4. 
We need to show that we may choose $a_1,\dots,a_\delta\in k^*$ such that  
$N_{\lambda_m}=N_m$ for $m=1,\dots,d$. 
Fix $h\in\{1,\dots,\delta\}$ and $m\in\{1,\dots,d\}$.
For convenience, put
\[
p:=p_h,\quad u:=u_h,\quad v:=v_h. 
\]
Let 
$E,F\in\Upsilon$ 
such that $p\in E\cap F$. Suppose that $u$ is the local 
equation of $E$ at $p$ and $v$ that of $F$. 
The completion $\wh\Ncal_{\lambda_m,p}$ of $\Ncal_{\lambda_m}$ at $p$ is 
generated by $(1/u)^{\lambda_{m,E}}(1/v)^{\lambda_{m,F}}$. In addition,
\begin{align*}
\wh\rho_{\lambda_m,\{E\}}((1/u)^{\lambda_{m,E}}(1/v)^{\lambda_{m,F}})=&
(t/u)^{\lambda_{m,E}}(1/v)^{\lambda_{m,F}}=
(1/a_h)^{\lambda_{m,E}}(1/v)^\gamma,\\
\wh\rho_{\lambda_m,\{F\}}((1/u)^{\lambda_{m,E}}(1/v)^{\lambda_{m,F}})=&
(1/u)^{\lambda_{m,E}}(t/v)^{\lambda_{m,F}}=
(1/a_h)^{\lambda_{m,F}}(1/u)^{-\gamma},
\end{align*}
where $\gamma:=\lambda_{m,F}-\lambda_{m,E}$.
For each $b\in k^*$ let
\[
\psi_b\: (1/v)^\gamma\wh\Ocal_{E,p}\oplus u^\gamma\wh\Ocal_{F,p}
\longrightarrow k
\]
be the map defined by sending $(f(1/v)^\gamma,gu^\gamma)$ to 
$f(p)-bg(p)$. Then $\wh N_{\lambda_m,p}=\text{Ker}(\psi_{a_h^\gamma})$. 
In addition, since $N_m$ is an invertible $\Ocal_{\wt C}$--submodule 
of $\ol N_{\lambda_m}$ such that $N_m\Ocal_{\ol C}=\ol N_{\lambda_m}$, there 
is $b_m\in k^*$ such that $\wh N_{m,p}=\text{Ker}(\psi_{b_m})$. We claim 
that we may choose $a_h\in k^*$ such that $a_h^\gamma=b_m$, and 
hence $\wh N_{\lambda_m,p}=\wh N_{m,p}$. In fact,
there's clearly $a_h\in k^*$ such that $a_h^\gamma=b_m$ if $\gamma\neq 0$. 
On the other hand, if $\gamma=0$ then 
$\wh N_{m,p}=\wh\Ocal_{\wt C,p}$ by Condition (6.5.1), and hence $b=1$. So 
any choice of $a_h\in k^*$ yields $a_h^\gamma=b_m$.

Now, let $m$ vary in $\{1,\dots,d\}$ and pick $b_m\in k^*$ for each $m$ 
as in the above paragraph.
By Condition~(6.5.1), 
$b_1^{\l_1}\cdots b_d^{\l_d}=1$ for every $\l\in\Z^d$ such that 
$\tau^{(\l)}_E=\tau^{(\l)}_F$. It follows that there is $a_h\in k^*$ 
such that $b_m=a_h^{\lambda_{m,F}-\lambda_{m,E}}$ for every $m$. 
Then $\wh N_{\lambda_m,p}=\wh N_{m,p}$ for 
each
$m=1,\dots,d$. 

Finally, let $h$ vary in $\{1,\dots,\delta\}$ and pick $a_h\in k^*$ for 
each $h$ as in the above paragraph. Now, for each 
$m=1,\dots,d$ the sheaves of fractional ideals $N_m$ and $N_{\lambda_m}$ 
coincide away from the 
reducible
nodes of $\wt C$ because
\[
N_m\Ocal_{\ol C}=\ol N_{\lambda_m}=N_{\lambda_m}\Ocal_{\ol C}.
\]
As we proved that $N_{\lambda_m,p}=N_{m,p}$ for each 
reducible
node $p$ of 
$\wt C$ and each $m=1,\dots,d$, we get 
$N_{\lambda_m}=N_m$ for $m=1,\dots,d$. 
\end{proof}

\begin{corollary}\setcounter{equation}{0}
Preserve~{\rm 6.3}. Let $\lambda\in\Z_\Upsilon$ and let $\mathbb P$ 
be the partition of $\Upsilon$ 
in
maximal $\lambda$--balanced subsets. Let 
$N$ be an invertible sheaf 
of $\wt C$. 
Then there is a 
regular smoothing 
$\wt\pi\:\wt S\to B$ of $\wt C$ 
such that
\[
N\cong\left.\Ocal_{\wt S}(\textstyle\sum_{E\in\Upsilon}\lambda_EE)
\right|_{\wt C}
\]
if and only if 
$N|_{\wt C_P}\cong N_{\lambda,P}$ 
for every $P\in\mathbb P$.
\end{corollary}

\begin{proof} The ``only if'' assertion is clear; see~6.4. 
Let's prove now the ``if'' assertion. Since
$N|_{\wt C_P}\cong N_{\lambda,P}$ 
for every 
$P\in\mathbb P$, we may 
assume that $N$ is 
an invertible 
$\Ocal_{\wt C}$--submodule 
of 
$\ol N_{\lambda,\mathbb P}$ satisfying 
$N\Ocal_{\ol C_{\mathbb P}}=\ol N_{\lambda,\mathbb P}$. To finish the 
proof we 
apply Theorem~6.5 with 
$d=1$, with $\lambda_1=\lambda$ and $N_1=N$. So we check that the conditions 
for Theorem~6.5 hold. First, since 
$N\Ocal_{\ol C_{\mathbb P}}=\ol N_{\lambda,\mathbb P}$, we have 
$N\Ocal_{\ol C}=\ol N_\lambda$ as well. Second, since 
$d=1$, we need only verify Condition~(6.5.1) for $\l_1=1$. 
Now, if 
$E,F\in\Upsilon$ are such that $\lambda_E=\lambda_F$, then 
there is $P\in\mathbb P$ such that ${E,F}\in P$. Since 
$N\Ocal_{\wt C_P}=N_{\lambda,P}$, we have 
$N\Ocal_{E\cup F}=N_{\lambda,{\{E,F\}}}$ 
as well. So Condition~(6.5.1) 
is verified. By Theorem~6.5, 
there is a regular smoothing $\wt\pi\:\wt S\to B$ of $\wt C$ such that 
$N\cong\Ocal_{\wt S}(\textstyle\sum_{E\in\Upsilon}\lambda_EE)|_{\wt C}$.
\end{proof}

\begin{remark}\setcounter{equation}{0} 
Preserve 6.3 and 6.4.
Main\`o calls an enriched structure 
the canonical generators of the image of 
$e_{\wt\pi}$, and view them
abstractly, without the additional structure of sheaves of fractional 
ideals. In \cite[Prop. 3.16]{Maino} she gives an intrinsic characterization 
of the enriched structures a nodal curve can have. The proof of 
Theorem~6.5 was inspired 
by the proof given to \cite[Prop. 3.16]{Maino}. In \cite{Maino}
Main\`o 
constructs and studies the moduli of stable curves with enriched structures. 
\end{remark}

\section{Regeneration, I}

\begin{theorem}\setcounter{equation}{0}
Preserve {\rm 4.1--2}. Let $L$ and $M$ be invertible sheaves on $\wt C$.
Then the following three statements hold.
\begin{enumerate} 
\item If {\rm (4.3.2a)} holds and $g_Y>0$,
then $L$ is the restriction of the 
canonical sheaf with focus on $X$ of 
some regular smoothing of $\wt C$ if and only if
\begin{equation}
L|_X\cong L_X,\quad L|_Y\cong L_Y,\quad L|_{Z_{p,j}}\cong
\begin{cases}
\Ocal_{Z_{p,j}}&\text{if $p\in I$ or $j\neq\rho_p$},\\
\Ocal_{Z_{p,j}}(1)&\text{if $p\not\in I$ and $j=\rho_p$.}
\end{cases}
\end{equation}
\item If {\rm (4.3.2b)} holds and $g_X>0$, 
then $M$ is the restriction of the 
canonical sheaf with focus on $Y$ of 
some regular smoothing of $\wt C$ if and only if
\begin{equation}
M|_X\cong M_X,\quad M|_Y\cong M_Y,\quad M|_{Z_{p,j}}\cong
\begin{cases}
\Ocal_{Z_{p,j}}&\text{if $p\in J$ or $j\neq\sigma_p$},\\
\Ocal_{Z_{p,j}}(1)&\text{if $p\not\in J$ and $j=\sigma_p$.}
\end{cases}
\end{equation}
\item If {\rm (4.3.2a--b)} hold and $g_Xg_Y>0$,
then $L$ and $M$ are the restrictions of the 
canonical sheaves with foci on $X$ and $Y$, respectively, of 
some regular smoothing of $\wt C$ if and only if {\rm (7.1.1--2)} hold and 
$(L^{\otimes\wt\beta}\otimes M^{\otimes\wt\alpha})|_{\wt C_{I\cap J}}\cong
\wt K_{I\cap J}$,
where 
\[
\wt C_{I\cap J}:=X\cup Y\cup(\bigcup_{p\in I\cap J}Z_p)\subseteq\wt C
\]
and 
\begin{equation}
\wt K_{I\cap J}:=
\wt\w|_{\wt C_{I\cap J}}^{\otimes(\wt\alpha+\wt\beta)}(
\sum_{p\not\in I\cap J}(\wt\beta\alpha_p-\wt\alpha\beta_p)(x_p-y_p)-
\wt\alpha\sum_{p\not\in J}x_p-\wt\beta\sum_{p\not\in I}y_p).
\end{equation}
\end{enumerate}
\end{theorem}

\begin{proof} Let $\Upsilon$ be the set of irreducible components of 
$\wt C$. Let 
$\lambda,\nu\in\Z_\Upsilon$
such that
\begin{equation}
\sum_{p\in\Delta} Z_p^{(\alpha_p,\rho_p)}+\gamma Y=
\sum_{E\in\Upsilon}\lambda_EE\quad\text{and}\quad
\sum_{p\in\Delta}\wh Z_p^{(\beta_p,\sigma_p)}+\epsilon X=
\sum_{E\in\Upsilon}\nu_EE.
\end{equation}

Let's prove Statement~1. 
Consider first the ``only if'' assertion. Let 
$\wt\pi$ be a regular smoothing of $\wt C$ and 
$\Lcal:=\w_{\wt\pi}(\sum_{p\in\Delta}Z_p^{(\alpha_p,\rho_p)}+\gamma Y)$,
where $\w_{\wt\pi}$ is the dualizing sheaf of $\wt\pi$.
By Theorem~5.2, 
$\Lcal$ is 
the canonical sheaf with focus on 
$X$. 
If $L\cong\Lcal|_{\widetilde C}$, then (7.1.1) holds 
by (5.1.2--3).

Consider now the ``if'' assertion. 
Preserve 6.3. 
Let $\mathbb P$ denote the partition of 
$\Upsilon$ 
in
maximal $\lambda$--balanced subsets. 
Since $g_Y>0$, by Lemma~3.1, $\alpha_p>0$ for every $p\in I$. 
Hence $\gamma>0$ as well. 
So $\lambda_X\neq\lambda_Y$, and hence 
$\wt C_P$ 
is treelike for 
every $P\in\mathbb P$. Let $N:=L\otimes\wt\w^{-1}$. 
Since (7.1.1) holds and 
$\wt C_P$ is treelike, $N|_{\wt C_P}\cong N_{\lambda,P}$ 
for every 
$P\in\mathbb P$. By Corollary~6.6, there is a regular smoothing 
$\wt\pi\:\wt S\to B$ 
of $\wt C$ such that $N\cong\Ncal|_{\wt C}$, where 
$\Ncal:=\Ocal_{\wt S}(\sum_{E\in\Upsilon}\lambda_EE)$. Let 
$\w_{\wt\pi}$ be the 
dualizing sheaf of $\wt\pi$. Since (7.1.4) 
holds, by Theorem~5.2 the canonical 
sheaf $\Lcal$ of $\wt\pi$ with focus on $X$ satisfies 
$\Lcal\cong\w_{\wt\pi}\otimes\Ncal$. So $L\cong\Lcal|_{\wt C}$. The proof 
of Statement~1 is complete. Statement~2 follows by analogy.

Let's prove Statement~3. For each pair 
$(\l,m)\in\Z\times\Z$
let 
$\tau^{(\l,m)}:=\l\lambda+m\nu$. 
Then $\tau^{(\l,m)}_X=\tau^{(\l,m)}_Y$ if and only if 
$(\l,m)\in\Z(\wt\beta,\wt\alpha)$. In addition, using that
\[
\gamma=\mu_p(\alpha_p+1)-\rho_p\text{ and }
\epsilon=\mu_p\beta_p+\sigma_p\text{ for every }p\in\Delta,
\]
it follows from (7.1.4) that
\begin{equation}
\sum_{E\in\Upsilon}\tau^{(\wt\beta,\wt\alpha)}_EE=
\frac{\gamma\epsilon}{\gcd(\gamma,\epsilon)}(X+Y+\sum_{p\in I\cap J}
\sum_{j=1}^{\mu_p-1}Z_{p,j})+D_1+D_2+D_3,
\end{equation}
where
\begin{align*}
D_1=&\sum_{p\in J\setminus I}
(\wt\alpha\epsilon+\wt\beta\alpha_p-\wt\alpha\beta_p)Z_{p,1}+
\sum_{p\in I\setminus J}
(\wt\beta\gamma+\wt\alpha\beta_p-\wt\beta\alpha_p)Z_{p,\mu_p-1},\\
D_2=&\sum_{p\not\in J}
(\wt\alpha\epsilon+\wt\beta\alpha_p-\wt\alpha(\beta_p+1))Z_{p,1}+
\sum_{p\not\in I}
(\wt\beta\gamma+\wt\alpha\beta_p-\wt\beta(\alpha_p+1))Z_{p,\mu_p-1},
\end{align*}
and where $D_3$ is a formal sum with integer coefficients 
of the rational curves 
$Z_{p,j}$ with $p\not\in I\cap J$ and $1<j<\mu_p-1$.

Consider first the ``only if'' assertion of 
Statement 3. 
Let 
$\wt\pi$ be a regular smoothing of $\wt C$, and $\w_{\wt\pi}$ its 
dualizing sheaf. By Theorem~5.2, the canonical 
sheaves $\Lcal$ and $\Mcal$ of $\wt\pi$ 
with foci on $X$ and $Y$, respectively, 
are given by (5.2.1--2). 
By (7.1.4),
\[
\Lcal^{\otimes\wt\beta}\otimes\Mcal^{\otimes\wt\alpha}\cong\w_{\wt\pi}
^{\otimes(\wt\alpha+\wt\beta)}
(\sum_{E\in\Upsilon}\tau^{(\wt\beta,\wt\alpha)}_EE).
\]
So, it follows from (7.1.5) that 
$\Lcal^{\otimes\wt\beta}\otimes\Mcal^{\otimes\wt\alpha}|_{\wt C_{I\cap J}}
\cong\wt K_{I\cap J}$.

Consider now the ``if'' 
assertion. Let $\wt L:=L\otimes{\wt\w}^{-1}$ and 
$\wt M:=M\otimes{\wt\w}^{-1}$. 
We shall use the set--up of~6.3. 
Let $\ol L:=\ol N_\lambda$ and $\ol M:=\ol N_\nu$. 
By (7.1.1--2),
we may view $\wt L$ and $\wt M$ as 
invertible $\Ocal_{\wt C}$--submodules of 
$\ol L$ and $\ol M$ 
such that 
$\wt L\Ocal_{\ol C}=\ol L$ and $\wt M\Ocal_{\ol C}=\ol M$. 
Multiplication by 
$c\in k^*_{\Upsilon}$ 
gives automorphisms of $\ol L$ and $\ol M$; let 
$\wt L_c$ and $\wt M_c$ denote the images of $\wt L$ and $\wt M$ under 
these automorphisms. 
It's clear that $\wt L_c\Ocal_{\ol C}=\ol L$ and 
$\wt M_c\Ocal_{\ol C}=\ol M$. 
For each $c\in k^*_\Upsilon$, each $p\in\Delta$ and each $j=0,\dots,\mu_p$ 
let $c_{p,j}:=c_{Z_{p,j}}$.

We'll choose 
$c,d\in k^*_\Upsilon$
such that $\wt L_c$ and $\wt M_d$ satisfy 
Condition (6.5.1) of Theorem~6.5. 
In other words, we'll pick 
$c,d\in k^*_\Upsilon$ 
such that 
\begin{equation}
\wt L_c^{\l}\wt M_d^{m}\Ocal_{E\cup F}=N_{\tau^{(\l,m)},\{E,F\}}
\end{equation}
for all 
$(\l,m)\in\Z\times\Z$
and all distinct $E,F\in\Upsilon$ such that 
$\tau^{(\l,m)}_E=\tau^{(\l,m)}_F$. 

In fact, we need only obtain (7.1.6) for $E$ and $F$ that intersect, that 
is, 
for $E$ and $F$ such that 
$\{E,F\}=\{Z_{p,j-1},Z_{p,j}\}$ 
for some $p\in\Delta$ and some $j=1,\dots,\mu_p$.
In addition, we need only obtain (7.1.6) for $(\l,m)$ in 
a $\Z$--basis of
\[
H_{E,F}:=\{(\l,m)\in
\Z\times\Z\,|\, \tau^{(\l,m)}_E=\tau^{(\l,m)}_F\}.
\]
Let $H_{p,j}:=H_{Z_{p,j-1},Z_{p,j}}$ for each $p\in\Delta$ and each 
$j=1,\dots,\mu_p$. Choose a $\Z$--basis $B_{p,j}$ of $H_{p,j}$ 
for each $p\in\Delta$ and each $j=1,\dots,\mu_p$. 
Note that 
$H_{p,j}=\Z\times\Z$ 
if and only if 
$\alpha_p=\beta_p=0$ and 
$\sigma_p<j\leq\rho_p$; 
in this case let 
$B_{p,j}:=\{(1,0),(0,1)\}$. In all other cases, 
$H_{p,j}$ has rank 1. If $p\in I\cap J$ then $(\wt\beta,\wt\alpha)$ 
is a basis of $H_{p,j}$ for every $j$; 
in this case let $B_{p,j}:=\{(\wt\beta,\wt\alpha)\}$. 

Consider distinct $E,F\in\Upsilon$ 
with non--empty intersection.
Then either 
$E\cup F$ is treelike or 
$\{E,F\}=\{X,Y\}$ and there are distinct $p,q\in\Delta$ such that 
$\mu_p=\mu_q=1$. In either case there is an isomorphism 
between 
$\wt L^{\l}\wt M^{m}\Ocal_{E\cup F}$ and $N_{\tau^{(\l,m)},\{E,F\}}$ 
for each $(\l,m)\in H_{E,F}$. In fact, 
such an isomorphism exists in the former case because both 
sheaves have the same restriction to $E$ and to $F$; and in the 
latter case because $H_{X,Y}=\Z(\wt\beta,\wt\alpha)$ and 
$L^{\otimes\wt\beta}\otimes M^{\otimes\wt\alpha}|_{X\cup Y}\cong
\wt K_{I\cap J}|_{X\cup Y}$.
So, for 
each
$(\l,m)\in H_{E,F}$ there is a unique $e^{(\l,m)}_{E,F}\in k^*$ 
such that 
\[
\wt L^{\l}\wt M^{m}\Ocal_{E\cup F}=(1,e^{(\l,m)}_{E,F})
N_{\tau^{(\l,m)},\{E,F\}}
\]
as subsheaves of 
$N_{\tau^{(\l,m)},E}\oplus N_{\tau^{(\l,m)},F}$. 
Let $e^{(\l,m)}_{p,j}:=e^{(\l,m)}_{Z_{p,j-1},Z_{p,j}}$ for each 
$p\in\Delta$ and each $j=1,\dots,\mu_p$.
Since 
$L^{\otimes\wt\beta}\otimes M^{\otimes\wt\alpha}|_{\wt C_{I\cap J}}
\cong\wt K_{I\cap J}$, 
it follows that 
\begin{equation}
e^{(\wt\beta,\wt\alpha)}_{p,1}e^{(\wt\beta,\wt\alpha)}_{p,2}\cdots 
e^{(\wt\beta,\wt\alpha)}_{p,\mu_p}=e^{(\wt\beta,\wt\alpha)}_{q,1}
e^{(\wt\beta,\wt\alpha)}_{q,2}\cdots e^{(\wt\beta,\wt\alpha)}_{q,\mu_q}
\text{ for all }p,q\in I\cap J.
\end{equation}

Therefore, finding 
$c,d\in k^*_\Upsilon$ 
such that (7.1.6) holds for all 
$(\l,m)\in\Z\times\Z$
and all distinct $E,F\in\Upsilon$ with 
$\tau^{(\l,m)}_E=\tau^{(\l,m)}_F$ is equivalent to finding 
$c,d\in k^*_\Upsilon$
such that, for each $p\in\Delta$,
\begin{equation}
c_{p,j-1}^\l d_{p,j-1}^m=e_{p,j}^{(\l,m)}
c_{p,j}^\l d_{p,j}^m
\text{ for each }j=1,\dots,\mu_p\text{ and each }(\l,m)\in B_{p,j}.
\end{equation}

Set $c_X:=1$ and $d_X:=1$. For 
$p\in I\cap J$ choose 
$c_{p,j},d_{p,j}\in k^*$ 
inductively for $j=1,\dots,\mu_p-1$ such that
\[
e_{p,j}^{(\wt\beta,\wt\alpha)}c_{p,j}^{\wt\beta}d_{p,j}^{\wt\alpha}=
c_{p,j-1}^{\wt\beta}d_{p,j-1}^{\wt\alpha}.
\]
By 
(7.1.7), 
we may choose $c_Y,d_Y\in k^*$ such that
\[
e_{p,\mu_p}^{(\wt\beta,\wt\alpha)}c_Y^{\wt\beta}d_Y^{\wt\alpha}=
c_{p,\mu_p-1}^{\wt\beta}d_{p,\mu_p-1}^{\wt\alpha}
\]
for every $p\in I\cap J$. If $I\cap J=\emptyset$, set $c_Y:=1$ and $d_Y:=1$. 
Then (7.1.8) is achieved for every $p\in I\cap J$.

If $p\in\Delta-(I\cap J)$ is such that $\alpha_p=\beta_p=0$ and 
$\sigma_p<\rho_p$, then 
choose $c_{p,j}$ inductively for $j=1,\dots,\rho_p$ such 
that
\[
e_{p,j}^{(1,0)}c_{p,j}=
c_{p,j-1},
\]
and choose $d_{p,j}$ inductively 
for $j=\mu_p-1,\dots,\sigma_p$ such that
\[
d_{p,j}=e_{p,j+1}^{(0,1)}d_{p,j+1}.
\]
In addition, set $c_{p,j}:=1$ for $j=\rho_p+1,\dots,\mu_p-1$ and 
$d_{p,j}:=1$ for 
$j=1,\dots,\sigma_p-1$.
Then (7.1.8) is achieved.

Let $p\in\Delta-(I\cap J)$ such that either $\alpha_p\neq 0$ or 
$\beta_p\neq 0$ or $\sigma_p\geq\rho_p$. Then $H_{p,j}$ has 
rank~1 
for 
every $j=1,\dots,\mu_p$. By (7.1.5), 
$\tau^{(\wt\beta,\wt\alpha)}_X=\tau^{(\wt\beta,\wt\alpha)}_Y$. 
Now, since $\wt\alpha$ and $\wt\beta$ are 
positive and $p\not\in I\cap J$,
it follows from (7.1.5) as well that either 
$\tau^{(\wt\beta,\wt\alpha)}_{Z_{p,1}}\neq\tau^{(\wt\beta,\wt\alpha)}_X$ or 
$\tau^{(\wt\beta,\wt\alpha)}_{Z_{p,\mu_p-1}}\neq
\tau^{(\wt\beta,\wt\alpha)}_Y$. So, there is 
$h_p\in\{1,\dots,\mu_p-1\}$ such that
\begin{equation}
H_{p,h_p}\cap H_{p,h_p+1}=0.
\end{equation}
For $j=1,\dots,h_p-1$ choose 
$c_{p,j},d_{p,j}\in k^*$ inductively such that
\[
e^{(\l,m)}_{p,j}c_{p,j}^\l d_{p,j}^m=c_{p,j-1}^\l d_{p,j-1}^m\text{\ where\ }
(\l,m)\in B_{p,j}.
\]
For $j=\mu_p-1,\dots,h_p+1$ choose $c_{p,j},d_{p,j}\in k^*$ inductively 
such that
\[
c_{p,j}^\l d_{p,j}^m=e^{(\l,m)}_{p,j+1}c_{p,j+1}^\l d_{p,j+1}^m
\text{\ where\ }(\l,m)\in B_{p,j+1}.
\]
Now, let $(\l_-,m_-)\in B_{p,h_p}$ and $(\l_+,m_+)\in B_{p,h_p+1}$. 
By (7.1.9), $(\l_-,m_-)$ and $(\l_+,m_+)$ are $\Z$--independent. So, 
we may choose $c_{p,h_p},d_{p,h_p}\in k^*$ such that
\[
e^{(\l_-,m_-)}_{p,h_p}c_{p,h_p}^{\l_-}d_{p,h_p}^{m_-}=
c_{p,h_p-1}^{\l_-}d_{p,h_p-1}^{m_-}\quad\text{and}\quad 
c_{p,h_p}^{\l_+}d_{p,h_p}^{m_+}=e^{(\l_+,m_+)}_{p,h_p+1}
c_{p,h_p+1}^{\l_+}d_{p,h_p+1}^{m_+}.
\]
Then (7.1.8) is achieved.

So, we obtained 
$c,d\in k^*_\Upsilon$ 
such that (7.1.6) holds for all 
$(\l,m)\in\Z\times\Z$
and all distinct $E,F\in\Upsilon$ such that 
$\tau^{(\l,m)}_E=\tau^{(\l,m)}_F$. By Theorem~6.5 there is a regular 
smoothing 
$\wt\pi\:\wt S\to B$ 
of $\wt C$ such that 
\[
\wt L_c\cong\Ocal_{\wt S}(\sum_{E\in\Upsilon}\lambda_EE)|_{\wt C}
\quad\text{and}\quad
\wt M_d\cong\Ocal_{\wt S}(\sum_{E\in\Upsilon}\nu_EE)|_{\wt C}.
\] 
Since $\wt L_c\cong L\otimes\wt\w^{-1}$ and 
$\wt M_d\cong M\otimes\wt\w^{-1}$, it follows from (7.1.4) and 
Theorem~5.2 that $L$ and $M$ are the restrictions of the 
canonical sheaves of $\wt\pi$ with foci on $X$ and $Y$, respectively.
\end{proof}

\section{Regeneration, II}

\begin{definition}\setcounter{equation}{0}
Preserve {\rm 4.1--2}. 
Let $V\subseteq H^0(L_X)$ and $W\subseteq H^0(M_Y)$ be 
vector subspaces. 

Suppose (4.3.2a) holds. Then $V$ is called \emph{$\mu$--smoothable} 
if $(V,L_X)$ is the 
limit canonical aspect with focus on $X$ of a regular smoothing of 
$\wt C$. Let $\V_{\mu,X}\subseteq\G_X$ 
be the subset of $\mu$--smoothable subspaces. 

Suppose (4.3.2b) holds. Then $W$ is called \emph{$\mu$--smoothable} 
if $(W,M_Y)$ is the 
limit canonical aspect with focus on $Y$ of a regular smoothing of 
$\wt C$. Let $\V_{\mu,Y}\subseteq\G_Y$ 
be the subset of $\mu$--smoothable subspaces.

Suppose 
(4.3.2a--b) 
hold. Then the pair 
$(V,W)$ is called \emph{$\mu$--smoothable} if 
$(V,L_X)$ and $(W,M_Y)$ are the limit canonical 
aspects with foci on $X$ and $Y$ of 
the same regular smoothing of $\wt C$. Let 
$\V_\mu\subseteq\V_{\mu,X}\times\V_{\mu,Y}$ 
be the subset of $\mu$--smoothable pairs. 
\end{definition}

\begin{theorem}\setcounter{equation}{0}
Preserve {\rm 4.1--2}. For each $F\subseteq\Delta$ 
let $C_F$ be the blow--up of $C$ along $\Delta-F$. 
Let $V\subseteq H^0(L_X)$ and 
$W\subseteq H^0(M_Y)$ be vector subspaces.
Then the following three statements hold.
\begin{enumerate}
\item If {\rm (4.3.2a)} holds and $g_Y>0$, 
then $V\in\V_{\mu,X}$ if and only if 
there is an invertible sheaf $L$ on $C_I$ such that 
$L|_X\cong L_X$ and $L|_Y\cong L_Y$ and such that 
$V$ is the image of the restriction map 
$H^0(L)\to H^0(L_X)$. 
\item If {\rm (4.3.2b)} holds and $g_X>0$, then 
$W\in\V_{\mu,Y}$ if and only if 
there is an invertible sheaf $M$ on $C_J$ such that 
$M|_X\cong M_X$ and $M|_Y\cong M_Y$ and such that 
$W$ is the image of the restriction map 
$H^0(M)\to H^0(M_Y)$.
\item If {\rm (4.3.2a--b)} hold and $g_Xg_Y>0$, then 
$(V,W)\in\V_\mu$ 
if and only if there are an invertible sheaf $L$ on $C_I$ 
as in Statement 1 and an 
invertible sheaf $M$ on $C_J$ 
as in Statement 2 such that 
their pull--backs 
$L_{I\cap J}$ and $M_{I\cap J}$ to $C_{I\cap J}$ satisfy
\begin{equation}
L_{I\cap J}^{\otimes\wt\beta}\otimes M_{I\cap J}^{\otimes\wt\alpha}\cong
K_{I\cap J},
\end{equation}
where 
\begin{equation}
K_{I\cap J}:=\w_{I\cap J}^{\otimes(\wt\alpha+\wt\beta)}
(\sum_{p\not\in I\cap J}
(\wt\beta\alpha_p-\wt\alpha\beta_p)(x_p-y_p)-\wt\alpha\sum_{p\not\in J}x_p-
\wt\beta\sum_{p\not\in I}y_p),
\end{equation}
where $\w_{I\cap J}$ is the pull--back to $C_{I\cap J}$ 
of the dualizing sheaf $\w$ of $C$.
\end{enumerate}
\end{theorem}

\begin{proof} Let $\psi\:\wt C\to C$ be 
the map contracting the subcurve $Z_p$ for every $p\in\Delta$. 
For each $F\subseteq\Delta$ let $\varphi_F\:C_F\to C$ be the blow--up map, 
\[
\wt C_F:=X\cup Y\cup(\textstyle\bigcup_{p\in F}Z_p)\subseteq\wt C,
\]
and let $\psi_F\:\wt C_F\to C_F$ be the map contracting 
the subcurve $Z_p$ for every $p\in F$.

Let's prove Statement~1. First, let $\wt\pi$ be a regular 
smoothing of $\wt C$ such that $(V,L_X)$ is the limit canonical 
aspect of $\wt\pi$ with focus on $X$. 
Let $\wt L$ be the restriction to $\wt C$ of the canonical 
sheaf of $\wt\pi$ with focus on $X$. By Proposition~7.1,
\begin{subequations}
\begin{align}
&\wt L|_X\cong L_X,\quad \wt L|_Y\cong L_Y,\\
&\wt L|_{Z_{p,j}}\cong
\begin{cases}
\Ocal_{Z_{p,j}}&\text{if $p\in I$ or $j\neq\rho_p$},\\
\Ocal_{Z_{p,j}}(1)&\text{if $p\not\in I$ and $j=\rho_p$.}
\end{cases}
\end{align}
\end{subequations}
By (8.2.3b), the sheaf $\psi_*\wt L$ is torsion--free, rank--1, and fails 
to be invertible precisely along the set $\Delta-I$. 
Hence there are an invertible 
sheaf $L$ on $C_I$ and an isomorphism 
$\lambda\:\psi_*\wt L\to\varphi_{I*}L$. 
Restricting $\lambda$ to $X$ and $Y$, and removing torsion, we obtain 
isomorphisms 
$\wt L|_X\to L|_X$ and $\wt L|_Y\to L|_Y$. 
So there is a commutative diagram,
\begin{equation}
\begin{CD}
H^0(\wt L) @>>> H^0(L)\\
@VVV @VVV\\
H^0(\wt L|_X) @>>> H^0(L|_X),
\end{CD}
\end{equation}
where the vertical maps are restriction maps, and the horizontal maps 
are the isomorphisms induced by $\lambda$. 
In addition, $L|_X\cong L_X$ and $L|_Y\cong L_Y$ by (8.2.3a).
Now, 
since (8.2.4) is commutative, the restriction maps $\tau\:H^0(L)\to H^0(L_X)$ and 
$\wt\tau\:H^0(\wt L)\to H^0(L_X)$ have the same image. 
Hence 
$V=\Image(\wt\tau)$ 
by Theorem~5.2, 
because $(V,L_X)$ is the limit canonical 
aspect of $\wt\pi$ with focus on 
$X$. So $V=\Image(\tau)$.

Conversely, let $L$ be an invertible sheaf on $C_I$ as in Statement~1. 
Then $\varphi_{I*}L$ is torsion--free, rank--1, and fails to be 
invertible precisely along $\Delta-I$. 
As
$\mu_p>1$ for 
each 
$p\in\Delta-I$, there are an invertible sheaf $\wt L$ on $\wt C$ such 
that (8.2.3b) holds and an isomorphism 
$\lambda\:\psi_*\wt L\to\varphi_{I*}L$. As before, $\lambda$ induces 
isomorphisms $\wt L|_X\cong L|_X$ and $\wt L|_Y\cong L|_Y$. So 
(8.2.3a) holds. By Theorem~7.1, there is 
a regular smoothing of $\wt C$ whose canonical sheaf with focus on $X$ 
restricts to $\wt L$ on $\wt C$. As before, 
$V$ is the image of the restriction map $H^0(\wt L)\to H^0(L_X)$. 
So $V$ is $\mu$--smoothable. The proof of Statement~1 is complete. 
Statement~2 follows by analogy.

Let's prove Statement~3. First, let $\wt\pi$ be a regular 
smoothing of $\wt C$ such that $(V,L_X)$ and $(W,M_Y)$ are the limit 
canonical aspects of $\wt\pi$ with foci on $X$ and $Y$, respectively. 
Let $\wt L$ and $\wt M$ be the restrictions to 
$\wt C$ of the canonical sheaves of $\wt\pi$ 
with foci on $X$ and $Y$, respectively. 
As observed in the proof of Statement~1, 
there are invertible sheaves $L$ and $M$ on 
$C_I$ and $C_J$, respectively, such that 
$\psi_*\wt L\cong\varphi_{I*}L$ and 
$\psi_*\wt M\cong\varphi_{J*}M$. In addition, $L$ and $M$ satisfy 
Statements~1~and~2, respectively.

Since $\psi_*\wt L\cong\varphi_{I*}L$, there is 
a map $\psi_I^*L\to\wt L|_{\wt C_I}$. This map is 
injective because it is injective generically on $X$ and $Y$ and 
because $L$ is invertible. So 
$\psi_I^*L$ and $\wt L|_{\wt C_I}$ are isomorphic 
because 
they 
restrict to 
isomorphic sheaves on the irreducible components of $\wt C_I$. 
Analogously, 
$\psi_J^*M\cong\wt M|_{\wt C_J}$. 
Let $L_{I\cap J}$ and $M_{I\cap J}$ be the 
pull--backs of $L$ and $M$ to $C_{I\cap J}$. 
Let $\wt K_{I\cap J}$ be given by (7.1.3). 
Since 
$(\wt L^{\otimes\wt\beta}\otimes\wt M^{\otimes\wt\alpha})|_{\wt C_{I\cap J}}
\cong\wt K_{I\cap J}$
by Theorem~7.1, 
we have
\begin{equation}
\psi_{I\cap J}^*(L^{\otimes\wt\beta}_{I\cap J}\otimes 
M_{I\cap J}^{\otimes\wt\alpha})\cong
\wt K_{I\cap J}.
\end{equation}
Now, since $\wt\w\cong\psi^*\w$, 
it follows that 
$\wt K_{I\cap J}\cong\psi_{I\cap J}^* K_{I\cap J}$. 
So (8.2.1) holds.

Conversely, let $L$ be an invertible sheaf on $C_I$ and $M$ an invertible 
sheaf on $C_J$ as in Statement~3. As observed in the proof of 
Statement~1, there are invertible sheaves $\wt L$ and $\wt M$ on $\wt C$ 
such that $\psi_*\wt L\cong\varphi_{I*}L$ and 
$\psi_*\wt M\cong\varphi_{J*}M$, and such that (8.2.3) and (8.2.6) below 
hold.
\begin{subequations}
\begin{align}
&\wt M|_X\cong M_X,\quad \wt M|_Y\cong M_Y,\\
&\wt M|_{Z_{p,j}}\cong
\begin{cases}
\Ocal_{Z_{p,j}}&\text{if $p\in J$ or $j\neq\sigma_p$},\\
\Ocal_{Z_{p,j}}(1)&\text{if $p\not\in J$ and $j=\sigma_p$.}
\end{cases}
\end{align}
\end{subequations}
Now, 
since $\psi_{I\cap J}^*K_{I\cap J}\cong\wt K_{I\cap J}$, it follows 
from (8.2.1) that (8.2.5) holds.
As shown above, 
$\psi_I^*L\cong\wt L|_{\wt C_I}$ and $\psi_J^*M\cong\wt M|_{\wt C_J}$. 
Hence 
$(\wt L^{\otimes\wt\beta}\otimes\wt M^{\otimes\wt\alpha})|_{\wt C_{I\cap J}}
\cong\wt K_{I\cap J}$.
By Theorem~7.1, there is a regular smoothing of $\wt C$ whose canonical 
sheaves with foci on $X$ and $Y$ restrict to $\wt L$ and $\wt M$, 
respectively. As shown in the proof of Statement~1, 
the restriction maps $H^0(\wt L)\to H^0(L_X)$ and 
$H^0(\wt M)\to H^0(M_Y)$ have images $V$ and $W$, respectively. 
So $(V,W)$ is $\mu$--smoothable.
\end{proof}

\begin{subsct}\setcounter{equation}{0}
\emph{Tori actions.} 
Preserve {\rm 4.1--2}. 
Denote by $C_{I\cap J}$ the blow--up of $C$ along the points of
$\Delta-(I\cap J)$. Fix isomorphisms $\zeta_{X,p}\:L_X(x_p)\to k$ and 
$\zeta_{Y,p}\:L_Y(y_p)\to k$ for each $p\in I$, and 
$\xi_{X,p}\:M_X(x_p)\to k$ and $\xi_{Y,p}\:M_Y(y_p)\to k$ for each 
$p\in J$. 
If $g_Xg_Y>0$ and $I\cap J\neq\emptyset$, choose them such that 
\[
\{(\zeta_{X,p}^{\otimes\wt\beta}\otimes\xi_{X,p}^{\otimes\wt\alpha},
\zeta_{Y,p}^{\otimes\wt\beta}\otimes\xi_{Y,p}^{\otimes\wt\alpha})\,|\, 
p\in I\cap J\}
\]
patch $L_X^{\otimes\wt\beta}\otimes M_X^{\otimes\wt\alpha}$ and 
$L_Y^{\otimes\wt\beta}\otimes M_Y^{\otimes\wt\alpha}$ 
to the sheaf $K_{I\cap J}$ given by~(8.2.2). 
Consider the corresponding evaluation maps,
\[
e_X\:H^0(L_X)\to k_I,\  e_Y\:H^0(L_Y)\to k_I,\  f_X\:H^0(M_X)\to k_J,\  
f_Y\:H^0(M_Y)\to k_J.
\]
Let $V:=\Image(e_Y)$ and $W:=\Image(f_X)$. 
Let $h_X:=\dim V$ and $h_Y:=\dim W$. 
Let 
\[
G_X:=\Grass_{h_X}(k_I)\quad\text{and}\quad G_Y:=\Grass_{h_Y}(k_J).
\]
Consider the natural actions of the tori 
$k^*_I$ and $k^*_J$ on $k_I$ and $k_J$, 
and their respective actions on $G_X$ and $G_Y$. Denote by $\O_V$ 
and $\O_W$ the orbits of $V$ and $W$ under these actions, 
and by $\psi_X\:k^*_I\to\O_V$ and $\psi_Y\:k^*_J\to\O_W$ the 
orbit maps. 
If $g_Xg_Y>0$, let
\begin{equation}
T:=\{(s,t)\in 
k^*_I\times k^*_J
\,|\, s_p^{\wt\beta}=t_p^{\wt\alpha}
\text{ for every }p\in I\cap J\},
\end{equation}
and denote by $\O$ the orbit of $(V,W)$ under the 
induced action of $T$ on $G_X \times G_Y$.
\end{subsct}

\begin{lemma}\setcounter{equation}{0} 
Preserve {\rm 4.1--2} and {\rm 8.3}.
Then the following three statements hold.
\begin{enumerate}
\item If {\rm (4.3.2a)} holds and $\alpha\neq 0$, 
then $e_X$ induces a closed 
embedding $\iota_X\:G_X \to \G_X$ such that 
$\iota_X(\O_V) = \V_{\mu,X}$.
\item If {\rm (4.3.2b)} holds and $\beta\neq 0$, then $f_Y$ induces a closed
embedding $\iota_Y\:G_Y \to \G_Y$ such that
$\iota_Y(\O_W) = \V_{\mu,Y}$.
\item If {\rm (4.3.2a--b)} hold, $\alpha \neq 0$ and $\beta \neq 0$, then 
$\V_\mu = (\iota_X \times \iota_Y)(\O)$, where $\iota_X$ and
$\iota_Y$ are the embeddings 
mentioned in Statements~1~and~2.
\end{enumerate}
\end{lemma}
\begin{proof} 
Let's prove Statement~1. 
Since $\alpha\neq 0$, also $g_Y>0$, and $\alpha_p>0$ for every 
$p\in I$ by Lemma~3.1. By Riemann--Roch, 
$h^0(L_X)=g+|\alpha|-g_Y$ and 
$e_X$ is surjective. Now, 
$e_Y$ is injective by (4.3.2a). So $h^0(L_X)-g=|I|-h_X$ 
by (5.1.5). Thus, 
taking inverse images by $e_X$ gives us a closed embedding 
$\iota_X\:G_X\to\G_X$. 
Since $g_Y>0$ and (4.3.2a) holds, the description of $\V_{\mu,X}$ given 
in Theorem~8.2 applies. Comparing this description with that of 
$\O_V$ given in~8.3, we get $\iota_X(\O_V)=\V_{\mu,X}$.
The proof of Statement~1 is complete. Statement~2 is proved analogously.
  
Let's prove now Statement~3. 
Since $\alpha$ and $\beta$ are non--zero, also $g_X$ and $g_Y$ are 
non--zero. 
Since (4.3.2a--b) hold as well, the description of $\V_\mu$ given 
in Theorem~8.2 applies. Comparing this description with that of 
$\O$ given in~8.3, we get 
$\V_\mu=(\iota_X\times\iota_Y)(\O)$. 
Statement~3 is proved.
\end{proof}

\begin{theorem}\setcounter{equation}{0}
Preserve {\rm 4.1--2}.
Then the following three statements hold.
\begin{enumerate}
\item If {\rm (4.3.2a)} holds, then $\V_{\mu,X}$
is locally closed in $\G_X$ and isomorphic 
to a torus. In addition, $\dim\V_{\mu,X}=|I|-1$ 
if $|\alpha|>g_Y$; otherwise, $\V_{\mu,X}=\{H^0(L_X)\}$.
\item If {\rm (4.3.2b)} holds, then $\V_{\mu,Y}$ 
is locally closed in $\G_Y$ and isomorphic 
to a torus. In addition, $\dim\V_{\mu,Y}=|J|-1$ 
if $|\beta|>g_X$; otherwise, $\V_{\mu,Y}=\{H^0(M_Y)\}$.
\item If {\rm (4.3.2a--b)} hold, then $\V_\mu$ is closed in 
$\V_{\mu,X}\times\V_{\mu,Y}$ and 
isomorphic to a torus. In addition,
\[
\dim\V_\mu=\begin{cases}
0&
\text{if $|\alpha|=g_Y$ and $|\beta|=g_X$,}\\
|I|-1&
\text{if $|\alpha|>g_Y$ and $|\beta|=g_X$,}\\
|J|-1&
\text{if $|\alpha|=g_Y$ and $|\beta|>g_X$,}\\
|I\cup J|-2&
\text{if $|\alpha|>g_Y$ and $|\beta|>g_X$ and $I\cap J=\emptyset$,}\\
|I\cup J|-1&
\text{if $|\alpha|>g_Y$ and $|\beta|>g_X$ and $I\cap J\neq\emptyset$.}
\end{cases}
\]
\end{enumerate}
\end{theorem}

\begin{proof} 
We shall use the set--up of 8.3. Let's prove Statement~1. 
Assume first that $|\alpha|=g_Y$. Then Theorem~5.2 says that 
$(H^0(L_X),L_X)$ is the limit canonical aspect with focus on $X$ of any 
regular smoothing of $\wt C$. Therefore $\V_{\mu,X}=\{H^0(L_X)\}$.

  Assume now that $|\alpha|>g_Y$. Then $\alpha\neq 0$ and $g_Y>0$. 
By (4.3.2a), all the Pl\"ucker coordinates of 
$V$ in $G_X$ are non--zero. In addition, 
$V\neq k^*_I$ because 
$|\alpha|>g_Y$. Since the Pl\"ucker coordinates of $V$ are non--zero, 
$\psi_X(t)=V$ if and only if $t^b=t^c$ for all subsets $b$ and $c$ of $I$ 
satisfying $|b|=|c|=h_X$. Since $h_X<|I|$, given $p,q\in I$ 
distinct, there is a subset $b\subseteq I$ with $|b|=h_X$ such that 
$p\in b$ but $q\not\in b$. Letting $c:=b-p+q$, we have $t^b=t^c$ if and 
only if $t_p=t_q$. It follows that the orbit map $\psi_X$ factors through 
an isomorphism 
$k^*_I/k^*\to\O_V$, where $k^*$ is viewed inside $k^*_I$ 
under the diagonal embedding. So $\O_V$ is isomorphic to 
a torus of dimension $|I|-1$. 

   Since $\alpha\neq 0$, Lemma~8.4 says that 
$\V_{\mu,X}$ is the image of $\O_V$ under 
a closed embedding $G_X\to\G_X$. 
Thus $\V_{\mu,X}$ is locally closed in $\G_X$ and isomorphic 
to a torus of dimension $|I|-1$. 
The proof of Statement~1 is complete. 
Statement~2 follows by analogy.

Let's prove Statement~3. If $|\beta|=g_X$
then 
$\V_\mu=\V_{\mu,X}\times\{H^0(M_Y)\}$. If $|\alpha|=g_Y$
then 
$\V_\mu=\{H^0(L_X)\}\times\V_{\mu,Y}$. So Statement~3 follows 
from Statement~1 in the first case, and from Statement~2 in the second case.

Assume now that $|\alpha|>g_Y$ and $|\beta|>g_X$. Then $\alpha$, $\beta$, 
$g_X$ and $g_Y$ are non--zero. View $k^*$ inside 
$k^*_I$ and $k^*_J$ 
under the diagonal embeddings. 
Let $T':=T\cap (k^*\times k^*)$,
where $T$ is given by (8.3.1).
Since $\wt\alpha$ and $\wt\beta$ are coprime, $T$ is a subtorus of dimension 
$|I\cup J|$ of 
$k^*_I\times k^*_J$. 
In addition, $T'$ is a one--dimensional 
subtorus of $T$ unless $I\cap J=\emptyset$, in which case 
$T'=k^*\times k^*$. 
Now, $\O\subseteq\O_V\times\O_W$ and the orbit map 
$\psi\:T\to\O$ is the restriction to $T$ of $\psi_X\times\psi_Y$. 
Since $\psi_X$ and $\psi_Y$ factor through isomorphisms 
$k^*_I/k^*\to\O_V$ and $k^*_J/k^*\to\O_W$, 
then $\psi$ factors through an isomorphism $T/T'\to\O$. So 
$\O$ is closed in $\O_V\times\O_W$ and isomorphic 
to a torus. In addition, $\dim\O=|I\cup J|-1$ unless
$I\cap J=\emptyset$, in which case $\dim\O=|I\cup J|-2$. 

Since $\alpha$ and $\beta$ are non--zero, Lemma~8.4 says that 
$\V_\mu$ is the image of $\O$ under a closed embedding 
$G_X\times G_Y\to\G_X\times\G_Y$ that sends $\O_V\times\O_W$ to 
$\V_{\mu,X}\times\V_{\mu,Y}$. Thus $\V_\mu$ is closed in 
$\V_{\mu,X}\times\V_{\mu,Y}$ and isomorphic to a torus of the dimension 
prescribed in Statement~3.
\end{proof}

\begin{subsct}\setcounter{equation}{0}
\emph{Extreme cases.} If (4.3.2a) holds, then 
\[
\V_{\mu,X}=\{V\in\G_X|\  V\supseteq H^0(L_X(-\sum_{p\in I}x_p))\text{ but } 
V\not\supseteq H^0(L_X(x_q-\sum_{p\in I}x_p))\text{ for any }q\in I\}
\]
in case $|\alpha|=g_Y+1$, and 
\[
\V_{\mu,X}=\{V\in\G_X|\  V\supseteq H^0(L_X(-\sum_{p\in I}x_p))\text{ but }
V\not\subseteq H^0(L_X(-x_q))\text{ for any }q\in I\}
\]
in case $|\alpha|=g_Y+|I|-1$. We prove here the latter case; the former 
is proved analogously.

So, suppose $|\alpha|=g_Y+|I|-1$. If $\alpha=0$ 
then $\V_{\mu,X}=\{H^0(L_X)\}$ by Theorem~8.5, which agrees with the 
above description of $\V_{\mu,X}$ in the case at hand.
Assume now that $\alpha\neq 0$. From now on we shall use the set-up of~8.3. 
By~(4.3.2a), all the Pl\"ucker coordinates of $V$ in $G_X$ 
are non--zero. Moreover, $h_X=1$. Then $\O_V$ parameterizes all 
one--dimensional subspaces $H\subseteq k_I$ 
such that
the restriction map 
$H\to k_q$ 
is an isomorphism 
for each $q \in I$. Now, Lemma~8.4 says that $\V_{\mu,X}$ 
is the image of $\O_V$ under the closed embedding $\iota_X\:G_X\to\G_X$ 
induced by $e_X$. The description of $\V_{\mu,X}$ given in the last 
paragraph follows easily now.

If (4.3.2b) holds, and either $|\beta|=g_X+1$ or $|\beta|=g_X+|J|-1$, 
then an analogous description holds for $\V_{\mu,Y}$. 
\end{subsct}

\begin{subsct}\setcounter{equation}{0}
\emph{Projections.} Preserve 4.1--2. Using Theorem~8.2 
we can describe $\V_{\mu,X}$ and $\V_{\mu,Y}$ in terms 
of central projections. In addition, we can describe 
$\V_\mu\subseteq\V_{\mu,X}\times\V_{\mu,Y}$ in terms of 
multi--linear algebra.

In fact, assume (4.3.2a) holds and $g_Y>0$. 
Then the natural map $\phi_Y\:Y\to\P(H^0(L_Y))$ 
is defined along 
$\{y_p\,|\, p\in I\}$, and the image $\{\phi_Y(y_p)\,|\, p\in I\}$ 
spans $\P(H^0(L_Y))$. 
Now, $\alpha_p>0$ for each $p\in I$ by Lemma~3.1. 
By Riemann--Roch, the natural map 
$\phi_X\:X\to\P(H^0(L_X))$ is defined 
along $\{x_p\,|\, p\in I\}$, and the 
set $\{\phi_X(x_p)\,|\, p\in I\}$ 
spans a projective subspace 
$\P_X\subseteq\mathbb P(H^0(L_X))$ of dimension $|I|-1$. Let 
$T_X$ be the set of all linear maps $\P_X\to\P(H^0(L_Y))$ 
sending $\phi_X(x_p)$ to $\phi_Y(y_p)$ for 
each
$p\in I$. Then $T_X$ is a 
torus of dimension $|I|-1$. For 
each
$t\in T_X$ let 
$\P_t\subseteq\P_X$ be the base locus of $t$, and let 
$(V_t,L_X)$ be the linear system given by projection with center $\P_t$. 
Then $V_t\in\V_{\mu,X}$ and the map 
$\lambda_X\:T_X\to\V_{\mu,X}$ sending $t$ to $V_t$ is an isomorphism.

Now, assume 
(4.3.2b) holds and $g_X>0$. The sheaves $M_X$ and $M_Y$ define natural maps
$\psi_X\:X\to\P(H^0(M_X))$ 
and 
$\psi_Y\:Y\to\P(H^0(M_Y))$. As above, 
$\psi_Y$ is defined along $\{y_p\,|\, p\in J\}$, and
the projective subspace $\P_Y\subseteq\mathbb P(H^0(M_Y))$ spanned by 
$\{\psi_Y(y_p)\,|\, p\in J\}$ 
has dimension $|J|-1$. 
In addition, $\psi_X$ is defined along $\{x_p\,|\, p\in J\}$, and 
$\{\psi_X(x_p)\,|\, p\in J\}$ spans $\P(H^0(M_X))$. 
Let $T_Y$ be the set of 
all linear maps $\P_Y\to\P(H^0(M_X))$ 
sending $\psi_Y(y_p)$ to $\psi_X(x_p)$ for 
each
$p\in J$. As above, we 
define an isomorphism $\lambda_Y\:T_Y\to\V_{\mu,Y}$.

Finally, assume 
(4.3.2a--b) hold and $g_Xg_Y>0$. Note that 
$\P_X=\P(V_X)$ and $\P_Y=\P(V_Y)$, where
\[
V_X:=\frac{H^0(L_X)}{H^0(L_X(-\sum_{p\in I}x_p))}\quad\text{and}\quad 
V_Y:=\frac{H^0(M_Y)}{H^0(M_Y(-\sum_{p\in J}y_p))}.
\]
Combining a $\wt\beta$--tuple with an $\wt\alpha$--tuple embedding, and 
composing 
with a Segre map, embed $\P_X\times\P(H^0(M_X))$ in $\ol\P_X$ and 
$\P(H^0(L_Y))\times\P_Y$ in $\ol\P_Y$, where
\begin{align*}
&\ol\P_X:=
\P(\text{Symm}_{\wt\beta}(V_X)\otimes\text{Symm}_{\wt\alpha}(H^0(M_X))),\\
&\ol\P_Y:=
\P(\text{Symm}_{\wt\beta}(H^0(L_Y))\otimes\text{Symm}_{\wt\alpha}(V_Y)).
\end{align*}
Let $C_{I\cap J}$ be the blow--up of $C$ along $\Delta-(I\cap J)$, and 
$K_{I\cap J}$ the 
invertible sheaf on $C_{I\cap J}$ given by (8.2.2). 
Let 
$\phi\:C_{I\cap J}\to\P(H^0(K_{I\cap J}))$ 
be the natural 
map, and 
$\ol\P_C$ the projective subspace of $\P(H^0(K_{I\cap J}))$ spanned by 
$\phi(I\cap J)$. In other words, $\ol\P_C=\P(\ol V)$, where
\[
\ol V:=\frac{H^0(K_{I\cap J})}{H^0(K_{I\cap J}(-\sum_{p\in I\cap J}p))}.
\]
Now, $K_{I\cap J}$ restricts to 
$L_X^{\otimes\wt\beta}\otimes M_X^{\otimes\wt\alpha}$ on $X$ and 
$L_Y^{\otimes\wt\beta}\otimes M_Y^{\otimes\wt\alpha}$ on $Y$. Thus, lifting 
sections, we obtain well--defined maps,
\[
\frac{H^0(L_X^{\otimes\wt\beta}\otimes M_X^{\otimes\wt\alpha})}
{H^0(L_X^{\otimes\wt\beta}\otimes M_X^{\otimes\wt\alpha}
(-\textstyle\sum_{p\in I}x_p))}\longrightarrow\ol V\quad\text{and}\quad
\frac{H^0(L_Y^{\otimes\wt\beta}\otimes M_Y^{\otimes\wt\alpha})}
{H^0(L_Y^{\otimes\wt\beta}\otimes M_Y^{\otimes\wt\alpha}
(-\textstyle\sum_{p\in J}y_p))}\longrightarrow\ol V,
\]
from which we derive natural linear maps 
$\tau_X\:\ol\P_C\to\ol\P_X$ and $\tau_Y\:\ol\P_C\to\ol\P_Y$, respectively.

Each $s\in T_X$ induces a linear map $\zeta_s\:\ol\P_X\to\ol\P$, 
where
\[
\ol\P:=
\P(\text{Symm}_{\wt\beta}(H^0(L_Y))\otimes\text{Symm}_{\wt\alpha}(H^0(M_X))).
\]
Also, each $t\in T_Y$ induces a linear map $\xi_t\:\ol\P_Y\to\ol\P$. 
Let
\[
T:=\{(s,t)\in T_X\times T_Y\,|\, \zeta_s\tau_X=\xi_t\tau_Y\}.
\]
Then $(\lambda_X\times\lambda_Y)(T)=\V_\mu$.
\end{subsct}

\section{Boundary of tori orbits on Grassmannians}

\begin{lemma}\setcounter{equation}{0} 
Let $I$ be a finite set
and let
$V\subseteq k_I$ be a non--zero 
vector
subspace. Let $h:=\dim V$ and 
$G:=\Grass_h(k_I)$. 
For each ordered tripartition 
$\mathbb I=(I',\ol I,I'')$ of $I$, let
\begin{equation}
V_{\mathbb I}:=k_{I'}+(k_{\ol I}\cap(V+k_{I''})).
\end{equation}
Let $\O$ denote the orbit of $V$ under the natural action of 
$k^*_I$ on $G$. If all Pl\"ucker coordinates of $V$ in $G$ are non--zero, 
then 
the
closure $\ol{\O}\subseteq G$ is the union of the 
orbits of the subspaces $V_{\mathbb I}\subseteq k_I$ obtained 
from all ordered 
tripartitions 
$\mathbb I=(I',\ol I,I'')$ of $I$ 
satisfying 
$|I'|<h\leq |I-I''|$.
\end{lemma}

\begin{proof} 
Let $\mathcal B:=\{b\subseteq I;\, |b|=h\}$. For each 
$b\in\mathcal B$ let $p_b$ denote the corresponding Pl\"ucker coordinate 
on $G$. In terms of Pl\"ucker coordinates, $s\in k^*_I$ acts on a point of 
$G$ with coordinates $(p_b\,|\, b\in\mathcal B)$ by taking it to the 
point with coordinates $(s^bp_b\,|\, b\in\mathcal B)$. 

Let $(\wt p_b\,|\, b\in\mathcal B)$ be Pl\"ucker coordinates of 
$V$ in $G$. By assumption, $\wt p_b\neq 0$ for 
every $b\in\mathcal B$. Let $Z\subseteq G$ be defined by the 
following 
equations on Pl\"ucker coordinates,
\begin{equation}
\wt p_{b_1}\wt p_{b_2}p_{b_3}p_{b_4}=\wt p_{b_3}\wt p_{b_4}p_{b_1}p_{b_2}
\quad\text{ for all }b_1,\dots,b_4\in\mathcal B\text{ with }b_1+b_2=b_3+b_4. 
\end{equation}
It's clear that 
$\O\subseteq Z$.

Let $W\in G$ with Pl\"ucker coordinates $(p_b\,|\, b\in\mathcal B)$. Let 
\[
\mathcal B_W:=\{b\in\mathcal B\,|\, p_b\neq 0\},\quad 
I'_W:=\bigcap_{b\in\mathcal B_W}b,\quad 
I''_W:=I-\bigcup_{b\in\mathcal B_W}b.
\]
Note that $|I'_W|\leq h\leq |I-I''_W|$. 
Set $\ol I_W:=I-(I'_W\cup I''_W)$.

We divide the proof in 
five 
steps.

\emph{Step~1}: Let $W\in Z$ with Pl\"ucker coordinates 
$(p_b\,|\, b\in\mathcal B)$. Then the following 
two
conditions hold.
\begin{subequations}
\begin{align}
&\mathcal B_W=\{b\in\mathcal B\,|\, I'_W\subseteq b\subseteq I-I''_W\}\\
&\text{There is }\ol s\in k^*_{\ol I_W}\text{ such that }
p_{b_1}\wt p_{b_2}=\ol s^{(b_1-b_2)}\wt p_{b_1}p_{b_2}
\text{ for all }b_1,b_2\in\mathcal B_W.
\end{align}
\end{subequations}

Indeed, let's check (9.1.3a) first.
If $b\in\mathcal B_W$
then 
$I'_W\subseteq b\subseteq I-I''_W$ by definition.

Conversely, suppose first that 
$|\mathcal B_W|=1$. Then $I'_W=I-I''_W$, 
and thus (9.1.3a) holds. Suppose now that $|\mathcal B_W|>1$. Let 
$b\in\mathcal B_W$. Then 
$I'_W\neq b\neq I-I''_W$. Let $i,j\in I$ such that 
$i\in b-I'_W$ and $j\in I-I''_W-b$. Using an induction argument, 
we need only show that 
$b-i+j\in\mathcal B_W$. Since $i\not\in I'_W$ and 
$j\not\in I''_W$, there are $b_1,b_2\in\mathcal B_W$ such that 
$i\not\in b_1$ and $j\in b_2$. We'll find 
$b_0\in\mathcal B_W$ such that $i\not\in b_0$ and $j\in b_0$. If 
$j\in b_1$ or $i\not\in b_2$, set $b_0:=b_1$ or $b_0:=b_2$, 
respectively. Suppose that $j\not\in b_1$ and $i\in b_2$. 
Then $|b_1\cap b_2|< h$, and hence there is 
$u\in b_1 - b_2$. Set 
$b_0:=b_1-u+j$. Then $i\not\in b_0$ and $j\in b_0$. 
Since $b_1+b_2=b_0+(b_2-j+u)$, 
Equations~(9.1.2) imply that $b_0\in\mathcal B_W$. 
Hence, since $b+b_0=(b-i+j)+(b_0-j+i)$, Equations~(9.1.2) imply 
that $b-i+j\in\mathcal B_{W}$ as well. Thus (9.1.3a) holds. 

Let's now check (9.1.3b). 
If $\ol I_W=\emptyset$ then $|\mathcal B_W|=1$, and hence (9.1.3b) holds 
trivially.
Assume that $\ol I_W\neq\emptyset$, and fix 
$\ol i\in\ol I_W$. Define 
$\ol s\in k^*_{\ol I_W}$ by letting $\ol s_{\ol i}:=1$ and 
\[
\ol s_i:=\frac{p_{b-\ol i+i}\wt p_b}{\wt p_{b-\ol i+i}p_b}
\]
for each $i\in\ol I_W-\ol i$, where $b\in\mathcal B_W$ is chosen 
such that $i\not\in b$ and $\ol i\in b$. 
Observe that such $b$ exists because (9.1.3a) holds. Moreover, 
$\ol s_i$ is independent of the choice of $b$, 
because if $b'\in\mathcal B_W$ is 
such that $i\not\in b'$ and $\ol i\in b'$, then 
$b+(b'-\ol i+i)=(b-\ol i+i)+b'$, and hence the claimed 
independence of $\ol s_i$ follows from Equations~(9.1.2).

If $i$ and $j$ are distinct elements of $\ol I_W$, then 
\begin{equation}
p_b\wt p_{b-j+i}=(\ol s_j/\ol s_i)\wt p_bp_{b-j+i}
\end{equation}
for every $b\in\mathcal B_W$ such that $i\not\in b$ and $j\in b$. 
Indeed, if either $\ol i=i$ or $\ol i=j$, then (9.1.4) 
holds by definition of $\ol s_j$ or $\ol s_i$. Suppose 
that $\ol i\in\ol I_W-i-j$. If $\ol i\in b$, 
then
\[
\frac{\ol s_i}{\ol s_j}=(\frac{p_{b-\ol i+i}\wt p_b}{\wt p_{b-\ol i+i}p_b})/
(\frac{p_{b-\ol i+i}\wt p_{b-j+i}}{\wt p_{b-\ol i+i}p_{b-j+i}})=
\frac{p_{b-j+i}\wt p_b}{p_b\wt p_{b-j+i}}.
\]
If $\ol i\not\in b$, then
\[
\frac{\ol s_i}{\ol s_j}=
(\frac{p_{b-j+i}\wt p_{b-j+\ol i}}{\wt p_{b-j+i}p_{b-j+\ol i}})/
(\frac{p_b\wt p_{b-j+\ol i}}{\wt p_bp_{b-j+\ol i}})=
\frac{p_{b-j+i}\wt p_b}{p_b\wt p_{b-j+i}}.
\]
Either way, (9.1.4) holds. Now, (9.1.3b) follows from (9.1.4) and an 
induction argument.
Step~1 is complete.

\emph{Step~2}: Let $W\in G$ with Pl\"ucker coordinates 
$(p_b\,|\, b\in\mathcal B)$. If (9.1.3a,b) hold, then $W\in\ol{\O}$. 
Indeed, let $\ol s\in k^*_{\ol I_W}$ as in (9.1.3b). Define 
$u\in\Z_I$ and $s\in k^*_I$ by letting for each 
$i\in I$,
\[
u_i:=
\begin{cases}
-1
&\text{if $i\in I'_W$,}\\
0
&\text{if $i\in\ol I_W$,}\\
1
&\text{if $i\in I''_W$.}
\end{cases}
\quad\text{and}\quad
s_i:=
\begin{cases}
1
& \text{if $i\in I'_W$,}\\
\ol s_i
& \text{if $i\in\ol I_W$,}\\
1
&  \text{if $i\in I''_W$.}
\end{cases}
\]
Define a map 
$\zeta\:k^*\to k^*_I$ 
by letting $\zeta(r)_i:=s_ir^{u_i}$ for 
each $i\in I$ and each $r\in k^*$. For each $r\in k^*$ let 
$W_r:=\zeta(r)V$ and $(p_b(r)\,|\, b\in\mathcal B)$ 
be Pl\"ucker coordinates of $W_r$. Because of our choice of $\zeta$, 
the power of $r$ 
in the expression for $p_b(r)$ is minimum exactly among those 
$b\in\mathcal B_W$. Let $\ol W\in G$ be 
the limit of $W_r$ as $r$ goes to 0, and 
$(\ol p_b\,|\, b\in\mathcal B)$ 
Pl\"ucker coordinates of $\ol W$. 
Of course, $\ol W\in\ol{\O}$. Now, 
$\ol p_b\neq 0$ if and only if $b\in\mathcal B_W$. Moreover, 
if 
$b_1,b_2\in\mathcal B_W$ then
\[
\frac{\ol p_{b_1}}{\ol p_{b_2}}=s^{(b_1-b_2)}\frac{\wt p_{b_1}}{\wt p_{b_2}}=
\ol s^{(b_1-b_2)}\frac{\wt p_{b_1}}{\wt p_{b_2}}=\frac{p_{b_1}}{p_{b_2}},
\]
where the last equality holds by (9.1.3b). So $\ol W=W$, 
and hence $W\in\ol{\O}$. Step 2 is complete.

\emph{Step~3}: If $\mathbb I=(I,\ol I,I'')$ is an ordered tripartition of 
$I$ such that $|I'|\leq h\leq |I-I''|$, then $\dim V_{\mathbb I}=h$ and 
$(p_b\,|\, b\in\mathcal B)$ are Pl\"ucker coordinates of $V_{\mathbb I}$, 
where
\[
p_b:=
\begin{cases}
\wt p_b&\text{if $I'\subseteq b\subseteq I-I''$},\\
0&\text{otherwise}.
\end{cases} 
\]
Indeed, let $v_1,\dots,v_h$ be a basis of $V$ in $k_I$. Write 
$I=\{1,\dots,\delta\}$ with 
\begin{align*}
I'&=\{1,\dots,|I'|\},\\  
\ol I&=\{|I'|+1,\dots,|I'+\ol I|\},\\  
I''&=\{|I'+\ol I|+1,\dots,|I|\}.
\end{align*}
Since no Pl\"ucker coordinate of $V$ is zero, we 
may assume that the matrix $(v_{i,j})_{1\leq i,j\leq h}$ is the 
identity. Since $|I'|\leq h\leq|I-I''|$, a basis of $V_{\mathbb I}$ is 
given by the vectors $\ol v_1,\dots,\ol v_h$ defined by
\[
\ol v_{i,j}:=
\begin{cases}
0&\text{if $j\in I''$},\\
0&\text{if $i\in I'$ and $j\in\ol I$},\\
v_{i,j}&\text{otherwise}.
\end{cases} 
\]
In particular, $\dim V_{\mathbb I}=h$. It follows as well that 
$(p_b\,|\, b\in\mathcal B)$ are Pl\"ucker coordinates of $V_{\mathbb I}$.
Step~3 is complete.

\emph{Step~4}: 
If $W\in\ol{\O}$, then there is an ordered 
tripartition $\mathbb I=(I',\ol I,I'')$ of $I$ such that 
$|I'|<h\leq |I-I''|$, and such that $W$ is in the 
orbit of $V_{\mathbb I}$ defined by (9.1.1). 
Indeed, let $W\in\ol{\O}$ with Pl\"ucker coordinates 
$(p_b\,|\, b\in\mathcal B)$. Since $\O \subseteq Z$, we have 
$W\in Z$, 
and hence (9.1.3a,b) hold by 
Step~1. Since (9.1.3b) holds, we may replace $W$ by a certain point in its 
orbit and assume that $p_b=\wt p_b$ for every $b\in\mathcal B_W$.
Hence, using Step~3, it follows from (9.1.3a) that $W=V_{\wt{\mathbb I}}$, 
where $\wt{\mathbb I}:=(I'_W,\ol I_W,I''_W)$. 

If $|I'_W|<h$ set $\mathbb I:=\wt{\mathbb I}$. 
Step~4 
is complete 
in this case. If $|I'_W|=h$
then $\ol I_W=\emptyset$, and thus 
$V_{\wt{\mathbb I}}=k_{I'_W}$. In addition, 
$I'_W\neq\emptyset$ because $h>0$. Let 
$\ol I\subseteq I'_W$ be any non--empty subset and put 
$\mathbb I:=(I'_W-\ol I,\ol I,I''_W)$. Now, 
$V_{\mathbb I}=k_{I'_W}$ since $h=|I-I''_W|$. 
Then $V_{\mathbb I}=V_{\wt{\mathbb I}}=W$.
Step~4 is complete.

\emph{Step~5}: 
If $\mathbb I=(I',\ol I,I'')$ is an 
ordered tripartition of $I$ such that 
$|I'|<h\leq|I-I''|$, 
then 
$V_{\mathbb I}\in\ol{\O}$. Indeed, 
let $(p_b\,|\, b\in\mathcal B)$ be Pl\"ucker coordinates of 
$V_{\mathbb I}$. 
By 
Step~3,
\[
\mathcal B_{V_{\mathbb I}}=
\{b\in\mathcal B\,|\,I'\subseteq b\subseteq I-I''\}.
\]
If $h<|I-I''|$
then $I'_{V_{\mathbb I}}=I'$ and 
$I''_{V_{\mathbb I}}=I''$. If $h=|I-I''|$
then 
$I'_{V_{\mathbb I}}=I'\cup\ol I$ and $I''_{V_{\mathbb I}}=I''$. 
In either case 
(9.1.3a) holds. 
By Step~3 as well,
(9.1.3b) holds 
trivially with $\ol s=1$. By Step~2, $V_{\mathbb I}\in\ol{\O}$, 
thus finishing 
Step~5.

Steps 4 and 5 show the lemma.
\end{proof}

\begin{lemma}\setcounter{equation}{0} 
Let $\Delta$ be a finite set, 
$I$ and $J$ non--empty subsets of $\Delta$. Let $V\subseteq k_I$ and 
$W\subseteq k_J$ be non--zero vector subspaces. Let $h_1:=\dim V$ and 
$G_1:=\Grass_{h_1}(k_I)$. 
Let $h_2:=\dim W$ 
and 
$G_2:=\Grass_{h_2}(k_J)$. 
For each ordered tripartition $\mathbb I=(I',\ol I,I'')$ of $I$ and 
each ordered tripartition $\mathbb J=(J',\ol J,J'')$ of $J$, let
\begin{equation}
V_{\mathbb I}:=k_{I'}+(k_{\ol I}\cap(V+k_{I''}))\quad\text{and}\quad
W_{\mathbb J}:=k_{J'}+(k_{\ol J}\cap(W+k_{J''})).
\end{equation}
Let $\lambda$ and $\tau$ be positive integers, and set
\[
T:=\{(s,t)\in k^*_I\times k^*_J\,|\, 
s_i^{\tau}t_j^{\lambda}=s_j^{\tau}t_i^{\lambda}
\text{ for all }i,j\in I\cap J\}.
\]
Let $\O$ denote the orbit of $(V,W)$ 
under the induced action of $T$ on $G_1\times G_2$. 
If all 
Pl\"ucker coordinates of $V$ and $W$ in $G_1$ and 
$G_2$ are non--zero, then 
the
closure 
$\ol{\O}\subseteq G_1\times G_2$ is the union of the 
orbits of 
the 
pairs $(V_{\mathbb I},W_{\mathbb J})$ obtained 
from all ordered tripartitions $\mathbb I=(I',\ol I,I'')$ of $I$ and 
$\mathbb J=(J',\ol J,J'')$ of $J$ 
satisfying 
\begin{equation}
|I'|<h_1\leq |I-I''|\quad\text{and}\quad
|J'|<h_2\leq |J-J''|,
\end{equation}
and such that the following two conditions hold.
\begin{subequations}
\begin{align}
&I'\cap J\not\subseteq J'\cap I\text{ or }
J''\cap I\not\subseteq I''\cap J\Longrightarrow
J\cap(I-I')\subseteq J''\cap I,\\
&J'\cap I\not\subseteq I'\cap J\text{ or }
I''\cap J\not\subseteq J''\cap I\Longrightarrow
I\cap (J-J')\subseteq I''\cap J.
\end{align}
\end{subequations}
\end{lemma}

\begin{proof} Let 
\[
\mathcal B:=\{b\subseteq I;\  |b|=h_1\}\quad\text{and}\quad
\mathcal C:=\{c\subseteq J;\  |c|=h_2\}.
\]
For each $b\in\mathcal B$ let $p_b$ 
denote the corresponding Pl\"ucker coordinate on $G_1$. For each 
$c\in\mathcal C$ let $q_c$ denote the corresponding 
Pl\"ucker coordinate on $G_2$. 

Let $\wt\nu:=(V,W)$. Let 
$(\wt p_b\,|\, b\in\mathcal B)$ be 
Pl\"ucker coordinates of 
$V$ in $G_1$ and $(\wt q_c\,|\, c\in\mathcal C)$ 
be Pl\"ucker 
coordinates of $W$ in $G_2$. By assumption, $\wt p_b\neq 0$ for 
every $b\in\mathcal B$ and $\wt q_c\neq 0$ for every $c\in\mathcal C$. 
Let $Z\subseteq G_1\times G_2$ be 
defined by the following equations on Pl\"ucker coordinates:
\begin{subequations}
\begin{align}
\wt p_{b_1}\wt p_{b_2}p_{b_3}p_{b_4}=\wt p_{b_3}\wt p_{b_4}p_{b_1}p_{b_2}
\quad&\text{for all }b_1,\dots,b_4\in\mathcal B\text{ with }
b_1+b_2=b_3+b_4,\\
\wt q_{c_1}\wt q_{c_2}q_{c_3}q_{c_4}=\wt q_{c_3}\wt q_{c_4}q_{c_1}q_{c_2}
\quad&\text{for all }c_1,\dots,c_4\in\mathcal C\text{ with }
c_1+c_2=c_3+c_4,\\
\wt p_{b_1}^{\tau}\wt q_{c_1}^{\lambda}p_{b_2}^{\tau}q_{c_2}^{\lambda}=
\wt p_{b_2}^{\tau}\wt q_{c_2}^{\lambda}p_{b_1}^{\tau}q_{c_1}^{\lambda}
\quad&\text{for all }b_1,b_2\in\mathcal B\text{ and }c_1,c_2\in\mathcal C
\text{ satisfying (9.2.5)},
\end{align}
\end{subequations}
where
\begin{equation}
\begin{cases}
b_1+c_1=b_2+c_2,&\\
b_1\cap(I\setminus J)=b_2\cap(I\setminus J),&\\
c_1\cap(J\setminus I)=c_2\cap(J\setminus I).
\end{cases}
\end{equation}

Given $\nu\in G_1\times G_2$ with 
Pl\"ucker coordinates $(p_b\,|\, b\in\mathcal B)$ and 
$(q_c\,|\, c\in\mathcal C)$, let 
\[
\mathcal B_\nu:=\{b\in\mathcal B\,|\, p_b\neq 0\}
\quad\text{and}\quad
\mathcal C_\nu:=\{c\in\mathcal C\,|\, q_c\neq 0\}.
\]
Put
\[
I'_\nu:=\bigcap_{b\in\mathcal B_\nu}b,\quad 
I''_\nu:=I-\bigcup_{b\in\mathcal B_\nu}b,\quad
J'_\nu:=\bigcap_{c\in\mathcal C_\nu}c,\quad   
J''_\nu:=J-\bigcup_{c\in\mathcal C_\nu}c.
\]
Note that $|I'_\nu|\leq h_1\leq |I-I''_\nu|$ and 
$|J'_\nu|\leq h_2\leq|J-J''_\nu|$. 
Set
\[
\ol I_\nu:=I-(I'_\nu\cup I''_\nu)\quad\text{and}\quad
\ol J_\nu:=J-(J'_\nu\cup J''_\nu),
\]
and let
\[
\ol T_\nu:=\{(s,t)\in k^*_{\ol I_\nu}\times k^*_{\ol J_\nu}\,|\, 
s_i^{\tau}t_j^{\lambda}=s_j^{\tau}t_i^{\lambda}
\text{ for all }i,j\in\ol I_\nu\cap\ol J_\nu\}.
\]

We divide the proof in five steps.

\emph{Step~1}: If $\nu\in\O$
then $\nu\in Z$. Indeed, let 
$\nu\in\O$ with Pl\"ucker coordinates $(p_b\,|\, b\in\mathcal B)$ and 
$(q_c\,|\, c\in\mathcal C)$. Since 
$\nu\in\O$, there is $(s,t)\in T$ such that 
$p_b=s^b\wt p_b$ for every $b\in\mathcal B$ and $q_c=t^c\wt q_c$ for 
every $c\in\mathcal C$. Clearly, (9.2.4a--b) hold. 
As for (9.2.4c), it holds if and 
only if 
$s^{\tau(b_1-b_2)}=t^{\lambda(c_2-c_1)}$ for all $b_1,b_2\in\mathcal B$ and 
$c_1,c_2\in\mathcal C$ satisfying (9.2.5). Since $(s,t)\in T$,
we need only show that 
\[
(b_1\cap J)-(b_2\cap J)=b_1-b_2=c_2-c_1=(c_2\cap I)-(c_1\cap I).
\]
But the three equalities above hold because 
(9.2.5) holds.
So $\nu\in Z$, finishing Step~1.

\emph{Step~2}: Let $\nu\in Z$ with Pl\"ucker coordinates 
$(p_b\,|\, b\in\mathcal B)$ and $(q_c\,|\, c\in\mathcal C)$. 
Then the following four conditions hold.
\begin{subequations}
\begin{align}
&\mathcal B_\nu=\{b\in\mathcal B\,|\, I'_\nu\subseteq b\subseteq I-I''_\nu\},\\
&\mathcal C_\nu=\{c\in\mathcal C\,|\, J'_\nu\subseteq c\subseteq J-J''_\nu\},\\
&\begin{cases}
I'_{\nu}\cap J\not\subseteq J'_{\nu}\cap I\text{ or }
J''_{\nu}\cap I\not\subseteq I''_{\nu}\cap J\Longrightarrow
J\cap(I-I'_{\nu})\subseteq J''_{\nu}\cap I,&\\
J'_{\nu}\cap I\not\subseteq I'_{\nu}\cap J\text{ or }
I''_{\nu}\cap J\not\subseteq J''_{\nu}\cap I\Longrightarrow
I\cap (J-J'_{\nu})\subseteq I''_{\nu}\cap J.
\end{cases}\\
&\text{There is }(\ol s,\ol t)\in\ol T_{\nu}
\text{ such that }
\begin{cases}
p_{b_1}\wt p_{b_2}=\ol s^{(b_1-b_2)}\wt p_{b_1}p_{b_2}
&
\text{for all }b_1,b_2\in\mathcal B_{\nu},\\
q_{c_1}\wt q_{c_2}=\ol t^{(c_1-c_2)}\wt q_{c_1}q_{c_2}
&
\text{for all }c_1,c_2\in\mathcal C_{\nu}.
\end{cases}
\end{align}
\end{subequations}

Indeed, (9.2.6a--b) follow as in the proof of Lemma~9.1. 

Let's check (9.2.6c) now. 
Suppose 
there are $i\in I'_{\nu}\cap(J-J'_{\nu})$ 
and $j\in(J-J''_{\nu})\cap(I-I'_{\nu})$. Since 
$i\in I'_\nu$ and $j\not\in I'_\nu$, there is $b\in\mathcal B_\nu$ such that 
$i\in b$ and $j\not\in b$. Since $i\not\in J'_\nu$ and 
$j\not\in J''_\nu$, and since (9.2.6b) holds, 
there is $c\in\mathcal C_{\nu}$ such that $i\not\in c$ and $j\in c$. 
Since $b+c=(b-i+j)+(c-j+i)$, Equations~(9.2.4c) imply that 
$b-i+j\in\mathcal B_{\nu}$. 
However, $i\in I'_{\nu}$, 
reaching
a contradiction. 
Analogously, if $J'_{\nu}\cap (I-I'_{\nu})\neq\emptyset$ 
then $I\cap(J-J'_{\nu})\subseteq I''_{\nu}\cap J$. 

Suppose now that there are $i\in J''_{\nu}\cap(I-I''_{\nu})$ and 
$j\in(J-J''_{\nu})\cap(I-I'_{\nu})$. Since $i\not\in I''_{\nu}$ and 
$j\not\in I'_{\nu}$, and since (9.2.6a) holds, 
there is $b\in\mathcal B_{\nu}$ such that $i\in b$ and $j\not\in b$. Since 
$i\in J''_\nu$ and $j\not\in J''_\nu$, there is 
$c\in\mathcal C_\nu$ such that $i\not\in c$ and $j\in c$. 
Since $b+c=(b-i+j)+(c-j+i)$, it follows from Equations~(9.2.4c) 
that $c-j+i\in\mathcal C_{\nu}$. 
However, $i\in J''_{\nu}$, 
reaching
a contradiction. 
Analogously, if $I''_{\nu}\cap(J-J''_{\nu})\neq\emptyset$ 
then $I\cap(J-J'_{\nu})\subseteq I''_{\nu}\cap J$. So (9.2.6c) is proved.

Let's check (9.2.6d). As in the proof of Lemma~9.1, there are 
$\ol s\in k^*_{\ol I_\nu}$ 
and $\ol t\in k^*_{\ol J_\nu}$ such that
\begin{equation}
\begin{cases}
p_{b_1}\wt p_{b_2}=\ol s^{(b_1-b_2)}\wt p_{b_1}p_{b_2}&
\text{ for all }b_1,b_2\in\mathcal B_\nu,\\
q_{c_1}\wt q_{c_2}=\ol t^{(c_1-c_2)}\wt q_{c_1}q_{c_2}&
\text{ for all }c_1,c_2\in\mathcal C_\nu.
\end{cases}
\end{equation}
If $i$ and $j$ are distinct elements of $\ol I_{\nu}\cap\ol J_{\nu}$, 
let $b\in\mathcal B_{\nu}$ and $c\in\mathcal C_{\nu}$ such that 
$i\not\in b\cup c$ and $j\in b\cap c$. Since $b+(c-j+i)=(b-j+i)+c$,
Equations~(9.2.4c) and 
(9.2.7) 
imply that 
$\ol s_i^\tau\ol t_j^\lambda=\ol s_j^\tau\ol t_i^\lambda$. 
So $(\ol s,\ol t)\in\ol T_{\nu}$, and hence (9.2.6d) holds. Step 2 is 
complete.

\emph{Step~3}: Let $\nu\in G_1\times G_2$ with 
Pl\"ucker coordinates 
$(p_b\,|\, b\in\mathcal B)$ and $(q_c\,|\, c\in\mathcal C)$. If 
(9.2.6a--d) hold, then $\nu\in\ol{\O}$. 
Indeed, for $u\in\Z_I$, 
$v\in\Z_J$, $s\in k^*_I$ 
and $t\in k^*_J$, define the map
\[
\zeta:=(\zeta_1,\zeta_2)\:k^*\longrightarrow k^*_I\times k^*_J,
\]
by letting $\zeta_1(r)_i:=s_ir^{u_i}$ for each $i\in I$ and 
$\zeta_2(r)_j:=t_jr^{v_j}$ for 
each $j\in J$. We'll choose $u$, $v$, $s$ and $t$ as indicated below.

Let $(\ol s,\ol t)\in\ol T_{\nu}$ as in (9.2.6d). 
For each $i\in\ol I_{\nu}$ 
set $u_i:=0$ and $s_i:=\ol s_i$. For each $j\in\ol J_{\nu}$ set 
$v_j:=0$ and $t_j:=\ol t_j$.

For each $i\in (I-\ol I_{\nu})\setminus J$ and each 
$j\in (J-\ol J_{\nu})\setminus I$
set $s_i:=1$ and $t_j:=1$, 
and
put
\[
u_i:=
\begin{cases}
-1
& \text{if $i\in I'_{\nu}$,}\\
1
&  \text{if $i\in I''_{\nu}$,}
\end{cases}
\quad\text{and}\quad v_j:=
\begin{cases}
-1
& \text{if $j\in J'_{\nu}$,}\\
1
&  \text{if $j\in J''_{\nu}$.}
\end{cases}
\]

Choose the remaining $s_i$ and $t_j$ in such a way that $(s,t)\in T$. 
As for the 
remaining $u_i$ and $v_j$ 
there are three cases to consider.

\emph{Case 1}: If
\begin{equation}
I\cap J=(I'_{\nu}\cap J'_{\nu})\cup(\ol I_{\nu}\cap\ol J_{\nu})\cup
(I''_{\nu}\cap J''_{\nu}),
\end{equation}
set 
$u_i:=-\lambda$ and $v_i:=-\tau$ if 
$i\in I'_{\nu}\cap J'_{\nu}$ and $u_i:=\lambda$ and $v_i:=\tau$ if
$i\in I''_{\nu}\cap J''_{\nu}$.
 
\emph{Case 2}: In case 
(9.2.8) 
does not hold, but
\begin{equation}
I\cap J=(I'_\nu\cap J'_\nu)\cup(\ol I_\nu\cap J'_\nu)\cup
(I''_\nu\cap J'_\nu)\cup(I''_\nu\cap\ol J_\nu)\cup(I''_\nu\cap J''_\nu),
\end{equation}
set
\begin{equation}
u_i:=
\begin{cases}
-\lambda
&	\text{ if }i\in I'_{\nu}\cap J'_{\nu},\\
\lambda
& 	\text{ if }i\in I''_{\nu}\cap J'_{\nu},\\
2\lambda
& 	\text{ if }i\in I''_{\nu}\cap\ol J_{\nu},\\
3\lambda
&	
\text{ if }i\in I''_{\nu}\cap J''_{\nu},
\end{cases}\quad\text{and}\quad v_i:=
\begin{cases}
-3\tau
&	\text{ if }i\in J'_{\nu}\cap I'_{\nu},\\
-2\tau
& 	\text{ if }i\in J'_{\nu}\cap\ol I_{\nu},\\
-\tau
& 	\text{ if }i\in J'_{\nu}\cap I''_{\nu},\\
\tau
& 	\text{ if }i\in J''_{\nu}\cap I''_{\nu}.
\end{cases}
\end{equation}

\emph{Case 3}: In case 
(9.2.8) 
does not hold, but
\begin{equation}
I\cap J=(I'_\nu\cap J'_\nu)\cup(I'_\nu\cap\ol J_\nu)\cup
(I'_\nu\cap J''_\nu)\cup(\ol I_\nu\cap J''_\nu)\cup(I''_\nu\cap J''_\nu),
\end{equation}
set $u_i$ and $v_i$ as in 
(9.2.10), 
but exchanging $I$ with $J$, 
$\lambda$ with $\tau$, and 
the $u_i$ with the $v_i$.

If 
(9.2.9) and (9.2.11) hold, then (9.2.8) 
holds as well. So 
the three cases above are independent. Moreover, 
one of the three cases occur. Indeed, if 
$J\cap(I-I'_\nu)\not\subseteq J''_\nu$ and 
$I\cap(J-J'_\nu)\not\subseteq I''_\nu$, 
then 
(9.2.8) 
occurs by (9.2.6c). On the other hand, if 
$J\cap(I-I'_{\nu})\subseteq J''_{\nu}$ then 
(9.2.11) 
holds, and if 
$I\cap(J-J'_{\nu})\subseteq I''_{\nu}$, then 
(9.2.9) 
holds.

We chose $\zeta$ in such a way that it factors through $T$. 
Indeed, since $(s,t)\in T$, it is enough to check that
$\tau(u_i-u_j)=\lambda(v_i-v_j)$ 
for all $i,j\in I\cap J$, what 
can 
easily be
done in each of the three cases above. 

Let $\nu_r:=\zeta(r)\wt\nu$ for each
$r\in k^*$, and let $(p_b(r)\,|\, b\in\mathcal B)$ and 
$(q_c(r)\,|\, c\in\mathcal C)$ be Pl\"ucker coordinates of 
$\nu_r$. Because of our choice of $\zeta$, the power of $r$ 
in the expression for $p_b(r)$ is minimum exactly among those 
$b\in\mathcal B_{\nu}$. Likewise, the power of $r$ 
in the expression for $q_c(r)$ is minimum exactly among those 
$c\in\mathcal C_{\nu}$. Let $\ol\nu\in G_1\times G_2$ be 
the limit of $\nu_r$ as $r$ goes to 0. As in the proof of Lemma~9.1, we 
have $\ol\nu=\nu$. Then $\nu\in\ol{\O}$, finishing Step~3.

\emph{Step~4}: Let $\nu\in\ol{\O}$. Then there are ordered 
tripartitions $\mathbb I=(I,\ol I,I'')$ and 
$\mathbb J=(J',\ol J,J'')$ of $I$ and $J$, respectively, 
such that (9.2.2--3) hold, and such that $\nu$ is in the 
orbit of $(V_{\mathbb I},W_{\mathbb J})$ defined by (9.2.1). 
Indeed, let $(p_b\,|\, b\in\mathcal B)$ and $(q_c\,|\, c\in\mathcal C)$ 
be the Pl\"ucker coordinates of $\nu\in\ol{\O}$. 
By Step~1, $\nu\in Z$. Hence (9.2.6a--d) hold by Step~2. Since (9.2.6d) 
holds, we may replace $\nu$ by a certain point in its orbit and assume that 
$p_b=\wt p_b$ for every $b\in\mathcal B_{\nu}$ and 
$q_c=\wt q_c$ for every $c\in\mathcal C_{\nu}$.
Let $\wt{\mathbb I}:=(I'_\nu,{\ol I}_\nu,I''_\nu)$ and 
$\wt{\mathbb J}:=(J'_\nu,{\ol J}_\nu,J''_\nu)$. As in the proof of 
Lemma~9.1, we have $\nu=(V_{\wt{\mathbb I}},W_{\wt{\mathbb J}})$. 

There are four cases to consider.

\emph{Case 1}: If $|I'_{\nu}|<h_1$ and $|J'_{\nu}|<h_2$, 
set $\mathbb I:=\wt{\mathbb I}$ and $\mathbb J:=\wt{\mathbb J}$. Clearly
(9.2.2) holds, 
whereas
(9.2.3a--b) hold by (9.2.6c). Since 
$\nu=(V_{\mathbb I},W_{\mathbb J})$, Step~4 is complete for Case 1.

\emph{Case 2}: If $|I'_\nu|=h_1$ and $|J'_\nu|<h_2$, then 
set $\mathbb I:=(I'_\nu-\ol I,\ol I,I''_\nu)$ for a certain 
non--empty subset $\ol I\subseteq I'_\nu$ to be defined below. Set 
$\mathbb J:=\wt{\mathbb J}$. As in the proof of Lemma~9.1, we 
have $\nu=(V_{\mathbb I},W_{\mathbb J})$.

Clearly
(9.2.2) holds. We need only choose $\ol I$ such that (9.2.3a--b) 
hold, and we'll do it according to whether 
(9.2.8), (9.2.9) or (9.2.11) occurs. If (9.2.8) or (9.2.9)
occurs, set $\ol I:=I'_{\nu}$. If 
(9.2.8) does not occur but (9.2.11) 
does, set
\[
\ol I:=\begin{cases}
I'_\nu\setminus J
&\text{if $I'_\nu\not\subseteq J$,}\\
I'_\nu\cap J''_\nu
&\text{if $I'_\nu\subseteq J$ and 
$I'_\nu\not\subseteq(J'_\nu\cup\ol J_\nu)$,}\\
I'_\nu\cap\ol J_\nu
&\text{if $I'_\nu\subseteq(J'_\nu\cup\ol J_\nu)$ and 
$I'_\nu\not\subseteq J'_\nu$,}\\
I'_\nu
&\text{if $I'_\nu\subseteq J'_\nu$.}
\end{cases}
\]
With any of the above choices (9.2.3a--b) hold 
for $\mathbb I$ and $\mathbb J$.

\emph{Case 3}: If $|I'_\nu|<h_1$ and $|J'_\nu|=h_2$, proceed as in 
Case 2, exchanging $I$ with $J$.

\emph{Case 4}: If $|I'_\nu|=h_1$ and $|J'_\nu|=h_2$, let 
$\ol I\subseteq I'_\nu$ and $\ol J\subseteq J'_\nu$ be non--empty subsets, 
and put
\[
\mathbb I:=(I'_\nu-\ol I,\ol I,I''_\nu)\quad\text{and}\quad
\mathbb J:=(J'_\nu-\ol J,\ol J,J''_\nu). 
\]
As in Case 2, we need only choose $\ol I$ and $\ol J$ such that (9.2.3a--b) 
hold. If (9.2.8) 
occurs, set $\ol I:=I'_\nu$ and $\ol J:=J'_\nu$. If 
(9.2.9) occurs but (9.2.8) 
does not, set 
$\ol I:=I'_\nu$ and $\ol J:=I''_\nu\cap J'_\nu$. If 
(9.2.11) occurs but (9.2.8) 
does not, set $\ol I:=I'_\nu\cap J''_\nu$ and $\ol J:=J'_\nu$. 
With any of the above choices
(9.2.3a--b) hold for $\mathbb I$ and $\mathbb J$.
Step~4 is complete.

\emph{Step~5}: If $\nu=(V_{\mathbb I},W_{\mathbb J})$, 
where $\mathbb I=(I',\ol I,I'')$ and $\mathbb J=(J',\ol J,J'')$ are 
ordered tripartitions of $I$ and $J$ such 
that 
(9.2.2--3) 
hold, then $\nu\in\ol{\O}$. Indeed, 
let $(p_b\,|\, b\in\mathcal B)$ and $(q_c\,|\, c\in\mathcal C)$ 
be Pl\"ucker coordinates of $\nu$. As in the proof of 
Lemma~9.1,
\[
\mathcal B_{\nu}=\{b\in\mathcal B|I'\subseteq b\subseteq I-I''\}.
\]
If $h_1<|I-I''|$
then $I'_\nu=I'$ and $I''_\nu=I''$. 
If $h_1=|I-I''|$
then $I'_\nu=I'\cup\ol I$ and $I''_\nu=I''$. 
In either case
(9.2.6a) holds. Analogously, (9.2.6b) holds. Now, (9.2.6c) is simply a 
restatement of (9.2.3) if $h_1<|I-I''|$ and $h_2<|J-J''|$. In the 
remaining cases, (9.2.6c) can be obtained from (9.2.3) as well. 
Finally, 
as in the proof of Lemma~9.1, Condition 
(9.2.6d) holds trivially with $(\ol s,\ol t)=1$. By Step~3, 
$\nu\in\ol{\O}$, thus finishing Step~5.

Steps 4 and 5 show the lemma.
\end{proof}

\section{The variety of limit canonical systems}

\begin{theorem}\setcounter{equation}{0} 
Preserve {\rm 4.1} and assume {\rm (4.3.1)}. 
Let $\mu,\mu'\in\Z_\Delta^+$. Then the following five assertions hold.
\begin{enumerate}
\item $\V_{\mu,X}=\V_{\mu',X}$ if and only if $\alpha_\mu=\alpha_{\mu'}$ 
and either $I_\mu=I_{\mu'}$ or $|\alpha_\mu|=g_Y$.
\item $\V_{\mu,Y}=\V_{\mu',Y}$ if and only if $\beta_\mu=\beta_{\mu'}$ 
and either $J_\mu=J_{\mu'}$ or $|\beta_\mu|=g_X$.
\item $\V_\mu=\V_{\mu'}$ if and only if $\V_{\mu,X}=\V_{\mu',X}$ and 
$\V_{\mu,Y}=\V_{\mu',Y}$.
\item Either $\V_\mu=\V_{\mu'}$ or $\V_\mu\cap\V_{\mu'}=\emptyset$.
\item If $\mu'=t\mu$ for a certain $t\in\Z^+$, then $\V_\mu=\V_{\mu'}$.
\end{enumerate}
\end{theorem}

\begin{proof} We take two preliminary steps.

\emph{Step~1}: If $\V_{\mu,X}\cap\V_{\mu',X}\neq\emptyset$, then 
$\alpha_\mu=\alpha_{\mu'}$ and either $I_\mu=I_{\mu'}$ or 
$|\alpha_\mu|=g_Y$. Indeed, let $V\in\G_X$, and assume that 
$V\in\V_{\mu,X}\cap\V_{\mu',X}$. Let
\[
L_X:=\w_X(\textstyle\sum_{p\in\Delta}(1+\alpha_{\mu,p})x_p)
\quad\text{and}\quad
L_Y:=\w_Y(\textstyle\sum_{p\in I_\mu}y_p-
\textstyle\sum_{p\in\Delta}\alpha_{\mu,p}y_p).
\]
Then $V\subseteq H^0(L_X)$. 
Now, $x_p$ is not a base point of 
$(V,L_X)$ for any $p\in\Delta$ by Theorem~5.2. 
So $g_Y-\alpha_{\mu,p}$ is the minimum order of 
$(V,\mathbb L)$ at $x_p$ for each $p\in\Delta$. Analogously, 
$g_Y-\alpha_{\mu',p}$ is the minimum order of 
$(V,\mathbb L)$ at $x_p$ for each $p\in\Delta$. So 
$\alpha_\mu=\alpha_{\mu'}$.

Assume now that 
$|\alpha_\mu|>g_Y$. Let $F\subseteq\Delta$. We claim 
that
\begin{equation}
V\supseteq H^0(L_X(-\textstyle\sum_{p\in F}x_p))
\end{equation}
if and only if $F\supseteq I_\mu$. In fact, 
$V\supseteq H^0(L_X(-\sum_{p\in I_\mu}x_p))$ by Theorem~5.2. So, 
if $F\supseteq I_\mu$ then (10.1.1) holds. Conversely, suppose 
(10.1.1) holds. Since $\alpha_\mu\neq 0$, 
by Riemann--Roch, 
$\{x_p\,|\, p\in\Delta\}$ imposes independent conditions on $H^0(L_X)$. 
So, since $V\supseteq H^0(L_X(-\sum_{p\in I_\mu}x_p))$ by Theorem~5.2, it 
follows from (10.1.1) that $V\supseteq H^0(L_X(-\sum_{p\in F'}x_p))$, 
where $F':=F\cap I_\mu$. We may thus assume that $F'=F$.

Suppose by contradiction that $F\neq I_\mu$. Replacing $F$ by a larger 
subset, if necessary, we may assume $|I_\mu-F|=1$.  
Since $\{x_p\,|\, p\in\Delta\}$ imposes 
independent conditions on $H^0(L_X)$, 
there is $s\in H^0(L_X)$ such that $s(x_p)=0$ for every $p\in F$ and 
$s(x_p)\neq 0$ for $p\in I_\mu-F$. So $s\in V$ by (10.1.1). 
Let 
$C_{I_\mu}$ 
be the blow--up of $C$ along 
$\Delta-I_\mu$. By Theorem~8.2, there is an invertible sheaf $L$ on 
$C_{I_\mu}$ 
such that $L|_X\cong L_X$ and 
$L|_Y\cong L_Y$,
and such that $V$ is the 
image of the restriction map $H^0(L)\to H^0(L_X)$. Let $\wt s\in H^0(L)$ 
lifting $s$ and $t:=\wt s|_Y$. So $t\in H^0(L_Y(-\sum_{p\in F}y_p))$. 
Since $|I_\mu-F|=1$ and $|\alpha|>g_Y$, it follows from (4.3.2a) that 
$h^0(L_Y(-\sum_{p\in F}y_p))=0$. So $t=0$, and hence $s(x_p)=0$ for 
every $p\in\Delta$. The contradiction proves our claim.

Now, $\alpha_{\mu}=\alpha_{\mu'}$. By analogy, (10.1.1) holds if and only if 
$F\supseteq I_{\mu'}$. So
$I_\mu=I_{\mu'}$. Step~1 is complete.

\emph{Step~2}: Assume that 
$\alpha_\mu=\alpha_{\mu'}$ and 
either $I_\mu=I_{\mu'}$ or 
$|\alpha_\mu|=g_Y$. 
Then 
$\V_{\mu,X}=\V_{\mu',X}$. 
Indeed, it follows from Theorem~8.2 and the hypothesis that 
$\V_{\mu,X}=\V_{\mu',X}$ if $|\alpha_\mu|>g_Y$. Now, if 
$|\alpha_\mu|=g_Y$ then
\[
\V_{\mu',X}=\V_{\mu,X}=\{H^0(L_X)\}
\]
by Theorem~8.5. Step~2 is complete.

Steps 1 and 2 prove Statement~1. Statement~2 follows by analogy. Let's 
prove Statement~3. By definition, $\V_{\mu,X}$ and $\V_{\mu,Y}$ are the 
images of $\V_\mu$ in $\G_X$ and $\G_Y$, respectively, for 
each $\mu\in\Z^\Delta_+$. Hence, if $\V_\mu=\V_{\mu'}$ then 
$\V_{\mu,X}=\V_{\mu',X}$ and $\V_{\mu,Y}=\V_{\mu',Y}$.

Conversely, assume 
that
$\V_{\mu,X}=\V_{\mu',X}$ and 
$\V_{\mu,Y}=\V_{\mu',Y}$. By Statements~1~and~2, 
we have
$(\alpha_\mu,\beta_\mu)=(\alpha_{\mu'},\beta_{\mu'})$. In addition, 
$I_\mu=I_{\mu'}$ when $|\alpha_\mu|>g_Y$ and $J_\mu=J_{\mu'}$ 
when $|\beta_\mu|>g_X$.
So, if $|\alpha_\mu|>g_Y$ and $|\beta_\mu|>g_X$ then 
$\V_\mu=\V_{\mu'}$ by Theorem~8.2. Now, if $|\alpha_\mu|=g_Y$ or 
$|\beta_\mu|=g_X$ 
then
\[
\V_\mu=\V_{\mu,X}\times\V_{\mu,Y}=\V_{\mu',X}\times\V_{\mu',Y}=\V_{\mu'},
\]
where the first and last equalities follow from Theorem~8.5. Statement~3 
is proved.

Let's prove Statement~4. Suppose that $\V_\mu\cap\V_{\mu'}\neq\emptyset$. 
Then $\V_{\mu,X}\cap\V_{\mu',X}\neq\emptyset$ as well. Thus 
$\V_{\mu,X}=\V_{\mu',X}$ by Step~1. Analogously, $\V_{\mu,Y}=\V_{\mu',Y}$. 
By Statement~3, $\V_\mu=\V_{\mu'}$.
Statement~4 is proved.

Finally, let's prove Statement~5. Suppose $\mu'=t\mu$ for a certain 
$t\in\Z^+$. Then the numerical 
data $(\alpha_\mu,I_\mu)$ 
associated with $\mu$ and $g_Y$ is equal to that associated with 
$\mu'$ and $g_Y$. Analogously, the numerical data 
$(\beta_\mu,J_\mu)$ associated with 
$\mu$ and $g_X$ is equal to that associated with 
$\mu'$ and $g_X$. Thus $\V_\mu=\V_{\mu'}$ by Statements~1,~2~and~3.
\end{proof}

\begin{definition}\setcounter{equation}{0} 
Preserve 4.1 and assume (4.3.1). For each $\mu\in\Q_\Delta^+$ let 
\[
\V_\mu:=\V_{t\mu},\quad\V_{\mu,X}:=\V_{t\mu,X},\quad\V_{\mu,Y}:=\V_{t\mu,Y},
\]
where $t$ is any positive integer such that $t\mu\in\Z_\Delta^+$. 
(By Theorem~10.1, $\V_\mu$ is well defined.) 
Give $\Q_\Delta^+$ the topology 
induced by the Euclidean topology of $\R_\Delta$. Let 
\[
\V:=\bigcup_{\mu\in\Q_\Delta^+}\V_\mu\subseteq\G.
\]
We call $\V$ the \emph{variety of limit canonical systems of $C$}. 
We call each $\V_\mu$ a \emph{stratum} of $\V$. 
\end{definition}

\begin{proposition}\setcounter{equation}{0} 
Preserve {\rm 4.1} and assume {\rm (4.3.1)}. 
For each $\mu\in\Q_\Delta^+$ let
\[
\mathbb U_\mu:=\{\mu'\in\Q_\Delta^+\,|\, \V_{\mu'}=\V_\mu\}.
\]
Then each $\mathbb U_\mu$ is convex and homogeneous, and the covering 
$\{\mathbb U_\mu\,|\, \mu\in\Q_\Delta^+\}$ of $\Q_\Delta^+$ is finite.
\end{proposition}

\begin{proof} By Theorem~10.1, $\V_{t\mu}=\V_\mu$ for each 
$\mu\in\Q_\Delta^+$ and each 
$t\in\Q^+$. 
In other words, each 
$\mathbb U_\mu$ is homogeneous.

Let $\mu,\mu'\in\Q_\Delta^+$. Clearly, $\mathbb U_\mu=\mathbb U_{\mu'}$ if and only 
if $\V_\mu=\V_{\mu'}$. So, by Theorem~10.1, 
$\mathbb U_\mu=\mathbb U_{\mu'}$ if and only if 
$(\alpha_\mu,\beta_\mu)=(\alpha_{\mu'},\beta_{\mu'})$,
\begin{equation}
I_\mu=I_{\mu'}\text{ when }|\alpha_\mu|>g_Y\quad\text{and}\quad
J_\mu=J_{\mu'}\text{ when }|\beta_\mu|>g_X.
\end{equation}
On the other hand, $0\leq\alpha_{\mu,p}\leq g_Y$ and 
$0\leq\beta_{\mu,p}\leq g_X$ 
for each $p\in\Delta$ and $I_\mu,J_\mu\subseteq\Delta$. So, as 
$\mu$ runs over $\Q_\Delta^+$, the associated numerical data 
$\alpha_\mu,\beta_\mu,I_\mu,J_\mu$ runs over a finite set. Therefore, 
the covering $\{\mathbb U_\mu\,|\, \mu\in\Q_\Delta^+\}$ of $\Q_\Delta^+$ 
is finite.

Assume now that $\mathbb U_\mu=\mathbb U_{\mu'}$. 
Let $t\in [0,1]\cap\Q$. Set 
\[
\mu(t):=t\mu+(1-t)\mu',\quad \rho(t):=t\rho_\mu+(1-t)\rho_{\mu'},\quad 
\sigma(t):=t\sigma_\mu+(1-t)\sigma_{\mu'},
\]
and put
\begin{align*}
I(t):=&
\{p\in\Delta\,|\, \rho(t)_p\geq\rho(t)_q\text{ for every }q\in\Delta\},\\
J(t):=&
\{p\in\Delta\,|\, \sigma(t)_p\leq\sigma(t)_q\text{ for every }q\in\Delta\}.
\end{align*}
Then $I(t)=I_\mu\cap I_{\mu'}$ if $I_\mu\cap I_{\mu'}\neq\emptyset$ and 
$J(t)=J_\mu\cap J_{\mu'}$ if $J_\mu\cap J_{\mu'}\neq\emptyset$. 
By (10.3.1),
\begin{equation}
I(t)=I_\mu\text{ when }|\alpha_\mu|>g_Y\quad\text{and}\quad
J(t)=J_\mu\text{ when }|\beta_\mu|>g_X.
\end{equation}
Now, since $(\alpha_\mu,\beta_\mu)=(\alpha_{\mu'},\beta_{\mu'})$, 
it follows that 
\[
(\alpha_{\mu(t)},\beta_{\mu(t)})=(\alpha_\mu,\beta_\mu)\quad\text{and}
\quad(I_{\mu(t)},J_{\mu(t)})=(I(t),J(t)).
\]
Then $\V_{\mu(t)}=\V_\mu$ by (10.3.2) and Theorem~10.1.
So $\mathbb U_\mu$ is convex.
\end{proof}

\begin{lemma}\setcounter{equation}{0} 
Preserve {\rm 4.1}. Then, for each $\mu\in\Q_\Delta^+$ there exists an open 
neighbourhood $U_\mu\subseteq\Q_\Delta^+$ of $\mu$ such 
that, for each open neighbourhood $U\subseteq U_\mu$ of $\mu$, the 
following four statements hold.
\begin{enumerate}
\item If $\ol\mu\in U$
then there are unique ordered tripartitions $(I',\ol I,I'')$ of $I_\mu$ and 
$(J',\ol J,J'')$ of $J_\mu$ satisfying the 
following six properties.
\begin{subequations}
\begin{equation}
\begin{align}
&g_Y+|I''|\leq|\alpha_\mu|<g_Y+|I_\mu-I'|,\\
&g_X+|J''|\leq|\beta_\mu|<g_X+|J_\mu-J'|,\\
&\ol I=I_{\ol\mu},\\
&\alpha_{\ol\mu,p}=\alpha_{\mu,p}-1\text{ if $p\in I''$ and }
\alpha_{\ol\mu,p}=\alpha_{\mu,p}\text{ otherwise,}\\
&\ol J=J_{\ol\mu},\\
&\beta_{\ol\mu,p}=\beta_{\mu,p}-1\text{ if $p\in J''$ and }
\beta_{\ol\mu,p}=\beta_{\mu,p}\text{ otherwise.}
\end{align}
\end{equation}
If $g_Xg_Y>0$, then $(I',\ol I,I'')$ and $(J',\ol J,J'')$ 
satisfy the two properties below as well:
\begin{equation}
\begin{align}
I'\cap J_\mu\not\subseteq J'\cap I_\mu\text{ or }
J''\cap I_\mu\not\subseteq I''\cap J_\mu\Longrightarrow
J_\mu\cap(I_\mu-I')\subseteq J''\cap I_\mu,&\\
J'\cap I_\mu\not\subseteq I'\cap J_\mu\text{ or }
I''\cap J_\mu\not\subseteq J''\cap I_\mu\Longrightarrow
I_\mu\cap(J_\mu-J')\subseteq I''\cap J_\mu.&
\end{align}
\end{equation}
\end{subequations}
\item If $g_Y>0$ and $(I',\ol I,I'')$ is 
an ordered tripartition of $I_\mu$ satisfying {\rm (10.4.1a)}, 
then there is $\ol\mu\in U$ satisfying {\rm (10.4.1c,d)}.
\item If $g_X>0$ and $(J',\ol J,J'')$ is 
an ordered tripartition of $J_\mu$ satisfying {\rm (10.4.1b)}, 
then there is $\ol\mu\in U$ satisfying {\rm (10.4.1e,f)}.
\item If $g_Xg_Y>0$ and $(I',\ol I,I'')$ and 
$(J',\ol J,J'')$ are ordered tripartitions of $I_\mu$ and 
$J_\mu$ satisfying {\rm (10.4.1a,b,g,h)}, 
then there is $\ol\mu\in U$ satisfying {\rm (10.4.1c--f)}.
\end{enumerate}
\end{lemma}

\begin{proof} To ease notation, let 
\[
(\alpha,\rho,I,\beta,\sigma,J):=
(\alpha_\mu,\rho_\mu,I_\mu,\beta_\mu,\sigma_\mu,J_\mu).
\]
Let $U_\mu\subseteq\Q_\Delta$ consist of the points of the form 
$\mu+\varepsilon$ for all $\varepsilon\in\Q_\Delta$ 
satisfying
\begin{equation}
|\varepsilon_p|<\frac{\min(\min^*
(\rho_q,\mu_q-\rho_q,\sigma_q,\mu_q-\sigma_q)\,\,|\, q\in\Delta)}
{3(1+\max(\max(\alpha_q,\beta_q)\,\,|\, q\in\Delta))}
\  \text{ for every $p\in\Delta$,}
\end{equation}
where ``$\min^*$'' means the smallest non--zero number among 
the four listed. If (10.4.2) holds, then $|\varepsilon_p|<\mu_p$ 
for every $p\in\Delta$, and hence $U_\mu\subseteq\Q_\Delta^+$.

We take three preliminary steps.

\emph{Step~1}: For each $\varepsilon\in\Q_\Delta$ there are 
ordered tripartitions $(I',\ol I, I'')$ of $I$ and 
$(J',\ol J,J'')$ of $J$ such that 
(10.4.1a,b) and (10.4.3a,b) below hold.
\begin{subequations}
\begin{align}
&\varepsilon_p\alpha_p<\varepsilon_q\alpha_q=
\varepsilon_r\alpha_r<\varepsilon_s\alpha_s
\text{ for all $p\in I'$, all $q,r\in\ol I$ and all 
$s\in I''$,}\\
&\varepsilon_p\beta_p<\varepsilon_q\beta_q=
\varepsilon_r\beta_r<\varepsilon_s\beta_s
\text{ for all $p\in J'$, all $q,r\in\ol J$ and all 
$s\in J''$.}
\end{align}
\end{subequations}
Indeed, for each $\tau\in\Q$ let 
\[
I'_\tau:=\{p\in I\,|\, \varepsilon_p\alpha_p<\tau\},\quad
\ol I_\tau:=\{p\in I\,|\, \varepsilon_p\alpha_p=\tau\},\quad
I''_\tau:=\{p\in I\,|\, \varepsilon_p\alpha_p>\tau\}.
\]
Since $g_Y\leq|\alpha|<g_Y+|I|$, 
there is a unique $\lambda_\varepsilon\in\Q$ such that
\[
g_Y+|I''_{\lambda_\varepsilon}|\leq|\alpha|<g_Y+|I-I'_{\lambda_\varepsilon}|.
\]
Let $(I',\ol I,I''):=(I'_{\lambda_\varepsilon},\ol I_{\lambda_\varepsilon},
I''_{\lambda_\varepsilon})$. Then (10.4.1a) and (10.4.3a) 
hold. By analogy, we construct $(J',\ol J,J'')$ satisfying (10.4.1b) and 
(10.4.3b). Step~1 is complete.

\emph{Step~2}: 
If $g_Xg_Y>0$ and 
there are 
ordered tripartitions $(I',\ol I, I'')$ of $I$ and 
$(J',\ol J,J'')$ of $J$ 
satisfying (10.4.3)
for a certain $\varepsilon\in\Q_\Delta$, then 
(10.4.1g,h) hold. Indeed, if $p,q\in I\cap J$ then
$\mu_p\alpha_p=\mu_q\alpha_q$ and 
$\mu_p\beta_p=\mu_q\beta_q$. Hence
\begin{equation}
\alpha_p\beta_q=\alpha_q\beta_p
\text{ for all $p,q\in I\cap J$.}
\end{equation}
Assume that $g_Xg_Y>0$. Let $(I',\ol I, I'')$ and $(J',\ol J,J'')$ be 
ordered tripartitions of $I$ and $J$ satisfying (10.4.3) for a certain 
$\varepsilon\in\Q_\Delta$.
If there are 
$p\in I'\cap(J-J')$
and $q\in(I-I')\cap(J-J'')$, then 
$\varepsilon_p\alpha_p<\varepsilon_q\alpha_q$ and
$\varepsilon_q\beta_q\leq\varepsilon_p\beta_p$ by 
(10.4.3). 
Now, $\alpha_p\beta_p>0$ because $g_Xg_Y>0$. 
Using (10.4.4) we get
\[
\varepsilon_p\alpha_p\beta_p<\varepsilon_q\alpha_q\beta_p=
\varepsilon_q\beta_q\alpha_p\leq\varepsilon_p\alpha_p\beta_p,
\]
reaching
a contradiction. If there are $p\in J''\cap (I-I'')$ and 
$q\in(I-I')\cap(J-J'')$, then 
$\varepsilon_p\alpha_p\leq\varepsilon_q\alpha_q$ and 
$\varepsilon_q\beta_q<\varepsilon_p\beta_p$ by 
(10.4.3). 
Using (10.4.4) as before, we 
get $\varepsilon_p\alpha_p\beta_p<\varepsilon_p\alpha_p\beta_p$, 
reaching
a 
contradiction. So (10.4.1g) holds. By analogy, (10.4.1h) holds. 
Step~2 is complete.

\emph{Step~3}: Let $\ol\mu:=\mu+\varepsilon$ for $\varepsilon\in\Q_\Delta$ 
such that (10.4.2) holds. Let $(I',\ol I,I'')$ and 
$(J',\ol J,J'')$ be ordered tripartitions of $I$ and $J$ 
such that (10.4.1a,b) hold. If (10.4.3a) holds, then so do 
(10.4.1c,d). If (10.4.3b) holds, then so do (10.4.1e,f). 
Indeed, define $\ol\alpha,\ol\beta\in\Z_\Delta$ by letting for each 
$p\in\Delta$,
\[
\ol\alpha_p:=
\begin{cases}
\alpha_p-1&\text{if $p\in I''$},\\
\alpha_p&\text{otherwise},
\end{cases}
\quad\text{and}\quad
\ol\beta_p:=
\begin{cases}
\beta_p-1&\text{if $p\in J''$},\\
\beta_p&\text{otherwise}.
\end{cases}
\]
By (10.4.1a,b),
\[
g_Y\leq|\ol\alpha|<g_Y+|\ol I|\quad\text{and}\quad 
g_X\leq|\ol\beta|<g_X+|\ol J|.
\]

If (10.4.3a) holds, let $\lambda:=\varepsilon_p\alpha_p$ for (any) 
$p\in\ol I$, and define $\ol\rho\in\Q_\Delta$ by letting for each 
$p\in\Delta$,
\[
\ol\rho_p:=
\begin{cases}
\ol\mu_p+\varepsilon_p\alpha_p-\lambda&\text{if $p\in I'$,}\\
\ol\mu_p&\text{if $p\in\overline I$,}\\
\varepsilon_p\alpha_p-\lambda&\text{if $p\in I''$,}\\
\rho_p+\varepsilon_p(1+\alpha_p)-\lambda&\text{if $p\not\in I$.}
\end{cases}
\]
Then 
\[
\ol\mu_p(\ol\alpha_p+1)-\ol\rho_p=\ol\mu_q(\ol\alpha_q+1)-\ol\rho_q
\quad\text{for all $p,q\in\Delta$.}
\]
In addition, it follows from (10.4.2) and (10.4.3a) that 
$0<\ol\rho_p\leq\ol\mu_p$ for every $p\in\Delta$, 
with $\ol\rho_p=\ol\mu_p$ if and only if $p\in\ol I$. So 
$(\ol\alpha,\ol\rho,\ol I)=(\alpha_{\ol\mu},\rho_{\ol\mu},I_{\ol\mu})$, 
proving (10.4.1c,d). 
 
If (10.4.3b) holds, let $\tau:=\varepsilon_p\beta_p$ for 
(any) 
$p\in\ol J$, and define $\ol\sigma\in\Q_\Delta$ by letting 
for each $p\in\Delta$,
\[
\ol\sigma_p:=
\begin{cases}
\tau-\varepsilon_p\beta_p&\text{if $p\in J'$,}\\
0&\text{if $p\in\overline J$,}\\
\ol\mu_p+\tau-\varepsilon_p\beta_p&\text{if $p\in J''$,}\\
\sigma_p+\tau-\varepsilon_p\beta_p&\text{if $p\not\in J$.}
\end{cases}
\]
As before, 
$(\ol\beta,\ol\sigma,\ol J)=(\beta_{\ol\mu},\sigma_{\ol\mu},J_{\ol\mu})$, 
proving (10.4.1e,f). Step 3 is complete.

Let's prove that $U:=U_\mu$ satisfies Statement~1. Indeed, let 
$\ol\mu:=\mu+\varepsilon$ for $\varepsilon\in\Q_\Delta$ 
such that (10.4.2) holds. By Step~1, there are 
ordered tripartitions $(I',\ol I, I'')$ of $I$ and 
$(J',\ol J,J'')$ of $J$ such that 
(10.4.1a,b) and 
(10.4.3a,b) 
hold. So 
(10.4.1c--f) 
hold by Step~3. If 
$g_Xg_Y>0$ then (10.4.1g,h) hold as well by Step~2. The uniqueness of the 
tripartitions follows from 
(10.4.1c--f).

Let $U\subseteq U_\mu$ be an open neighbourhood of $\mu$. 
Let's prove Statement~2. Assume $g_Y>0$, and 
let $(I',\ol I,I'')$ be an 
ordered tripartition of $I$ 
satisfying (10.4.1a).  Define $\upsilon\in\Q_\Delta$ by letting for each 
$p\in\Delta$,
\begin{equation}
\upsilon_p:=
\begin{cases}
0
&\text{ if $p\in I'\cup(\Delta-I)$,}\\
\mu_p
&\text{ if $p\in\ol I$,}\\
2\mu_p
&\text{ if $p\in I''$.}
\end{cases} 
\end{equation}
Let $\varepsilon:=t\upsilon$, where $t\in\Q^+$ is such that 
$\ol\mu:=\mu+\varepsilon$ lies in $U$. Now, 
$\alpha_p>0$ for each $p\in I$ because $g_Y>0$. Then 
(10.4.3a) follows from 
(10.4.5). So (10.4.1c,d) 
hold by Step~3. 
Analogously, Statement~3 holds.

Let's prove Statement~4. Assume that $g_Xg_Y>0$. Let 
$(I',\ol I,I'')$ and $(J',\ol J,J'')$ be 
ordered tripartitions of $I$ and $J$ 
satisfying (10.4.1a,b,g,h). Define $\upsilon\in\Q_\Delta$ 
according to the following 
three cases.

\emph{Case 1}: 
If
\begin{equation}
I\cap J=(I'\cap J')\cup(\ol I\cap\ol J)\cup(I''\cap J''),
\end{equation}
let for each $p\in\Delta$,
\[
\upsilon_p:=
\begin{cases}
0
&\text{if $p\in I'\cup J'\cup(\Delta-(I\cup J))$,}\\
\mu_p
&\text{if $p\in\ol I\cup\ol J$,}\\
2\mu_p
&\text{if $p\in I''\cup J''$.}
\end{cases}
\]

\emph{Case 2}: 
If (10.4.6) does not hold, but
\[
I\cap J=(I'\cap J')\cup(\ol I\cap J')\cup(I''\cap J')\cup(I''\cap\ol J)\cup
(I''\cap J''),
\]
let for each $p\in\Delta$,
\[
\upsilon_p:=
\begin{cases}
0
&\text{if $p\in I'\cup(J'\setminus I)\cup(\Delta-(I\cup J))$,}\\
\mu_p
&\text{if $p\in\ol I$,}\\
2\mu_p
&\text{if $p\in I''\cap J'$,}\\
3\mu_p
&\text{if $p\in\ol J$,}\\
4\mu_p
&\text{if $p\in J''\cup(I''\setminus J)$.}
\end{cases}
\]

\emph{Case 3}: 
If (10.4.6) does not hold, but
\[
I\cap J=(I'\cap J')\cup(I'\cap\ol J)\cup(I'\cap J'')\cup(\ol I\cap J'')\cup
(I''\cap J''),
\]
define $\upsilon_p$ for each $p\in\Delta$ as in Case~2,
but now exchanging $I$ with $J$.

As in the proof of Lemma~9.2, it
follows from (10.4.1g,h) that the three cases above are independent, 
and one of them occurs.

Let $\varepsilon:=t\upsilon$, where $t\in\Q^+$ is such that 
$\ol\mu:=\mu+\varepsilon$ lies in $U$. Now, 
$\alpha_p>0$ for each $p\in I$ and $\beta_p>0$ for each $p\in J$ because 
$g_Xg_Y>0$. Then 
Conditions~(10.4.3a,b) 
hold in each of the three cases above. 
So (10.4.1c--f) hold by Step~3.
\end{proof}

\begin{theorem}\setcounter{equation}{0} 
Preserve {\rm 4.1} and 
assume {\rm (4.3.1)}. Then, for each $\mu\in\Q_\Delta^+$, 
there is an open neighbourhood $U_\mu\subseteq\Q_\Delta^+$ of $\mu$ 
such that 
the closure $\ol\V_\mu\subseteq\G$ satisfies
\[
\ol\V_\mu=\bigcup_{\ol\mu\in U}\V_{\ol\mu}
\]
for every open neighbourhood $U\subseteq U_\mu$ of $\mu$.
\end{theorem}

\begin{proof} Let $\mu\in\Q_\Delta^+$. 
By Theorem~10.1, we may assume that $\mu\in\Z_\Delta^+$. 
Let $U_\mu$ be the open neighborhood of $\mu$ given by Lemma~10.4. 
To ease notation, let
\[
(\alpha,\rho,I,\beta,\sigma,J):=
(\alpha_\mu,\rho_\mu,I_\mu,\beta_\mu,\sigma_\mu,J_\mu).
\]
In addition, if $g_Xg_Y>0$ let 
$(\wt\alpha,\wt\beta):=(\wt\alpha_\mu,\wt\beta_\mu)$.
Let 
$L_X$, $L_Y$, $M_X$ and $M_Y$ 
be the invertible sheaves given by (4.2.1). 
As in 8.3, let $C_{I\cap J}$ denote the blow--up of $C$ along 
$\Delta-(I\cap J)$. Fix isomorphisms $\zeta_{X,p}\:L_X(x_p)\to k$ and 
$\zeta_{Y,p}\:L_Y(y_p)\to k$ for each $p\in I$, and 
$\xi_{X,p}\:M_X(x_p)\to k$ and $\xi_{Y,p}\:M_Y(y_p)\to k$ for each 
$p\in J$. 
If $g_Xg_Y>0$ and $I\cap J\neq\emptyset$, choose them such that 
\[
\{(\zeta_{X,p}^{\otimes\wt\beta}\otimes\xi_{X,p}^{\otimes\wt\alpha},
\zeta_{Y,p}^{\otimes\wt\beta}\otimes\xi_{Y,p}^{\otimes\wt\alpha})\,|\, 
p\in I\cap J\}
\]
patch $L_X^{\otimes\wt\beta}\otimes M_X^{\otimes\wt\alpha}$ and 
$L_Y^{\otimes\wt\beta}\otimes M_Y^{\otimes\wt\alpha}$ 
to the sheaf $K_{I\cap J}$ given by~(8.2.2). 
Consider the corresponding evaluation maps,
\[
e_X\:H^0(L_X)\to k_I,\  e_Y\:H^0(L_Y)\to k_I,\  f_X\:H^0(M_X)\to k_J,\  
f_Y\:H^0(M_Y)\to k_J.
\]
Let $V:=\Image(e_Y)$ and $W:=\Image(f_X)$. 
Let $h_X:=\dim V$ and $h_Y:=\dim W$. 
Let 
\[
G_X:=\Grass_{h_X}(k_I)\quad\text{and}\quad G_Y:=\Grass_{h_Y}(k_J).
\]
Consider the natural actions of the tori 
$k^*_I$ and $k^*_J$ on $k_I$ and $k_J$, 
and their respective actions on $G_X$ and $G_Y$. Let $\O_V$ 
and $\O_W$ denote the orbits of $V$ and $W$ under these 
actions.

Let $U\subseteq U_\mu$ be an open neighborhood of $\mu$. 
We divide the proof in three steps.

\emph{Step~1}: We'll show that 
the closure $\ol\V_{\mu,X}\subseteq\G_X$ satisfies
\[
\ol\V_{\mu,X}=\bigcup_{\ol\mu\in U}\V_{\ol\mu,X}
\]
if $g_Y>0$. 
In fact, assume that $g_Y>0$. Then $\alpha\neq 0$. 
By Lemma~8.4, the map 
$e_X$ induces a closed embedding $\iota_X\:G_X\to\G_X$ such
that $\iota_X(\O_V) = \V_{\mu,X}$.

For each ordered tripartition $\mathbb I=(I',\ol I,I'')$ of $I$ 
let 
\[
V_{\mathbb I}:=k_{I'}+(k_{\ol I}\cap(V+k_{I''})),
\]
and denote by $\O_{V_{\mathbb I}}$ the orbit 
of $V_{\mathbb I}$ in $G_X$ under the action of $k^*_I$. 
By (4.3.1), all 
the Pl\"ucker coordinates of $V$ in $G_X$ are non--zero. Thus, by Lemma~9.1, 
the closure $\ol{\O}_V\subseteq G_X$ is the union of the orbits 
$\O_{V_{\mathbb I}}$ 
obtained 
from all ordered tripartitions $\mathbb I=(I',\ol I,I'')$ of $I$ such that 
\begin{equation}
|I'|<h_X\leq|I-I''|.
\end{equation}

Since 
$h_X=g_Y-|\alpha|+|I|$, Condition (10.5.1) is equivalent to (10.4.1a).
By Lemma~10.4, each ordered tripartition 
$\mathbb I=(I',\ol I,I'')$ of $I$ satisfying (10.5.1) satisfies (10.4.1c,d) 
for a certain $\ol\mu\in U$. 
Conversely, for each $\ol\mu\in U$
there is an ordered tripartition $\mathbb I=(I',\ol I,I'')$ 
of $I$ satisfying (10.4.1c,d) and (10.5.1).

Now, let $\mathbb I=(I',\ol I,I'')$ be an ordered tripartition of $I$ 
satisfying (10.5.1), and let $\ol\mu\in U$ such that (10.4.1c,d) hold. Let 
\begin{equation}
\ol L_X:=L_X(-\textstyle\sum_{p\in I''}x_p)\quad\text{and}\quad 
\ol L_Y:=L_Y(-\textstyle\sum_{p\in I'}y_p).
\end{equation}
Let $\ol e_X$ and $\ol e_Y$ be the following compositions,
\[
\begin{CD}
\ol e_X\:H^0(\ol L_X) @>>> H^0(L_X) @>e_X>> k_I @>>> k_{\ol I},\\
\ol e_Y\:H^0(\ol L_Y) @>>> H^0(L_Y) @>e_Y>> k_I @>>> k_{\ol I},
\end{CD}
\]
where in each composition the first map is the natural injection, and 
the third map is the natural surjection. Equivalently, $\ol e_X$ and 
$\ol e_Y$ are the evaluation maps corresponding to the isomorphisms 
$\ol\zeta_{X,p}\:\ol L_X(x_p)\to k$ and 
$\ol\zeta_{Y,p}\:\ol L_Y(y_p)\to k$ 
induced by $\zeta_{X,p}$ and $\zeta_{Y,p}$ for 
each $p\in\ol I$. Let $\ol V:=\Image(\ol e_Y)$ 
and $\ol h_X:=h_X-|I'|$. 
Since $\ol h_X>0$ by (10.5.1), 
from (4.3.1) we get $\dim\ol V=\ol h_X$. In addition, 
$V_{\mathbb I}=\ol V+k_{I'}$. 
Let $\ol G_X:=\Grass_{\ol h_X}(k_{\ol I})$, 
and denote by 
$\ol\iota_X\:\ol G_X\to\mathbb G_X$ the closed embedding induced 
by $\ol e_X$, 
as in Lemma~8.4. 
Let $\O_{\ol V}$ denote 
the orbit of $\ol V$ under the natural action 
of $k^*_{\ol I}$. Since $V_{\mathbb I}=\ol V+k_{I'}$, there is a bijection 
$\psi_X\:\O_{\ol V}\to \O_{V_{\mathbb I}}$ taking a 
subspace $H\in \O_{\ol V}$ to $H+k_{I'}$. Note that 
\[
\iota_X(\O_{V_{\mathbb I}})=\iota_X(\psi_X(\O_{\ol V}))=
\ol\iota_X(\O_{\ol V}).
\]

By (10.4.1c,d),
\begin{equation}
\ol L_X=\w_X(\textstyle{\sum}_{p\in\Delta}(1+\alpha_{\ol\mu,p})x_p)
\quad\text{and}\quad
\ol L_Y=\w_Y(\textstyle{\sum}_{p\in I_{\ol\mu}}y_p-
\textstyle{\sum}_{p\in\Delta}\alpha_{\ol\mu,p}y_p).
\end{equation}
So, it follows from Lemma~8.4 that
$\ol\iota_X(\O_{\ol V})=\V_{\ol\mu,X}$, and hence 
$\iota_X(\O_{V_{\mathbb I}})=\V_{\ol\mu,X}$. 
Thus, by Lemma~9.1,
\[
\ol\V_{\mu,X}=\bigcup_{\ol\mu\in U}\V_{\ol\mu,X},
\]
finishing Step~1.

\emph{Step~2}: We'll show that the theorem holds if $g_Xg_Y=0$. In fact, 
assume that $g_X=0$. By Theorem~8.5,
\[
\V_{\mu',Y}=\{H^0(\w_Y(\textstyle{\sum}_{p\in\Delta}y_p))\}
\]
and $\V_{\mu'}=\V_{\mu',X}\times\V_{\mu',Y}$ for every $\mu'\in\Q_\Delta^+$. 
If $g_Y=0$ as well, then $\V_{\mu'}=\V_{\mu}$ for every 
$\mu'\in\Q_\Delta^+$, and 
hence the theorem holds. On the other hand, if 
$g_Y>0$ then 
the closure $\ol\V_\mu\subseteq\G$ satisfies
\[
\ol\V_\mu=\ol \V_{\mu,X}\times\V_{\mu,Y}=
\bigcup_{\ol\mu\in U}\V_{\ol\mu,X}\times\V_{\mu,Y}=
\bigcup_{\ol\mu\in U}\V_{\ol\mu,X}\times\V_{\ol\mu,Y}=
\bigcup_{\ol\mu\in U}\V_{\ol\mu},
\]
where the second equality holds 
by Step~1. Now, the case where $g_Y=0$ is analogous. Step~2 is complete.

\emph{Step~3}: We'll show that the theorem holds if $g_Xg_Y>0$. 
Since $g_X>0$, by Lemma~8.4, the map $f_Y$ induces 
a closed embedding $\iota_Y\: G_Y\to\mathbb G_Y$ 
such that $\iota(\O_W)=\V_{\mu,Y}$. Let
\begin{equation*}
T:=\{(s,t)\in 
k^*_I\times k^*_J
\,|\, s_p^{\wt\beta}=t_p^{\wt\alpha}
\text{ for all }p\in I\cap J\}.
\end{equation*}
Let $\O$ be the orbit of $(V,W)$ under the 
induced action of $T$ on $G_X \times G_Y$. By Lemma~8.4, 
we have $\iota(\O)=\V_\mu$,
where $\iota:=\iota_X\times\iota_Y\:G_X\times G_Y\to\mathbb G$. 

For each pair of 
ordered tripartitions $\mathbb I=(I',\ol I,I'')$ of $I$ and 
$\mathbb J=(J',\ol J,J'')$ of $J$, define 
the subspaces $V_{\mathbb I}\subseteq k_I$ and 
$W_{\mathbb J}\subseteq k_J$ by (9.2.1). By (4.3.1), all the
Pl\"ucker coordinates 
of $V$ in $G_X$ and of $W$ in $G_Y$ are non--zero. Thus, by Lemma~9.2, 
the closure $\ol{\O}\subseteq G_X\times G_Y$ is the union of the 
orbits under $T$ of the 
pairs $(V_{\mathbb I},W_{\mathbb J})$ 
obtained from all ordered tripartitions $\mathbb I=(I',\ol I,I'')$ of $I$ 
and $\mathbb J=(J',\ol J,J'')$ of $J$ 
satisfying (9.2.2,3).

As in Step~1, Condition (9.2.2) is equivalent to Conditions (10.4.1a,b). 
Since (9.2.3) is equal to (10.4.1g,h), by Lemma~10.4, 
for each pair consisting of an ordered tripartition 
$(I',\ol I,I'')$ of $I$ and an ordered tripartition 
$(J',\ol J,J'')$ of $J$ satisfying (9.2.2,3), there is a certain 
$\ol\mu\in U$ satisfying 
(10.4.1c--f). 
Conversely, for each $\ol\mu\in U$ 
there are ordered tripartitions $(I',\ol I,I'')$ of $I$ 
and $(J',\ol J,J'')$ of $J$ 
satisfying (10.4.1a--h).

Now, let $\mathbb I=(I',\ol I,I'')$ be an ordered tripartition of $I$ and 
$\mathbb J=(J',\ol J,J'')$ an ordered tripartition of $J$ satisfying 
(9.2.2,3), and let $\ol\mu\in U$ 
such that (10.4.1c--f) hold. Let $\ol L_X$ and $\ol L_Y$ be given by 
(10.5.2), and put
\begin{equation}
\ol M_X:=M_X(-\textstyle{\sum}_{p\in J'}x_p)\quad\text{and}\quad
\ol M_Y:=M_Y(-\textstyle{\sum}_{p\in J''}y_p).
\end{equation}
Let $\ol e_X\: H^0(\ol L_X)\to k_{\ol I}$ and 
$\ol e_Y\: H^0(\ol L_Y)\to k_{\ol I}$ be the evaluation maps corresponding 
to the isomorphisms $\ol\zeta_{X,p}\:\ol L_X(x_p)\to k$ and 
$\ol\zeta_{Y,p}\:\ol L_Y(y_p)\to k$ induced by $\zeta_{X,p}$ and 
$\zeta_{Y,p}$ for each $p\in\ol I$. Let 
$\ol f_X\: H^0(\ol M_X)\to k_{\ol J}$ and 
$\ol f_Y\: H^0(\ol M_Y)\to k_{\ol J}$ 
be the evaluation maps corresponding to the isomorphisms 
$\ol\xi_{X,p}\:\ol M_X(x_p)\to k$ and 
$\ol\xi_{Y,p}\:\ol M_Y(y_p)\to k$ 
induced by $\xi_{X,p}$ and $\xi_{Y,p}$ for each $p\in\ol J$. 

Let $\ol V:=\Image(\ol e_Y)$ and 
$\ol h_X:=h_X-|I'|$. 
In addition, let $\ol G_X:=\Grass_{\ol h_X}(k_{\ol I})$, and 
denote by $\ol\iota_X\:\ol G_X\to\mathbb G_X$ 
the closed embedding induced by $\ol e_X$, as in Lemma~8.4.
Let $\ol W:=\Image(\ol f_X)$ and $\ol h_Y:=h_Y-|J'|$. Then 
$\dim\ol W=\ol h_Y$ and 
$W_{\mathbb J}=\ol W+k_{J'}$,
as in Step~1. 
Let $\ol G_Y:=\Grass_{\ol h_Y}(k_{\ol J})$, and denote by 
$\ol\iota_Y\:\ol G_Y\to\G_Y$ the closed embedding induced by $\ol f_Y$. 
Let $O$ denote the orbit of $(V_{\mathbb I},W_{\mathbb J})$ in $G$ under 
the action of $T$. 
Let $\O^{-}$ 
denote the orbit of $(\ol V,\ol W)$ under the 
natural action of
\[
\ol T:=\{(s,t)\in k^*_{\ol I}\times k^*_{\ol J}\,|\, 
s_p^{\wt\beta}=t_p^{\wt\alpha}
\text{ for every }p\in\ol I\cap\ol J\}
\]
on $\ol G:=\ol G_X\times\ol G_Y$. Since $V_{\mathbb I}=\ol V+k_{I'}$ and 
$W_{\mathbb J}=\ol W+k_{J'}$, there is a bijection 
$\psi\:\O^{-}\to O$ taking a pair 
$(H_X,H_Y)\in \O^{-}$ to 
$(H_X+k_{I'},H_Y+k_{J'})\in O$. 
Let $\ol\iota := \ol\iota_X\times\ol\iota_Y\:\ol G_X\times\ol G_Y\to\G$.
Note that
\begin{equation}
\iota(O)=\iota(\psi(\O^{-}))=\ol\iota(\O^{-}).
\end{equation}

By (10.4.1c,d), we have (10.5.3). In addition, by (10.4.1e,f),
\[
\ol M_X=\w_X(\textstyle{\sum}_{p\in J_{\ol\mu}}x_p-
\textstyle{\sum}_{p\in\Delta}\beta_{\ol\mu,p}x_p)
\quad\text{and}\quad
\ol M_Y=\w_Y(\textstyle{\sum}_{p\in\Delta}(1+\beta_{\ol\mu,p})y_p).
\]

We will apply Lemma~8.4 to show that $\ol\iota(\O^{-})=\V_{\ol\mu}$. First, 
note that $\alpha_{\ol\mu}\neq 0$. Indeed, $\alpha_{\ol\mu,p}=\alpha_p$ for 
every $p\in\ol I$ by (10.4.1d). Now, $\ol I$ is non--empty and 
contained in $I$, and 
$\alpha_p>0$ for every $p\in I$. Hence $\alpha_{\ol\mu,p}>0$ for every 
$p\in\ol I$. Analogously, $\beta_{\ol\mu}\neq 0$. So, Lemma~8.4 
applies. If $I_{\ol\mu}\cap J_{\ol\mu}=\emptyset$, we get 
$\ol\iota(\O^{-})=\V_{\ol\mu}$ immediately. 

Assume $I_{\ol\mu}\cap J_{\ol\mu}\neq\emptyset$, and let 
$q\in I_{\ol\mu}\cap J_{\ol\mu}$. Then
\[
\wt\alpha_{\ol\mu}=
\frac{\alpha_{\ol\mu,q}}{\gcd(\alpha_{\ol\mu,q},\beta_{\ol\mu,q})}=
\frac{\alpha_q}{\gcd(\alpha_q,\beta_q)}=\wt\alpha,
\]
where the middle equality follows from (10.4.1c--f). Analogously, 
$\wt\beta_{\ol\mu}=\wt\beta$. Since
\[
\{(\zeta_{X,p}^{\otimes\wt\beta}\otimes\xi_{X,p}^{\otimes\wt\alpha},
\zeta_{Y,p}^{\otimes\wt\beta}\otimes\xi_{Y,p}^{\otimes\wt\alpha})\,|\, 
p\in I\cap J\}
\]
patch $L_X^{\otimes\wt\beta}\otimes M_X^{\otimes\wt\alpha}$ and 
$L_Y^{\otimes\wt\beta}\otimes M_Y^{\otimes\wt\alpha}$ 
to the sheaf $K_{I\cap J}$ given by 
(8.2.2), it follows from (10.5.2) and (10.5.4) that 
\[
\{(\ol\zeta_{X,p}^{\otimes\wt\beta_{\ol\mu}}\otimes
\ol\xi_{X,p}^{\otimes\wt\alpha_{\ol\mu}},
\ol\zeta_{Y,p}^{\otimes\wt\beta_{\ol\mu}}\otimes
\ol\xi_{Y,p}^{\otimes\wt\alpha_{\ol\mu}})\,|\, 
p\in I_{\ol\mu}\cap J_{\ol\mu}\}
\]
patch 
$\ol L_X^{\otimes\wt\beta_{\ol\mu}}\otimes
\ol M_X^{\otimes\wt\alpha_{\ol\mu}}$ and 
$\ol L_Y^{\otimes\wt\beta_{\ol\mu}}\otimes
\ol M_Y^{\otimes\wt\alpha_{\ol\mu}}$ 
to the sheaf
\[
\ol K:=\ol K_{I\cap J}(-\wt\beta_{\ol\mu}(\textstyle{\sum}_{p\in I''}x_p+
\textstyle{\sum}_{p\in I'}y_p)
-\wt\alpha_{\ol\mu}(\textstyle{\sum}_{p\in J''}y_p+
\textstyle{\sum}_{p\in J'}x_p)),
\]
where $\ol K_{I\cap J}$ is the pull--back of $K_{I\cap J}$ to the blowup 
$\ol C$ of $C$ along $\Delta-(I_{\ol\mu}\cap J_{\ol\mu})$.
Using Equation (8.2.2) and (10.4.1c--f), we get
\[
\ol K\cong
\ol\w^{\otimes(\wt\alpha_{\ol\mu}+\wt\beta_{\ol\mu})}
\big(\sum_{p\not\in I_{\ol\mu}\cap J_{\ol\mu}}
(\wt\beta_{\ol\mu}\alpha_{\ol\mu,p}-\wt\alpha_{\ol\mu}\beta_{\ol\mu,p})
(x_p-y_p)-\wt\alpha_{\ol\mu}\sum_{p\not\in J_{\ol\mu}}x_p
-\wt\beta_{\ol\mu}\sum_{p\not\in I_{\ol\mu}}y_p\big),
\]
where $\ol\w$ is the pull--back to $\ol C$ of the dualizing 
sheaf $\w$ of $C$. 
Applying Lemma~8.4, we get $\ol\iota(\O^{-})=\V_{\ol\mu}$. 

Since $\ol\iota(\O^{-})=\V_{\ol\mu}$, we get $\iota(O)=\V_{\ol\mu}$ from 
(10.5.5). Thus, by Lemma~9.2,
\[
\ol\V_\mu=\bigcup_{\ol\mu\in U}\V_{\ol\mu},
\]
finishing Step~3.
\end{proof}

\begin{theorem}\setcounter{equation}{0} 
Preserve {\rm 4.1} and assume 
{\rm (4.3.1)}. Then the variety 
$\V$
of limit canonical systems 
of $C$
is 
projective and connected. 
In addition, $\V$ is of pure dimension 
$\delta-1$, unless $g_X=g_Y=0$; in the exceptional case, 
\begin{equation}
\V=\{(H^0(\w_X(\textstyle\sum_{p\in\Delta}x_p)),
H^0(\w_Y(\textstyle\sum_{p\in\Delta}y_p)))\}.
\end{equation}
\end{theorem}
 
\begin{proof} By Proposition~10.3, the covering 
$\{\V_\mu\,\,|\, \mu\in\Q_\Delta^+\}$ of $\V$ is finite. 
Hence, it follows from 
Theorem~10.5 
that $\V$ is projective. 

Let's prove that $\V$ is connected. Let $\mu_1,\mu_2\in\Q_\Delta^+$. 
For each $t\in[0,1]\cap\Q$ let 
\[
\mu(t):=(1-t)\mu_1+t\mu_2\quad\text{and}\quad
I_t:=\{s\in[0,1]\cap\Q\,|\,\V_{\mu(s)}=\V_{\mu(t)}\}
\]
By Proposition~10.3, each $I_t$ is an interval, and finitely many 
$t_1,\dots,t_n\in[0,1]\cap\Q$ suffice to get a covering 
$I_{t_1},\dots,I_{t_n}$ of $[0,1]\cap\Q$. Assume 
$t_1<\dots<t_n$. 
Then $\V_{\mu(t_1)}=\V_{\mu_1}$ and $\V_{\mu(t_n)}=\V_{\mu_2}$. 
By Theorem~10.5, for each $i=1,\dots,n-1$ either 
$\V_{\mu(t_i)}\subseteq\ol\V_{\mu(t_{i+1})}$ or 
$\ol\V_{\mu(t_i)}\supseteq\V_{\mu(t_{i+1})}$; 
in either case
$\ol\V_{\mu(t_i)}\cap\ol\V_{\mu(t_{i+1})}\neq\emptyset$. Since each 
$\V_{\mu(t)}$ is irreducible by Theorem~8.5, it follows that 
$\V_{\mu_1}$ and $\V_{\mu_2}$ lie in the same connected component of 
$\V$. So $\V$ is connected.

If $g_X=g_Y=0$, then $\alpha_\mu=\beta_\mu=0$ for 
each
$\mu\in\Q_\Delta^+$. 
Hence (10.6.1) follows from Theorem~8.5.

Assume now that either $g_X>0$ or $g_Y>0$. Let's show 
that $\V$ is of pure dimension $\delta-1$. Let 
$\mu\in\Z_\Delta^+$. 
Then $\dim\V_\mu\leq\delta-1$ by Theorem~8.5. So, it's enough to show that 
either $\dim\V_\mu=\delta-1$ or there is 
$\mu'\in\Q_\Delta^+$ such that $\V_\mu\subseteq\ol\V_{\mu'}-\V_{\mu'}$. 
To ease notation, let
\[
(\alpha,\rho,I,\beta,\sigma,J,\epsilon):=
(\alpha_\mu,\rho_\mu,I_\mu,\beta_\mu,\sigma_\mu,J_\mu,\epsilon_\mu).
\]
There are two cases to consider.

\emph{Case 1}: Assume that there is $p\in\Delta$ such that
\begin{equation}
p\not\in I\text{ if }|\alpha|>g_Y\text{ and }
p\not\in J\text{ if }|\beta|>g_X.
\end{equation}
Let $p\in\Delta$ satisfying (10.6.2), and put
\[
t':=\min\Big(\frac{\rho_p}{1+\alpha_p},
\frac{\mu_p-\sigma_p}{1+\beta_p}\Big).
\]
Then $t'>0$. For each $t\in\Q$ such that $0\leq t\leq t'$ let 
$\mu(t):=\mu-tp$. We claim that $t'<\mu_p$. Indeed, if 
$p\not\in I\cap J$ then either $\rho_p<\mu_p$ or $\sigma_p>0$, and 
hence $t'<\mu_p$. On the other hand, 
since $g_X>0$ or $g_Y>0$, 
if $p\in I\cap J$
then 
$\alpha_p>0$ or $\beta_p>0$, and hence 
$t'\leq\mu_p/2$. Therefore
$\mu(t)\in\Q_\Delta^+$ if $0\leq t\leq t'$. 

Now, if $0\leq t<t'$ then 
$(\alpha_{\mu(t)},\beta_{\mu(t)})=(\alpha,\beta)$. Moreover, it follows 
from (10.6.2) that $I_{\mu(t)}=I$ if $|\alpha|>g_Y$ and 
$J_{\mu(t)}=J$ if $|\beta|>g_X$. By Theorem~10.1, $\V_{\mu(t)}=\V_\mu$ 
if $0\leq t<t'$. Thus $\V_{\mu}\subseteq\ol\V_{\mu'}$ by Theorem~10.5, 
where $\mu':=\mu(t')$. On the other hand, either 
$\alpha_{\mu',p}=\alpha_p+1$ or $\beta_{\mu',p}=\beta_p+1$. By 
Theorem~10.1, $\V_{\mu'}\cap\V_\mu=\emptyset$. So 
$\V_\mu\subseteq\ol\V_{\mu'}-\V_{\mu'}$, finishing the proof.

\emph{Case 2}: Assume now that there is no $p\in\Delta$ satisfying 
(10.6.2). In addition, assume that $\dim\V_\mu\neq\delta-1$. We claim that 
\begin{equation}
|\alpha|>g_Y,\quad|\beta|>g_X,\quad I\cup J=\Delta \quad \text{and} \quad I\cap J=\emptyset.
\end{equation}
Indeed, suppose that $|\alpha|=g_Y$. Since there is no $p\in\Delta$ 
satisfying (10.6.2), we must have $|\beta|>g_X$. Then $\dim\V_\mu=|J|-1$ by 
Theorem~8.5. So $J\neq\Delta$ because $\dim\V_\mu<\delta-1$. Now, 
(10.6.2) holds for $p\in\Delta-J$, 
reaching
a contradiction. Thus $|\alpha|>g_Y$. 
Analogously, 
$|\beta|>g_X$. Now, 
(10.6.2) holds for 
$p\in\Delta-(I\cup J)$.
Thus $I\cup J=\Delta$. 
Finally, if $I\cap J\neq\emptyset$ then 
$\dim\V_\mu=\delta-1$ by Theorem~8.5. Thus $I\cap J=\emptyset$, and 
(10.6.3) holds.

Define $\lambda\in\Q_\Delta$ by letting $\lambda_p:=\mu_p$ if $p\in J$ and 
$\lambda_p:=0$ otherwise. Put
\[
t':=\min\Big(\min\Big(\frac{\sigma_p}{\epsilon}\,\Big|\, p\not\in J\Big),
\min\Big(\frac{\rho_p}{\mu_p(1+\alpha_p)}\,\Big|\, p\in J\Big)\Big).
\]
Then $t'>0$. Let $\mu(t):=\mu-t\lambda$ for each $t\in\Q$ such that 
$0\leq t\leq t'$. 
Now, $\rho_p<\mu_p$ for every $p\in J$ because $I\cap J=\emptyset$. So 
$t'<1$, and hence $\mu(t)\in\Q_\Delta^+$ if $0\leq t\leq t'$.
Let $\mu':=\mu(t')$. Then either 
$\alpha_{\mu',p}=\alpha_{\mu,p}+1$ for some $p\in J$ or $J_{\mu'}\supset J$. 
Since $|\beta|>g_X$, in either case $\V_\mu \cap \V_{\mu'}=\emptyset$ by 
Theorem~10.1. On the other hand,
\[
(\alpha_{\mu(t)},\beta_{\mu(t)},I_{\mu(t)},J_{\mu(t)})=(\alpha,\beta,I,J)
\]
if $0\leq t<t'$. Hence $\V_\mu\subseteq\ol\V_{\mu'}-\V_{\mu'}$ as in Case 1.
\end{proof}

\begin{theorem}
\setcounter{equation}{0} Preserve {\rm 4.1} and assume 
{\rm (4.3.1)}. For each $\nu=(V_1,V_2)\in\mathbb G$, let 
\[
W_\nu:=R_{\nu,X}+R_{\nu,Y}+\sum_{p\in\Delta}g(\delta-2)p,
\]
where $R_{\nu,X}$ and $R_{\nu,Y}$ are the ramification divisors of the 
linear systems $(V_1,\mathbb L)$ and $(V_2,\mathbb M)$, respectively. 
If $\V$ is the variety of limit canonical systems of $C$, then
$\{W_\nu\,|\, \nu\in\V\}$ is the set of limit Weierstrass divisors of 
smoothings of $C$.
\end{theorem}

\begin{proof} Let $\mu\in\Z_\Delta^+$, and denote by $\wt C$ the 
curve gotten from $C$ by splitting 
the branches of $C$ at each $p\in\Delta$ 
and connecting them by a chain of $\mu_p-1$ rational smooth 
curves; see~4.2.
Let
\begin{equation}
L_X:=\w_X(\textstyle{\sum}_{p\in\Delta}(1+\alpha_{\mu,p})x_p)
\quad\text{and}\quad
M_Y:=\w_Y(\textstyle{\sum}_{p\in\Delta}(1+\beta_{\mu,p})y_p).
\end{equation}
Let $\nu:=(V_1,V_2)\in\V_\mu$. By definition, $(V_1,L_X)$ and $(V_2,M_Y)$ 
are the limit canonical aspects with foci on $X$ and $Y$ of a 
certain regular smoothing 
$\wt\pi$ of $\wt C$. Let $\pi$ be the induced smoothing of 
$C$; see~2.7.
By Theorem~5.2, the limit Weierstrass scheme $W$ of $\pi$ satisfies
\begin{equation}
[W]=R_1+R_2+\sum_{p\in\Delta}g(g-1-\alpha_{\mu,p}-\beta_{\mu,p})p,
\end{equation}
where $R_1$ and $R_2$ are the ramification divisors of $(V_1,L_X)$ and 
$(V_2,M_Y)$. Now, since 
\[
\mathbb L=L_X(\textstyle{\sum}_{p\in\Delta}(g_Y-\alpha_{\mu,p})x_p)
\quad\text{and}\quad
\mathbb M=M_Y(\textstyle{\sum}_{p\in\Delta}(g_X-\beta_{\mu,p})y_p),
\]
we have
\[
R_1=R_{\nu,X}-g\textstyle{\sum}_{p\in\Delta}(g_Y-\alpha_{\mu,p})p
\quad\text{and}\quad
R_2=R_{\nu,Y}-g\textstyle{\sum}_{p\in\Delta}(g_X-\beta_{\mu,p})p.
\]
Using the above expressions in (10.7.2), since $g=g_X+g_Y+\delta-1$, 
we get $[W]=W_\nu$. So $W_\nu$ is a limit Weierstrass divisor.

Conversely, if $\pi$ is a smoothing of $C$, then there is 
$\mu\in\Z_\Delta^+$ such that $\pi$ induces a regular 
smoothing $\wt\pi$ of the curve $\wt C$ obtained from $C$ by splitting 
the branches of $C$ at each $p\in\Delta$ 
and connecting them by a chain of $\mu_p-1$ rational smooth curves; 
see~2.7 and~4.2. Let $L_X$ and $M_Y$ be as in (10.7.1). By Theorem~5.2 
there are vector subspaces $V_1\subseteq H^0(L_X)$ and 
$V_2\subseteq H^0(M_Y)$ such that $(V_1,L_X)$ and $(V_2,M_Y)$ are 
the limit canonical aspects of $\wt\pi$ 
with foci on $X$ and $Y$, respectively. Let $\nu:=(V_1,V_2)\in\G$. 
By definition, $\nu\in\V_\mu$. Now, proceed as in the above paragraph 
to conclude that the limit Weierstrass scheme $W$ of $\pi$ satisfies 
$[W]=W_\nu$.
\end{proof}

\section{Miscellany}

\begin{proposition}\setcounter{equation}{0} Preserve {\rm 4.1} and assume 
{\rm (4.3.1)}. If $g_X>0$ or $g_Y>0$ then there is a 
$(\delta-1)$--dimensional family of limit Weierstrass divisors 
on 
$C$.
\end{proposition}

\begin{proof} Assume that $g_Y>0$. Then we may choose $\alpha\in\Z_\Delta^+$ 
such that $|\alpha|=g_Y+\delta-1$. Define $\mu\in\Q_\Delta^+$ by letting 
$\mu_p:=1/\alpha_p$ for each $p\in\Delta$. 
Then $\alpha_\mu=\alpha$ and $I_\mu=\Delta$. 
Let 
\[
L_X:=\w_X(\textstyle{\sum}_{p\in\Delta}(1+\alpha_p)x_p)\quad\text{and}\quad
\ol L_X:=L_X(-\textstyle{\sum}_{p\in\Delta}x_p).
\]
Then $h^0(L_X)=g+\delta-1$ and $h^0(\ol L_X)=g-1$. 
For each $q\in X$ let 
\[
V_q:=H^0(L_X(-gq))+H^0(\ol L_X)\subseteq H^0(L_X).
\]
Let $U\subseteq X$ be the dense open subset of points $q\in X$ such
that $h^0(L_X(-gq))=\delta-1$ and $h^0(\ol L_X(-gq))=0$. Then 
$\dim V_q=h^0(L_X)-1$ for each $q\in U$.

We claim that there are $q_1,\dots,q_{\delta-1}\in U$ such that 
\begin{subequations}
\begin{align}
&\dim(V_{q_1}\cap\cdots\cap V_{q_j})=h^0(L_X)-j,\\
&V_{q_1}\cap\cdots\cap V_{q_j}\not\subseteq H^0(L_X(-x_p))
\text{ for any $p\in\Delta$,}
\end{align}
\end{subequations}
for $j=0,\dots,\delta-1$. Indeed, 
let's prove the claim 
by induction on $j$. 
First, (11.1.1a,b) hold trivially for $j=0$.
Now, 
let $j\in\{0,\dots,\delta-2\}$ and 
suppose that there are 
$q_1,\dots,q_j\in U$ 
such that (11.1.1a,b) 
hold. 
Let $V_j:=V_{q_1}\cap\cdots\cap V_{q_j}$. 
Then there is 
$q\in U$ 
such that
\begin{equation}
\dim(V_j\cap H^0(L_X(-x_p-gq)))=\delta-j-2\quad\text{for every $p\in\Delta$}.
\end{equation}
By (11.1.1a) and (11.1.2), 
\[
h^0(L_X(-x_p-gq))=\delta-2\quad\text{and}\quad V_j+H^0(L_X(-x_p-gq))=H^0(L_X)
\]
for each $p\in\Delta$.
Now, if $V_j\subseteq V_q$ then $V_j+H^0(L_X(-gq))\subseteq V_q$, and 
hence 
$V_q=H^0(L_X)$. However, $\dim V_q=h^0(L_X)-1$ because $q\in U$, reaching 
a contradiction. So 
$V_j\not\subseteq V_q$, 
and hence $\dim(V_j\cap V_q)=h^0(L_X)-j-1$.

In addition,
if $V_j\cap V_q\subseteq H^0(L_X(-x_p))$ for a certain 
$p\in\Delta$, then 
\[
V_j\cap V_q=V_j\cap (H^0(L_X(-x_p-gq))+H^0(\ol L_X)).
\]
Now, $h^0(\ol L_X(-gq))=0$ because $q\in U$. Moreover, 
$V_j\supseteq H^0(\ol L_X)$. So
\[
\dim(V_j\cap H^0(L_X(-x_p-gq)))=\dim(V_j\cap V_q)-h^0(\ol L_X).
\]
Since $\dim(V_j\cap V_q)=h^0(L_X)-j-1$, it follows that
\[
\dim(V_j\cap H^0(L_X(-x_p-gq)))=h^0(L_X)-j-1-(g-1)=\delta-j-1,
\]
contradicting (11.1.2). Hence $V_j\cap V_q\not\subseteq H^0(L_X(-x_p))$ 
for any $p\in\Delta$. The induction proof 
of our claim
is complete.

For each $D=(q_1,\dots,q_{\delta-1})\in U^{\delta-1}$ let 
$V_D:=V_{q_1}\cap\cdots\cap V_{q_{\delta-1}}$. The following conditions on 
$D\in U^{\delta-1}$,
\begin{subequations}
\begin{align}
&\dim V_D=g,\\
&V_D\not\subseteq H^0(L_X(-x_p))\text{ for any $p\in\Delta$,}
\end{align}
\end{subequations}
define an open subset $W\subseteq U^{\delta-1}$. By our claim, 
$W\neq\emptyset$.

Let $D:=(q_1,\dots,q_{\delta-1})\in W$. Since
$V_{q_j}\supseteq H^0(\ol L_X)$ 
for each $j=1,\dots,\delta-1$, also $V_D\supseteq H^0(\ol L_X)$.  
As shown in~8.6, Conditions~(11.1.3) imply
that $V_D\in\V_{\mu,X}$. 
Since $V_D+H^0(L_X(-gq_j))\subseteq V_{q_j}$, and since 
$\dim V_{q_j}=g+\delta-2$ and $h^0(L_X(-gq_j))=\delta-1$, 
it follows that $V_D\cap H^0(L_X(-gq_j))\neq 0$ for each 
$j=1,\dots,\delta-1$. 
So $q_1,\dots,q_{\delta-1}$ are ramification points 
of $(V_D,L_X)$, and hence of $(V_D,\mathbb L)$. Since $V_D\in\V_{\mu,X}$ 
it follows from Theorem~10.7 that $q_1,\dots,q_{\delta-1}$ sit on the 
support of a limit Weierstrass divisor 
on
$C$. Since $\dim W=\delta-1$, 
there must be a $(\delta-1)$--dimensional 
family of limit Weierstrass divisors.
\end{proof}

\begin{theorem}\setcounter{equation}{0} Preserve {\rm 4.1} and assume 
{\rm (4.3.1)}. Let $\V$ be the variety of limit canonical systems of $C$, 
and denote by $N(\V)$ the number of its irreducible components. 
Set 
\[
n_\delta(h):=\binom{h+\delta-1}{\delta}-\binom{h}{\delta}
\quad\text{and}\quad
g_{i,j}:=\gcd(g_X+i,g_Y+j)
\]
for 
all
$h,i,j\in\Z^+$. Assume that $\delta>1$. 
Then the following two statements hold.
\begin{enumerate}
\item If $g_Xg_Y=0$ or $g_X=g_Y$ then $N(\V)=n_\delta(\max(g_X,g_Y))$, 
unless $g_X=g_Y=0$; in the exceptional case, $N(\V)=1$.
\item If $g_Xg_Y>0$
then 
\[
N(\V)\geq n_\delta(g_X)+n_\delta(g_Y)-\sum_{i,j=1}^{\delta-1}
\binom{g_{i,j}-1}{\delta-1}.
\]
\item $\V$ is irreducible if and only if $g_X\leq 1$ and $g_Y\leq 1$.
\end{enumerate}
\end{theorem}

\begin{proof} If $g_X=g_Y=0$ then $\V$ is simply a point by 
Theorem~10.6, thus proving Statements~1~and~3 in this case. 
We may assume from now on that 
$g_Y>0$. (The case where $g_X>0$ is completely analogous.) 

Given $i\in\Z^+$, the number of $\alpha\in\Z^+_\Delta$ such
$|\alpha|=i$ is $\binom{i-1}{\delta-1}$. Thus the number 
of $\alpha\in\Z^+_\Delta$ such that $i<|\alpha|<i+\delta$ is 
$n_\delta(i)$.

If $h\in\Z^+$ then $n_\delta(h)=1$ if and only if $h=1$. Indeed, 
$n_\delta(1)=1$ because $\delta>1$. In addition, if $h>1$ then 
the number of $\alpha\in\Z^+_\Delta$ such that $|\alpha|=h+\delta-1$ is 
at least $\delta$, and hence $n_\delta(h)>1$ because $\delta>1$.

Let $N_X(\V)$ be the number of strata $\V_\mu$ of $\V$ 
with $\dim\V_{\mu,X}=\delta-1$. By Theorem~8.5, $\dim\V_{\mu,X}=\delta-1$ 
if and only if $I_\mu=\Delta$ and $|\alpha_\mu|>g_Y$. Now, 
given $\alpha\in\Z^+_\Delta$ satisfying $g_Y<|\alpha|<g_Y+\delta$, there is 
$\mu\in\Q_\Delta^+$ such that $I_\mu=\Delta$ and $\alpha_\mu=\alpha$. 
(For instance, pick $\mu\in\Q_\Delta^+$ given by $\mu_p:=1/\alpha_p$ 
for each $p\in\Delta$.) If also $\mu'\in\Q_\Delta^+$ satisfies 
$I_{\mu'}=\Delta$ and 
$\alpha_{\mu'}=\alpha$ then there is $t\in\Q^+$ such that 
$\mu'=t\mu$, and hence $\V_{\mu'}=\V_\mu$ by Theorem~10.1. It follows that
sending $\V_\mu$ to $\alpha_\mu$ gives a 1--1 correspondence between the set 
of strata $\V_\mu$ of $\V$ with $\dim\V_{\mu,X}=\delta-1$ and the set of 
$\alpha\in\Z^+_\Delta$ with $g_Y<|\alpha|<g_Y+\delta$. Thus 
$N_X(\V)=n_\delta(g_Y)$.

If $g_X=0$ then $N(\V)=N_X(\V)$, thus proving Statement~1 in this case. 
Moreover, $N(\V)=1$ if and only if $g_Y=1$, thus proving 
Statement~3 in this case. 

Assume from now on that $g_X>0$. By analogy, 
$N_Y(\V)=n_\delta(g_X)$, where $N_Y(\V)$ is the number of
strata 
$\V_\mu$ of $\V$ with $\dim\V_{\mu,Y}=\delta-1$. 

Fix $i,j\in\{1,\dots,\delta-1\}$. Let 
$\alpha,\beta\in\Z_\Delta^+$ such that
\begin{equation}
|\alpha|=g_Y+i\quad\text{and}\quad|\beta|=g_X+j.
\end{equation}
Then there is $\mu\in\Q_\Delta^+$ such that 
\begin{equation}
\mu_p\alpha_p=\mu_q\alpha_q\quad\text{and}\quad\mu_p\beta_p=\mu_q\beta_q
\quad\text{for all $p,q\in\Delta$}
\end{equation}
if and only if
\begin{equation} 
\frac{\alpha_p}{\beta_p}=\frac{g_Y+i}{g_X+j}\quad
\text{ for every $p\in\Delta$.}
\end{equation}
Now, $\alpha,\beta\in\Z_\Delta^+$ satisfy (11.2.1,3) if and only 
if there is $\tau\in\Z_\Delta^+$ such that
\[
|\tau|=g_{i,j},\quad\alpha=\frac{g_Y+i}{g_{i,j}}\tau,\quad
\beta=\frac{g_X+j}{g_{i,j}}\tau.
\]
So,
there is a 1--1 correspondence between the set of pairs 
$(\alpha,\beta)\in\Z_\Delta^+\times\Z_\Delta^+$ satisfying (11.2.1) and 
(11.2.2) for a certain $\mu\in\Q_\Delta^+$ and the set of 
$\tau\in\Z_\Delta^+$ such that $|\tau|=g_{i,j}$. The number of 
$\tau\in\Z_\Delta^+$ such that $|\tau|=g_{i,j}$ is 
$\binom{g_{i,j}-1}{\delta-1}$. 

Let $N$ be the number of strata $\V_\mu$ of $\V$ such that 
\begin{equation}
\dim\V_{\mu,X}=\dim\V_{\mu,Y}=\delta-1.
\end{equation}
Then $N(\V)\geq N_X(\V)+N_Y(\V)-N$. Now, by Theorem~8.5, Condition (11.2.4) 
holds if and only if
\[
I_\mu=J_\mu=\Delta,\quad|\alpha_\mu|>g_Y,\quad|\beta_\mu|>g_X.
\]
It follows from Theorem~10.1 that 
sending $\V_\mu$ to $(\alpha_\mu,\beta_\mu)$ gives a 1--1 correspondence 
between the set of strata $\V_\mu$ of $\V$ satisfying (11.2.4) and the set 
of pairs $(\alpha,\beta)\in\Z_\Delta^+\times\Z_\Delta^+$ satisfying 
(11.2.1) for certain $i,j\in\{1,\dots,\delta-1\}$ and (11.2.2) for 
a certain $\mu\in\Q_\Delta^+$. Thus 
\[
N=\sum_{i,j=1}^{\delta-1}\binom{g_{i,j}-1}{\delta-1}.
\]
Since $N(\V)\geq N_X(\V)+N_Y(\V)-N$, we proved Statement~2.

Without loss of generality, assume from now on that $g_Y\geq g_X$. 
Clearly $N(\V)\geq N_X(\V)$, and hence 
$N(\V)\geq n_\delta(g_Y)$. So, 
$N(\V)=1$ only if $n_\delta(g_Y)=1$, and hence only if $g_Y=1$. 
It follows that $\V$ is irreducible only if $g_X=g_Y=1$.

Assume from now on that $g_X=g_Y$. Then 
$(\alpha_\mu,I_\mu)=(\beta_\mu,J_\mu)$ for each $\mu\in\Q_\Delta^+$. It 
follows from Theorem~8.5 that $\dim\V_{\mu,X}=\delta-1$ if and only if 
$\dim\V_\mu=\delta-1$. Hence $N(\V)=N_X(\V)$, thus proving 
Statement~1. Now, if $g_X=g_Y=1$ then $N(\V)=1$ by Statement~1.
\end{proof}

\begin{example}\setcounter{equation}{0} Preserve {\rm 4.1} and assume 
{\rm (4.3.1)}. 
Assume that $g_X=g_Y=1$. 
Then $g=\delta+1$. 
Choose a basis for $H^0(\w^{\otimes 2})$ and consider the corresponding 
bicanonical map
$\phi\:C\to\mathbb P^{3\delta-1}$. 
Let $H$ (resp. $H_X$, resp. $H_Y$) be the 
subspace
of 
$\mathbb P^{3\delta-1}$ 
spanned by 
the image of 
$\Delta$ (resp. $X$, resp. 
$Y$) under $\phi$.
Then $\dim H=\delta-1$ and 
\[
\dim H_X=\dim H_Y=2\delta-1.
\]
Moreover, $H_X\cap H_Y=H$. We may view each hyperplane $L\subseteq H$ as a 
subspace of both 
$H_X$ and $H_Y$,
and 
compose
$\phi|_X\:X\to H_X$ and $\phi|_Y\:Y\to H_Y$ with the 
projections from $H_X$ and $H_Y$
centered at $L$. Since 
$\w_X$ and $\w_Y$ are trivial, these 
compositions
correspond to 
vector subspaces $V_{L,X}\subseteq H^0(\w_X(2\Delta))$ and 
$V_{L,Y}\subseteq H^0(\w_Y(2\Delta))$ of dimension $g$. 
It follows from the discussion in~8.7 that
\[
\V=\{(V_{L,X},V_{L,Y})\,|\, L\text{ is a hyperplane of }H\}.
\]
\end{example}

\begin{subsct}\setcounter{equation}{0}
\emph{Graphic representation.} Preserve {\rm 4.1} and assume 
{\rm (4.3.1)}. For each $\lambda\in\Z_\Delta^+$ and each non--empty 
subset $K\subseteq\Delta$ let
\begin{align*}
\mathbb U_{X,\lambda,K}:=
&\{\mu\in\Q_\Delta^+\,|\, \alpha_\mu=\lambda\text{ and }I_\mu=K\},\\
\mathbb U_{Y,\lambda,K}:=
&\{\mu\in\Q_\Delta^+\,|\, \beta_\mu=\lambda\text{ and }J_\mu=K\}.
\end{align*}
Then $\mathbb U_{X,\lambda,K}$ and $\mathbb U_{Y,\lambda,K}$ 
are homogeneous.
Now,
$\mu\in\mathbb U_{X,\lambda,K}$ 
if and only if Condition~(14.4.1) below holds.
\begin{equation}
\text{There's $\rho\in\Q_\Delta^+$ such that } 
\begin{cases}
0<\rho_p\leq\mu_p \text{ for every $p\in\Delta$,}&\\
\rho_p=\mu_p \text{ if and only if $p\in K$,}&\\
\mu_p(\lambda_p+1)-\rho_p=\mu_q(\lambda_q+1)-\rho_q
\text{ for all $p,q\in\Delta$.}&
\end{cases}
\end{equation}
Condition (11.4.1) can be restated as
\begin{equation}
\frac{\mu_p}{\mu_q}\in 
\begin{cases}
\left(\frac{\lambda_q}{\lambda_p+1},\frac{\lambda_q}{\lambda_p}\right]
&\text{if $q\in K$,}\\
\left[\frac{\lambda_q}{\lambda_p},\frac{\lambda_q+1}{\lambda_p}\right)
&\text{if $p\in K$,}\\
\left(\frac{\lambda_q}{\lambda_p+1},\frac{\lambda_q+1}{\lambda_p}\right)
&\text{if $p,q\not\in K$.}
\end{cases}
\end{equation}
Similar conditions define the points $\mu\in\mathbb U_{Y,\lambda,K}$. 
Condition~(11.4.2) shows 
that $\mathbb U_{X,\lambda,K}$ is convex. 
Analogously, $\mathbb U_{Y,\lambda,K}$ is convex as well. 

For each $\mu\in\Q_\Delta^+$ let
\begin{align*}
\mathbb U_{X,\mu}:=&\{\mu'\in\Q_\Delta^+\,|\, \V_{\mu',X}=\V_{\mu,X}\},\\
\mathbb U_{Y,\mu}:=&\{\mu'\in\Q_\Delta^+\,|\, \V_{\mu',Y}=\V_{\mu,Y}\},\\
\mathbb U_\mu:=&\{\mu'\in\Q_\Delta^+\,|\, \V_{\mu'}=\V_\mu\}.
\end{align*} 
Then $\mathbb U_\mu=\mathbb U_{X,\mu}\cap\mathbb U_{Y,\mu}$ by 
Theorem~10.1. In addition, 
\begin{subequations}
\begin{align}
&\mathbb U_{X,\mu}=
\begin{cases}
\mathbb U_{X,\alpha_\mu,I_\mu}
&\text{if $|\alpha_\mu|>g_Y$,}\\
\textstyle{\bigcup}_{I\subseteq\Delta}
\mathbb U_{X,\alpha_\mu,I}
&\text{if $|\alpha_\mu|=g_Y$,}
\end{cases}\\
&\mathbb U_{Y,\mu}=
\begin{cases}
\mathbb U_{Y,\beta_\mu,J_\mu}
&\text{if $|\beta_\mu|>g_X$,}\\
\textstyle{\bigcup}_{J\subseteq\Delta}
\mathbb U_{Y,\beta_\mu,J}
&\text{if $|\beta_\mu|=g_X$.}
\end{cases}
\end{align}
\end{subequations}
In particular, $\mathbb U_\mu$, $\mathbb U_{X,\mu}$ and 
$\mathbb U_{Y,\mu}$ are convex and homogeneous.

Using (14.4.2,3) we can represent graphically 
the covering of $\Q_\Delta^+$ by the 
strata $\mathbb U_\mu$, and then study their interrelations. 
Since there is a 1--1 correspondence between the $\mathbb U_\mu$ and the 
$\V_\mu$ this graphic representation is useful in the study of 
the variety of limit canonical systems $\V$. 
For instance, using Theorem~10.5 
and the graphic representation
we can determine whether two irreducible components of $\V$ intersect, and 
what the intersection is. We can also compute the number of irreducible 
components of 
$\V$ using Theorem~8.5 and the graphic representation.
See 11.5 and 11.6 below.

As the $\mathbb U_\mu$ are homogeneous, it's more 
convenient to fix $q\in\Delta$ and consider the map 
\begin{equation}
\Q_\Delta^+\to\Q_{\Delta-\{q\}}^+,\text{ given by $\mu\mapsto\ol\mu$, 
where $\ol\mu_p:=\mu_p/\mu_q$ for 
each
$p\in\Delta-\{q\}$.}
\end{equation}
So, instead of depicting the covering 
of $\Q_\Delta^+$ by the $\mathbb U_\mu$, we can depict the covering of 
$\Q_{\Delta-\{q\}}^+$ by the corresponding images of the $\mathbb U_\mu$. 
\end{subsct}

\begin{subsct}\setcounter{equation}{0}
\emph{Curves with $\delta=2$.} Preserve {\rm 4.1} and assume 
{\rm (4.3.1)}. 
Assume $\delta=2$. 
By Theorem~8.5, a stratum 
$\V_\mu$ of the variety of limit canonical 
systems
$\V$ has dimension 1 if and only if either 
$|\alpha_\mu|=g_Y+1$ or $|\beta_\mu|=g_X+1$. 
Fix $q\in\Delta$, and identify $\Q^+_{\Delta-\{q\}}$ with $\Q^+$. 
Under the map (11.4.4), the 
strata $\mathbb U_\mu$ are sent to the intervals in $\Q^+$ depicted in 
Figure~3, where we assume that $g_Y>g_X>0$.
\setlength{\unitlength}{0.240900pt}
\begin{center}
{
\footnotesize
\begin{picture}(1450,250)(-50,0)
\put(1230,100){\raisebox{-.10pt}{\makebox(0,0){$\times$}}}
\put(1230,130){\makebox(20,20)[c]{$g^{\,}_Y$}}
\put(1080,100){\raisebox{-.10pt}{\makebox(0,0){$\ast$}}}
\put(1070,40){\makebox(20,45)[c]{$g^{\,}_X$}}
\put(920,100){\raisebox{-.10pt}{\makebox(0,0){$\times$}}}
\put(920,130){\makebox(0,50)[c]{$\frac{g^{\,}_Y-1}{2}$}}
\put(480,100){\raisebox{-.10pt}{\makebox(0,0){$\ast$}}}
\put(470,40){\makebox(15,20)[c]{$\frac{2}{g^{\,}_{\!X}-1}$}}
\put(320,100){\raisebox{-.10pt}{\makebox(0,0){$\times$}}}
\put(320,130){\makebox(0,50)[c]{$\frac{2}{g^{\,}_Y-1}$}}
\put(180,100){\raisebox{-.10pt}{\makebox(0,0){$\ast$}}}
\put(170,40){\makebox(15,20)[c]{$\frac{1}{g^{\,}_{\!X}}$}}
\put(97,100){\raisebox{-.10pt}{\makebox(0,0)[c]{$\times$}}}
\put(97,130){\makebox(0,50)[c]{$\frac{1}{g^{\,}_Y}$}}
\put(0,90){\line(0,1){20}}
\put(-10,130){0}
\put(0,100){\line(1,0){650}}
\multiput(685,100)(16,0){5}{\line(1,0){6}}
\put(800,100){\line(1,0){500}}
\multiput(1335,100)(16,0){5}{\line(1,0){6}}
\end{picture}
}
\ \\
\text{\textsc{Figure~3. }Graphic representation for $\delta=2$.}
\addtocounter{figure}{1}
\end{center}
\vskip0.2cm
The points marked as ``$\times$'' (resp. ``$\ast$'') correspond 
to the strata $\V_\mu$ of $\V$ with 
$|\alpha_\mu|=g_Y+1$ (resp. $|\beta_\mu|=g_X+1$).
The number of distinct marked points on the line in Figure~3 is
the number $N(\V)$ of irreducible components of $\V$. Thus
\[
N(\V)=g_X+g_Y-\gcd(g_X+1,g_Y+1)+1.
\]
So,
if $\delta=2$ then equality 
holds in Statement~2 of Theorem~11.2.

The open intervals on the line in Figure~3 correspond 
to the strata $\V_\mu$ of $\V$ with $|\alpha_\mu|=g_Y$ 
and $|\beta_\mu|=g_X$, and hence represent points of $\mathbb V$. 
By Theorem~10.5, the stratum of $\V$ corresponding to a certain interval in 
Figure~3 is the intersection of the closures in $\mathbb G$ of the 
one--dimensional strata of $\V$ corresponding to each end of 
the interval. Therefore, $\V$ is as described in Figure~4.
Each irreducible component of $\V$ is a projective line, unless it 
corresponds to a point on the line in 
Figure~3 that is simultaneously marked as ``$\times$'' and ``$\ast$''. 
In this case, the irreducible component is isomorphic to a rational 
curve of bidegree 
\[
\frac{(g_X+1,g_Y+1)}{\gcd(g_X+1,g_Y+1)}
\]
in the quadric $\mathbb P^1\times\mathbb P^1$.
\begin{figure}[ht]
\begin{center}
\includegraphics[angle=0]{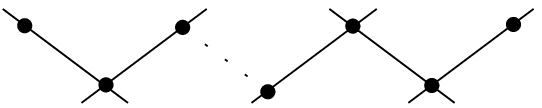}
\caption{The variety of limit canonical systems for $\delta=2$.}
\end{center}
\end{figure}
\end{subsct}

\begin{subsct}\setcounter{equation}{0}
\emph{Curves with $\delta=3$.} Preserve {\rm 4.1} and assume 
{\rm (4.3.1)}. 
Assume $\delta=3$. Fix $r\in\Delta$, and identify $\Q^+_{\Delta-\{r\}}$ 
with $\Q^+\times\Q^+$.
Under the map (11.4.4), the 
strata $\mathbb U_\mu$ are sent to convex regions in $\Q^+\times\Q^+$ 
of dimensions 
varying from $0$ to $2$. These dimensions are ``dual'' to those 
of the corresponding strata 
$\V_\mu$
of the variety of limit canonical systems $\V$. 
More precisely, the image of 
$\mathbb U_\mu$ in $\Q^+\times\Q^+$ has dimension $i$ if and only if 
$\dim\V_\mu=2-i$. 

We depict the decomposition of $\Q^+\times\Q^+$ for 
$g_X=g_Y=3$ in Figure~5a, and for $g_X=2$ and $g_Y=4$ in Figure~5b. In 
Figure~5b the solid lines depict the decomposition of $\Q^+\times\Q^+$ by 
the images under (11.4.4) of the strata $\mathbb U_{\mu,X}$, whereas 
the dashed lines depict the decomposition given by the strata 
$\mathbb U_{\mu,Y}$.
\begin{figure}[ht]
\begin{center}
\vspace{0.2in}
\mbox{
\includegraphics[width=2.2in, viewport= 87 162 522 562]{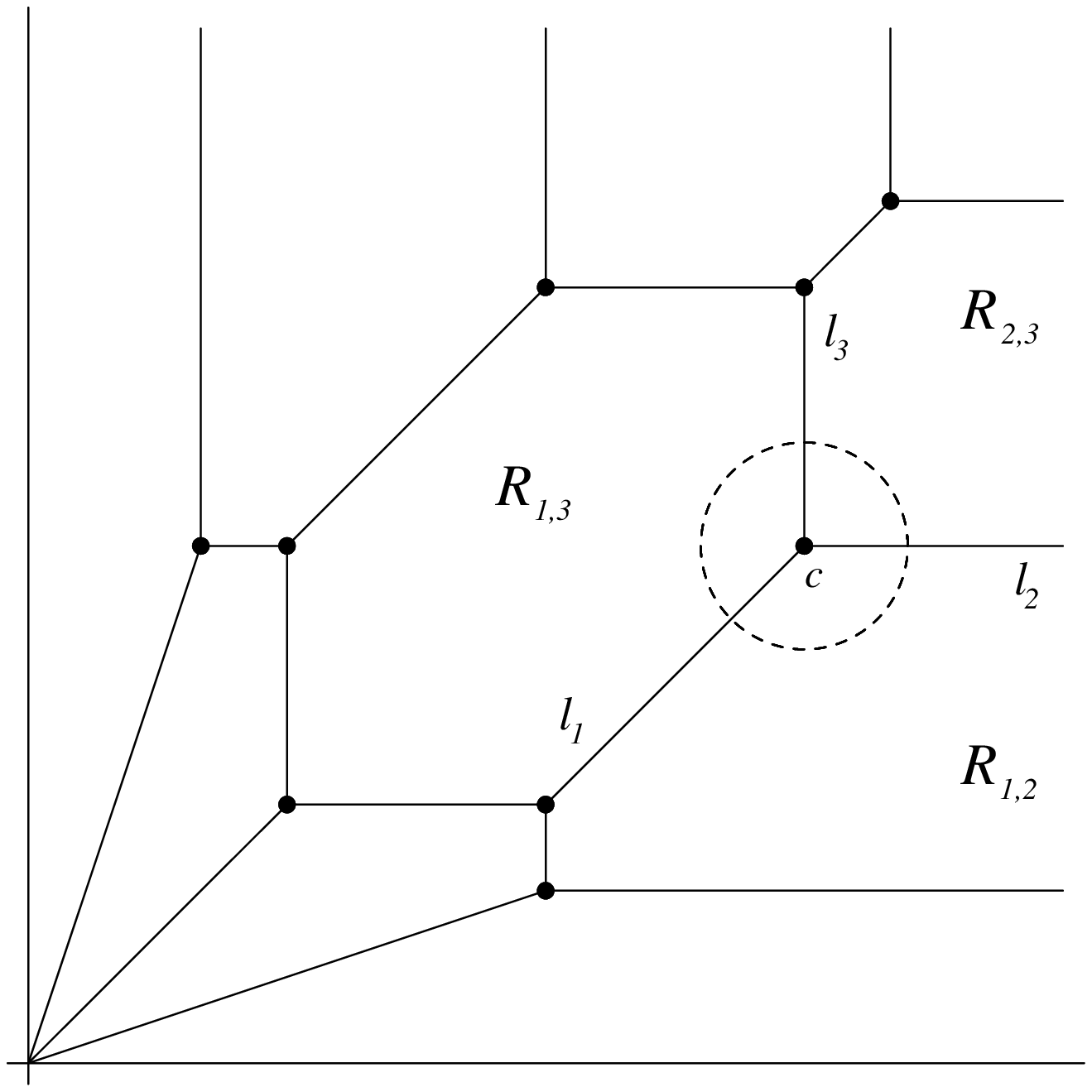}
\mbox{\hspace{1.5cm}}
\includegraphics[width=2.2in, viewport= 87 162 522 562]{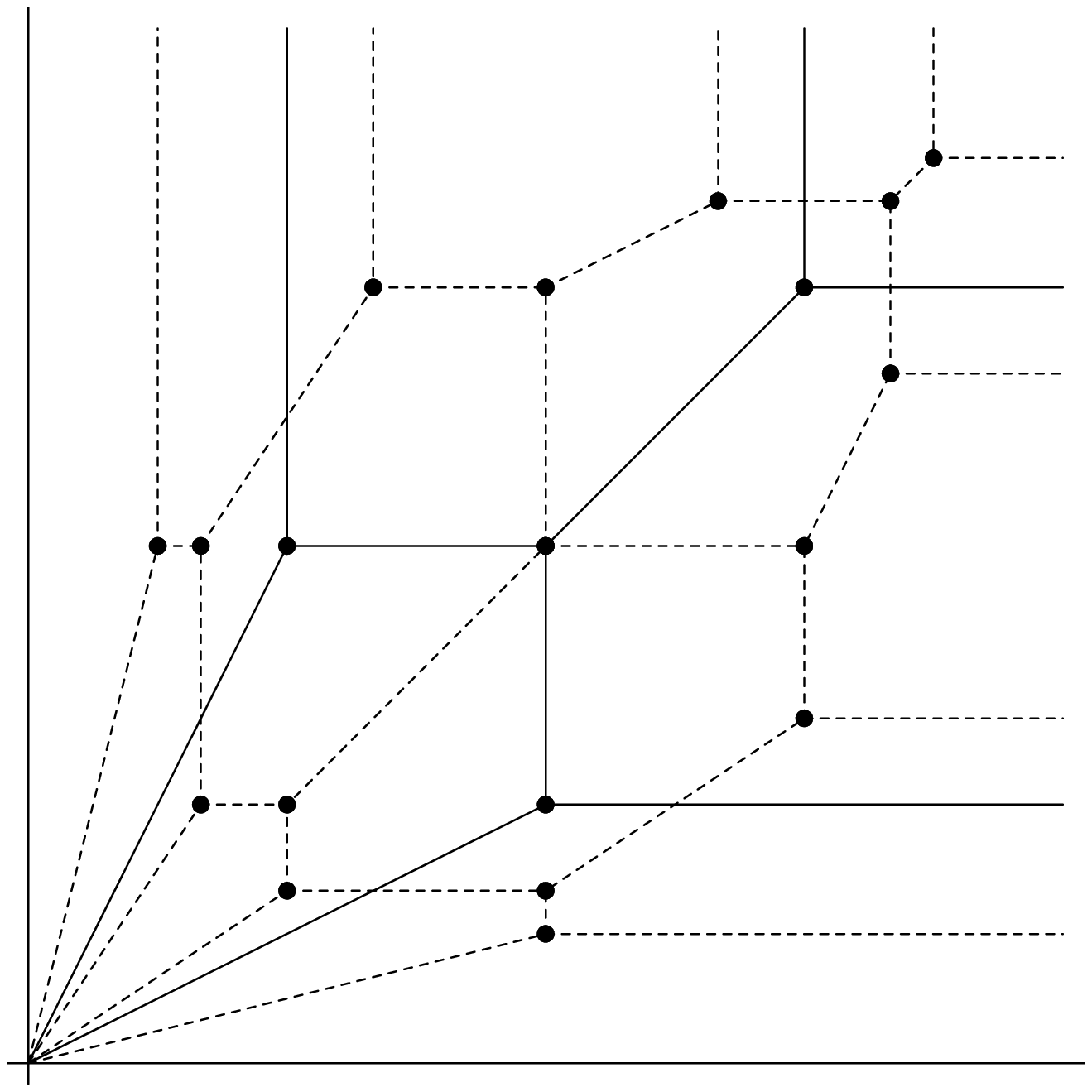}
}
\\[-6pt]
\mbox{\hspace{0.65in}}\ \ \textsc{Figure}~5a: $g_X=g_Y=3$.
\hfill{\textsc{Figure}~5b: $g_X=2\text{ and }g_Y=4$.\mbox{\hspace{0.5in}}}
\end{center}
\end{figure}

The open disc in Figure~5a corresponds to 
one of the neighborhoods $U_\mu$ whose existence 
was claimed by Theorem~10.5. More precisely, let $\mu\in\Q_\Delta^+$ 
whose image under the map (11.4.4) is the center $c$ of the disc, and 
let $U_\mu\subseteq\Q_\Delta^+$ be 
the 
inverse image of the disc under this same map. Then
the closure $\ol\V_\mu\subseteq\G$ satisfies
\[
\ol\V_\mu=\bigcup_{\ol\mu\in U_\mu}\V_{\ol\mu}.
\]
More precisely, $\ol\V_\mu$ is an irreducible component of $\V$ and 
\[
\ol\V_\mu-\V_\mu=\ol\V_{\mu_1}\cup\ol\V_{\mu_2}\cup\ol\V_{\mu_3},
\]
where $\mu_1,\mu_2,\mu_3\in\Q_\Delta^+$ are (any) points having images under 
(11.4.4) in the interiors of the three segments of lines 
$l_1,l_2,l_3$ meeting at $c$, as shown in Figure~5a. Each 
$\ol\V_{\mu_i}$ is irreducible of dimension 1. In addition, if 
$1\leq i<j\leq 3$ then $\ol\V_{\mu_i}\cap\ol\V_{\mu_j}=\V_{\mu_{i,j}}$ where 
$\mu_{i,j}\in\Q_\Delta^+$ is (any) point having image under (11.4.4) in 
the open region $R_{i,j}$ whose boundary contains $l_i$ and 
$l_j$, as shown in Figure~5a. Each $\V_{\mu_{i,j}}$ is a point.

Let $N(\V)$ be the number of irreducible components of $\V$. If 
$g_X=g_Y=3$ then $N(\V)$
is the number of marked points in Figure~5a. So $N(\V)=9$, as shown in 
Theorem~11.2. If $g_X=2$ and $g_Y=4$ then $N(\V)$ is the number of 
points in Figure~5b that are either marked or in the intersection of a 
solid line and a dashed line. There are $19$ marked points and 
$6$ points of intersection of a solid line and a dashed line in 
Figure~5b. Thus $N(\V)=25$. Note that in this case 
Theorem~11.2 says simply that $N(\V)\geq 19$. 
\end{subsct}

\vspace{0.5cm}

{\smallsc Instituto de Matem\'atica Pura e Aplicada, 
Estrada Dona Castorina 110, 22460-320 Rio de Janeiro RJ, Brazil}

{\smallsl E-mail address: \small\verb?esteves@impa.br?}

\vskip0.4cm

{\smallsc Universidade Federal Fluminense, 
Rua M\'{a}rio Santos Braga, s/n, Valonguinho, 24020-005 Niter\'{o}i RJ, 
Brazil}

{\midrm New address: \smallsc Massachusetts Institute of Technology, 
Department of Mathematics, 77~Massachusetts Avenue, 02139 Cambridge MA, 
U.S.A.}

{\smallsl E-mail address: \small\verb?nivaldo@impa.br?}

\end{document}